\documentclass{article}
\usepackage[utf8]{inputenc}
\usepackage{authblk}

\newcommand{\red}[1]{{\color{red}{#1}}} 

\usepackage{amsmath, mathtools}
\DeclareMathOperator*{\argmin}{argmin}

\usepackage{amssymb,amsbsy}
\usepackage{amsfonts}
\usepackage{amsthm}
\usepackage{bbm}
\usepackage[mathscr]{eucal}
\usepackage{geometry}
\usepackage{authblk}
\usepackage{float}
\usepackage{comment}
\usepackage{enumitem}
\usepackage{soul}
\usepackage{fourier}
\usepackage{multicol}
\usepackage{makecell}

\usepackage[normalem]{ulem}

\usepackage{hyperref}
\hypersetup{
	colorlinks = true,
	urlcolor = {magenta},
	citecolor = {red},
	linkcolor = {blue}
}
\usepackage[capitalize]{cleveref}
\usepackage{doi}
\usepackage[square,numbers]{natbib}

\usepackage{xcolor}
\usepackage{todonotes}

\newcommand{\RNum}[1]{\uppercase\expandafter{\romannumeral #1\relax}}

\DeclareMathAlphabet\mathbfcal{OMS}{cmsy}{b}{n}
\newcommand{\ccT}{\mathbfcal{T}}
\newcommand{\bbt}{\mathbbm{t}}
\newcommand{\pr}{\mathbb{P}}								

\newcommand{\bb}{\begin{equation}}
\newcommand{\ee}{\end{equation}}

\newcommand{\Prob}[1]{\pr\left(#1\right)}					

\newcommand{\CProb}[2]{\pr\left(#1 \Bigm| #2\right)}	

\newcommand{\E}{\mathbb{E}}								

\newcommand{\Exp}[1]{\E\left[#1\right]}					

\newcommand{\CExp}[2]{\E\left[\left.#1\right|#2\right]}	

\newcommand{\plim}{\ensuremath{\stackrel{\pr}{\rightarrow}}}	
\newcommand{\dlim}{\ensuremath{\stackrel{d}{\rightarrow}}}		

\newcommand{\1}{\mathbbm{1}}								
\newcommand{\ind}[1]{\1_{\{#1\}}}	
\newcommand{\indE}[1]{\1_{#1}}					

\newcommand{\NM}[1]{\todo[inline, color=orange]{\color{black}Neeladri: #1 }}
\newcommand{\BC}[1]{\todo[inline, color=green!75!black]{\color{black}Beno\^it: #1 }}
\newcommand{\MU}[1]{\todo[inline, color=yellow]{\color{black}Meltem: #1 }}

\newcommand\numberthis{\addtocounter{equation}{1}\tag{\theequation}}

\allowdisplaybreaks

\newtheorem{theorem}{Theorem}[section]
\newtheorem{lemma}[theorem]{Lemma}

\newtheorem{proposition}[theorem]{Proposition}
\newtheorem{corollary}[theorem]{Corollary}
\newtheorem{conjecture}[theorem]{Conjecture}

\theoremstyle{definition}
\newtheorem{definition}[theorem]{Definition}

\newtheorem{remark}[theorem]{Remark}

\numberwithin{equation}{section}

\newcommand{\cB}{\mathcal B}

\newcommand{\bT}{\mathbf T}

\newcommand{\cC}{\mathcal C}
\newcommand{\cD}{\mathcal D}
\newcommand{\cE}{\mathcal E}

\newcommand{\cP}{\mathcal P}

\newcommand{\cT}{\mathcal T}
\newcommand{\cW}{\mathcal W}

\newcommand{\cY}{\mathcal Y}
\newcommand{\cZ}{\mathcal Z}
\newcommand{\T}{\mathbb T}

\newcommand{\C}{\mathbb C}

\newcommand{\Z}{\mathbb Z}
\newcommand{\N}{\mathbb N}
\newcommand{\R}{\mathbb R}

\newcommand*{\be}{\begin{equation}}
\newcommand*{\ba}{\begin{aligned}}
	\newcommand*{\ea}{\end{aligned}}

\newcommand{\toindis}{\overset d\longrightarrow}
\newcommand{\toinp}{\overset{\mathbb P}{\longrightarrow}}

\newcommand{\eps}{\epsilon}

\newcommand{\floor}[1]{\left\lfloor #1\right\rfloor}
\newcommand{\ceil}[1]{\left\lceil #1 \right\rceil}

\definecolor{clou}{rgb}{0.8,0.25,0.5125}

\newcommand{\invisible}[1]{}
\newcommand{\numitem}[1]{$\mathbf{(#1)}$}

\title{On exponentially height-penalized random trees}
\author[1]{Louigi Addario-Berry}
\affil[1]{McGill University. Email: louigi.addario@mcgill.ca}
\author[2]{Benoît Corsini}
\affil[2]{National University of Singapore. Email: bcorsini@nus.edu.sg}
\author[3]{Neeladri Maitra}
\affil[3]{University of Illinois, Urbana-Champaign. Email: nmaitra@illinois.edu}
\author[4]{Meltem \"Unel}
\affil[4]{Universit\'e Sorbonne Paris-Nord. Email: unel@lipn.univ-paris13.fr}
\date{}

\begin{document}
\maketitle

\begin{abstract}
    Given $n \in \N$ and $\mu \in \R$, a {\em $\mu$-height-biased tree of size $n$} is a random plane tree $\bT_n$ with $n$ vertices with law given by $\pr(\bT_n=t) \propto e^{-\mu h(t)}$, where $t$ ranges over fixed plane trees with $n$ vertices, and $h(t)$ is the height of $t$. 

    Fix a sequence $(\mu_n)_{n \ge 1}$ of real numbers, and for $n \ge 1$ let $\bT_n$ be a $\mu$-height-biased tree of size $n$. In \cite{durhuus2023trees}, the authors describe the asymptotic behaviour of $h(\bT_n)$ when $\mu_n \equiv \mu \in \R$ is fixed.
    In this work, we extend their results to arbitrary sequences of positive parameters depending on $n$.
    Most notably, we show that such a tree behaves like a height-biased Continuum Random Tree (CRT) when $\mu_n$ is of order $1/\sqrt{n}$; that its height is asymptotically $(2\pi^2n/\mu_n)^{1/3}$ when $\mu_n$ is of larger order than $1/\sqrt{n}$ and of smaller order than $n$; and that its height converges to a fixed constant when $\mu_n$ is of order at least $n$, with some random jumps under specific conditions on $\mu_n$.

    We additionally prove various results on second order behaviours, and large deviation principles for the height, for different regimes of $\mu_n$.
    Finally, we describe new statistics of these trees, covering their widths, their root degrees, and the local structure around their roots.
\end{abstract}

\tableofcontents

\section{Introduction}\label{sec:intro_new}

Trees have been at the heart of active research for many decades in combinatorics, probability, and theoretical computer science.  Driven in part by the mathematical analysis of algorithms, asymptotic methods from combinatorics have been employed for the enumeration of different classes of trees. These results, being important in their own right, have also been essential in defining various probability measures on the different tree ensembles, hence for the probabilistic study of random tree models. 
A large class of random tree models is given by the so-called \emph{simply generated trees}, introduced by Meir and Moon in \cite{meir1978altitude}. These models offer a \emph{local} picture by construction, in the sense that the weight of a tree is proportional to the weight factor associated to each vertex depending on the degree of the latter. A well studied sub-class of such trees is given by the genealogical trees of the Bienaym\'e processes (also known as Galton--Watson processes). The local and scaling limits coming from different variants of this model have been studied extensively and correspond to large universality classes \cite{abraham2015introduction,gall2011scaling}. The uniform random plane tree model, which undergirds the models studied in this work, turns out to be an instance of the Bienaymé model with the offspring distribution being geometric of parameter $\frac{1}{2}$, and conditioned on the number of vertices.

In \cite{durhuus2023trees},  the authors studied a random tree model with \emph{non-local} weight factors where the height of the tree -- a global statistic -- is the argument of the weight function: for $n \in \N$ and $\mu \in \R$, the random plane tree $\bT_n$ with $n$ vertices with law given by $\pr(\bT_n=T_n) \propto e^{-\mu h(T_n)}$, where $T_n$ ranges over fixed plane trees with $n$ vertices, and $h(T_n)$ is the height of $T_n$. This model exhibits a phase transition at $\mu=0$ where the local limit passes from a single-spine phase to a multi-spine one. The multi-spine tree, called the Poisson tree, obtained when $\mu>0$, does not fall into a well-known universality class where we find the local limits of large Bienaymé trees. We note that the Poisson tree can also be obtained by conditioning the uniform plane trees on their width growing as a prescribed sequence with a certain growth-rate as a function of the size \cite{abraham2020very}. Hence, along with the interesting properties of this limiting object, biasing by height turns out to be a rather "natural" way to obtain trees with rather unusual properties.

In this work, we extend the results of \cite{durhuus2023trees} to the case where $\mu=\mu_n$ is a positive sequence that can depend on $n$. Before presenting a wide range of results on different properties of these \emph{$\mu$-height-biased trees}, let us {first set some basic terminology regarding trees and} introduce the model formally. 

\paragraph{Notations and conventions for trees.} All trees arising in this paper are finite. We consider \emph{rooted plane trees} which we also simply refer to as \emph{trees} throughout the paper. {Thus, the trees that we consider {have} a distinguished vertex which we call the \emph{root} of the tree, and for any vertex, there is a left-to-right ordering of its children.}
We write $\bbt_n$ for the set of rooted plane trees with $n$ nodes.
It is well-known (see e.g., \cite[Chapter 1]{drmota2009random}) that $|\bbt_n|=C_{n-1}$ where $C_n$ is the $n$-th Catalan number:
\begin{align}\label{eq:catalan}
    C_n=\frac{1}{n+1}\binom{2n}{n}\,.
\end{align}

{For a tree $T$, we slightly abuse notation and also use $T$ to denote its vertex set. For any tree $T$, the  \emph{height} or \emph{depth} of a vertex $v \in T$ is its distance from the root, i.e., the number of edges on the unique path connecting $v$ with the root; we denote this quantity by $h(v)$. Note that by definition the root has height zero.} We also define its \emph{height} $h(T)$ of the tree $T$ to be the largest depth of any of its nodes, {i.e.,
\begin{align*}
    h(T)=\max\Big\{h(v):v\in T\Big\}\,.
\end{align*}
}

{For a tree $T$, 
and non-negative integer $k$, the {\em $k$'th generation} of $T$ is the set of vertices at distance exactly $k$ from the root of $T$. We write $\cZ_k=\cZ_k(T)$ for the size of the $k$'th generation, and $w(T)$ for the {\em width}  of $T$, which is defined to be the size of the largest generation, i.e.,
\begin{align}\label{eq:def width}
    w(T)=\max\Big\{\cZ_k:k\geq 0\Big\}\,.
\end{align}
}

\paragraph{Model definition.}

The model of random trees introduced in \cite{durhuus2023trees} is defined as follows.
For any $n\in\N$ and $\mu\in\R$, let $\mathbf{T}_n^{\mu}$ {be the random element of $\bbt_n$} with distribution given by 
\begin{align}\label{eq:T dist}
    \Prob{\mathbf{T}_n^{\mu}=T}=\frac{e^{-{\mu} h(T)}}{Z_n^{\mu}}\,,
\end{align}
where
\begin{align}\label{eq:part fun}
    Z_n^{\mu}=\sum_{T\in\bbt_n}e^{-{\mu} h(T)}
\end{align}
is a normalizing constant, referred to as the \emph{partition function}. We refer to $\bT^\mu_n$ as a \emph{$\mu$-height-biased tree of size $n$}. {Throughout the paper, by $\varnothing$ we denote the root of the random tree $\bT_n$. When $\mu\equiv 0$ in \eqref{eq:T dist}, we obtain the uniform measure, and we denote the corresponding uniform random tree of size $n$ by $\T_n$ for the rest of this paper.} 

In~\cite{durhuus2023trees}, the authors studied local properties of this model when $\mu$ is a fixed real, not depending on $n$, and observed three different {asymptotic} behaviours {as $n \to \infty$}, depending on whether $\mu$ was equal to $0$, positive, or negative.
{The} case $\mu=0$ is the usual uniform plane tree whose asymptotic local limit is Kesten's tree associated to the geometric offspring distribution of parameter $1/2$, see~\cite{kesten1986subdiffusive}, which we refer to as simply Kesten's tree (KT) for the rest of the paper. 
In the present work, our goal is to {extend these results to} {consider} the case where $\mu=\mu_n{>0}$ can depend on $n$. Due to the large range of new behaviours arising from {$\mu_n$ to depend on $n$}, {in this work we} focus on the case that $\mu_n>0$, and leave the case $\mu_n<0$ for future work.) 
For the rest of the paper, we usually remove the subscript $n$ from $\mu$ and the superscript $\mu$ from $\bT_n^\mu$ and $Z_n^\mu$, {simply writing $\bT_n$ and $Z_n$,} except when we need to make the distinction between two distributions with different $\mu$-parameters. {After introducing some standard notation, in the sections that follow, we state our main results.}

\paragraph{Notation.}
{We write $x_n\asymp y_n$, $x_n\ll y_n$, $x_n\gg y_n$, $x_n=O(y_n)$, and $x_n=\Omega( y_n)$, respectively, to indicate that $x_n/y_n$ is bounded away from $0$ and $\infty$, that $x_n/y_n$ converges to $0$, that $x_n/y_n$ diverges to $\infty$, that $x_n/y_n$ is bounded, and that $x_n/y_n$ is bounded away from $0$. Sometimes it is convenient for us to write $x_n=o(y_n)$ for $x_n\ll y_n$. We also write $x_n\sim y_n$ to say that $x_n/y_n$ converges to $1$. For $x_n$ and $y_n$ possibly random sequences, we write $x_n\asymp_{\mathbb{P}}y_n$ to indicate that the ratio $x_n/y_n$ is \textit{bounded in probability} from above and below, i.e., for any $\eps>0$, one can obtain a $K=K(\eps)$ so that $\Prob{K^{-1}\leq|x_n/y_n|\leq K}\geq 1-\eps$. For any function $f:S\to \R$ for some set $S$, we denote $\|f\|:=\sup_{s\in S}|f(s)|$ the supremum value of $f$, possibly infinite. By $\mathrm{Geo}(p)$ we denote the law of a discrete random variable having the geometric distribution with success parameter $p$.}

\subsection{Global properties in the Brownian regime}
In this section, we state results on \emph{global} properties of the tree $\mathbf{T}_n$ when $\mu=O( 1/\sqrt{n})$.
Note that in this regime, we are allowing $\mu<0$. This is the only section of the paper where we allow {$\mu$ to take negative values}.
We term this the \emph{Brownian regime} due to the strong relation between our model in this regime and Brownian motions, in particular through the \emph{Continuum Random Tree (CRT)}.

Our first result, on the scaling limit in the Brownian regime, establishes the convergence of the rescaled tree, seen as a metric space, toward a height-biased CRT. {For more precise definitions about the objects appearing below, and the topology appearing in the statement of Theorem \ref{thm:scaling limit}, we refer the reader to Section \ref{sec:prelims}.}

{Let $\mathbf{e}=(\mathbf{e}(t):t \in [0,1])$ be the standard Brownian excursion on $[0,1]$. Consider the random function $\mathbf{e}_\alpha=(\mathbf{e}_\alpha(t):t \in [0,1])$ whose law $\rho_\alpha$ has Radon-Nikodým derivative with respect to the law $\rho_0$ of $\mathbf{e}$ given by
\begin{align*}
    \frac{d\rho_\alpha}{d \rho_0}(f):=\frac{e^{-\alpha \|f\|}}{\int e^{-\alpha \| g \| }d\rho_0 (g)}\, ,
\end{align*}
and where we write $\|g\|:=\max_{t \in [0,1]}|g(t)|$ for any $g \in C([0,1],\R)$, the set of all continuous functions from $[0,1]$ to $\R$. Define $(\ccT_\alpha,d_\alpha,\sigma_\alpha)$ be the measured random continuum tree encoded by the function $\mathbf{e}_\alpha$.
}

\begin{theorem}[Scaling limit in the Brownian regime]\label{thm:scaling limit}
    Let $\mu=\mu_n$ satisfy $\lim_{n\rightarrow\infty}\mu\sqrt{n}=\alpha\in\R$.
    Let $\mathbf{T}_n=\bT^\mu$ be a $\mu$-height-biased tree of size $n$, let $d_n$ denote the graph distance on $\mathbf{T}_n$ with each edge having length 1, and let $\sigma_n$ denote the uniform measure on the  vertices of ~$\bT_n$. Then
    \begin{align*}
        \left(\mathbf{T}_n,\frac{d_n}{\sqrt{2n}}, {\sigma_n}\right) \toindis \big(\ccT_{\alpha},d_{\alpha}\,, \sigma_{\alpha}\big),
    \end{align*}
    where the convergence in distribution above is with respect to the Gromov-Hausdorff-Prokhorov (GHP) topology.
\end{theorem}

The proof of this theorem can be found in Section~\ref{sec:global_result_proofs}.
Thus, the global geometry of the random tree $\mathbf{T}_n$ is completely described by appropriately tilting
{the measured continuum random tree encoded by $\mathbf{e}$}. We next complement the above scaling limit with a characterization of the \emph{width} of the tree, as defined in~\eqref{eq:def width}, {whose proof can again be found in} ~\cref{sec:global_result_proofs}.

{For the next result, let us introduce the random function $B_\alpha=(B_\alpha(t):t \in [0,1]) \in C([0,1],\R)$, whose law $\xi_\alpha$ has Radon-Nikodým derivative with respect to {the law $\rho_0$ of}  $\mathbf{e}$ given by
\begin{align*}
\frac{d \xi_{\alpha}}{d\rho_0}(f)=\frac{e^{-\alpha \Bar{f}}}{\int e^{-\alpha \Bar{g}} d\rho_0(g)}\, , \numberthis \label{eq:law_balpha_xialpha}
        \end{align*}
where for $f\in C([0,1],\R)$ we use the notation 
\begin{align*}
    \Bar{f}:=\frac{1}{2}\int_{0}^1 \frac{dt}{f(t)}.
\end{align*}
}

\begin{remark}[On the process $B_\alpha$ being well-defined]
    In general, for a continuous function $f$, the integral $\int_0^1\frac{ds}{f(s)}$ need not  be well-defined, since there can be a positive measure subset of $[0,1]$ where $f$ vanishes. However, it is classical that this does not happen a.s. for the Brownian excursion $\mathbf{e}$ (see, e.g., \cite{chassaing2000height} and the references therein), so that $\int_0^1\frac{ds}{\mathbf{e}(s)}$ is well defined as a random variable. In particular, the law $\rho_0$ of $\mathbf{e}$ is supported by the set of functions $f$ for which the integral $\int_0^1\frac{ds}{f(s)}$ makes sense. Thus, the change of measure that defines $\xi_\alpha$ makes sense and therefore the process $B_\alpha$ is well-defined.
\end{remark}

\begin{theorem}[Width in the Brownian regime]\label{thm:width_scaling_Brownian}
    Let $\mu=\mu_n$ satisfy $\lim_{n\rightarrow\infty}\mu\sqrt{n}=\alpha\in\R$.
    Let $\mathbf{T}_n$ be a $\mu$-height-biased tree of size $n$.
    Then the sequence $(w(\mathbf{T}_n))_{n \geq 1}$ satisfies
        \begin{align*}
            \frac{w(\mathbf{T}_n)}{\sqrt{n}} \toindis \|B_\alpha\|.
        \end{align*}
\end{theorem}

\subsection{Global properties in the non-Brownian regime}
The previous section identifies the global behavior of the tree $\mathbf{T}_n$ in the Brownian regime. We next complement the previous results by providing results on global properties of the tree $\mathbf{T}_n$ in the \emph{non-Brownian regime}, when $\mu\gg 1/\sqrt{n}$.
To state the results concisely, it is useful to further subdivide the non-Brownian regime into two parts.
\begin{itemize}
    \item When $1/\sqrt{n}\ll \mu \ll n$, we are in the \emph{intermediate regime}.
    \item When $\mu= \Omega(n)$, we are in the \emph{extreme regime}.
\end{itemize}

We start by stating the results for the asymptotic behaviour of the height of $\mathbf{T}_n$ in the intermediate regime.

\begin{theorem}[Asymptotic height in the intermediate regime]\label{thm:inter lln}
    Let $\mu=\mu_n$ be such that $1/\sqrt{n}\ll\mu\ll n$.
    Let $\mathbf{T}_n$ be a $\mu$-height-biased tree of size $n$.
    Then
    \begin{align*}
        \frac{h(\mathbf{T}_n)}{\left(\frac{2\pi^2n}{\mu}\right)^{1/3}}\toinp 1\,.
    \end{align*}
\end{theorem}
The proof of Theorem~\ref{thm:inter lln} can be found in Section~\ref{sec:lln}.
Note that $(2\pi^2n/\mu)^{1/3}$ ranges from being of order $n^{1/2}$ when $\mu \asymp n^{-1/2}$, to being of order $1$ when $\mu \asymp n$.
We now extend this result to include the fluctuations around the height.
As we present below, the fluctuations actually have a wide range of behaviour: from a standard central limit theorem when $\mu\ll n^{1/4}$, to a discrete central limit theorem when $\mu\asymp n^{1/4}$, and finally to \emph{Bernoulli fluctuations} when $\mu\gg n^{1/4}$, the latter behavior also extending into the extreme regime.

\begin{theorem}[Gaussian fluctuations for the height]\label{thm:inter clt}
    Let $\mu=\mu_n$ be such that $1/\sqrt{n}\ll\mu\ll n^{1/4}$.
    Let $\mathbf{T}_n$ be a $\mu$-height-biased tree of size $n$.
    Then
    \begin{align*}
        \frac{h(\mathbf{T}_n)-\left(\frac{2\pi^2n}{\mu}\right)^{1/3}}{\sqrt{\frac{1}{3}\left(\frac{2\pi^2n}{\mu^4}\right)^{1/3}}}\toindis\mathcal{N}(0,1)\,,
    \end{align*}
    where $\mathcal{N}(0,1)$ is a standard normal random variable.
\end{theorem}

The proof of Theorem~\ref{thm:inter clt} can be found in Section~\ref{sec:clt}.
Observe that when $\mu$ is of order $n^{1/4}$, the denominator, which corresponds to the standard deviation of the height, becomes of constant order.
In that regime, the previous result almost directly holds, except that the convergence cannot be towards a normal random variable, since $h(\mathbf{T}_n)$ is an integer, but towards the integer equivalent of a normal random variable. To state this result, for any $t \in (2,\infty)$ and $x>0$ define the function
\begin{align}
    \lambda_x(t):=xt+\log\left(1+\tan^2 \left(\frac{\pi}{t}\right) \right)\label{eq:func_lambdax_minvalue}
\end{align}
and let $t_x:= \argmin\lambda_{x}(t)$ be its minimizing argument.
It is straightforward to check that $\lambda_x$ is convex and thus admits a unique minimizer. We refer the reader to Proposition~\ref{prop:asympt lx} in Appendix~\ref{sec:on part fun} for further properties of $\lambda_x$ and $t_x$.

\begin{theorem}[Discrete central limit theorem for the height]\label{thm:discrete clt}
    Let $\mu=\mu_n$ satisfy $\lim_{n\rightarrow\infty}\mu/n^{1/4}=\gamma>0$. Assume that $t_{\mu/n}-\floor{t_{\mu/n}}$ admits a limit $\delta\in[0,1]$.
    Then the asymptotic of the height $h(\mathbf{T}_n)$ is given by
    \begin{align*}
        h(\mathbf{T}_n)-\floor{t_{\mu/n}}+2\toindis X_{\gamma,\delta},
    \end{align*}
    where $X_{\gamma,\delta}$ is an integer-valued random variable distributed as
    \begin{align*}
        \Prob{X_{\gamma,\delta}=k\big.}=\frac{1}{C_{\gamma,\delta}}\exp\left(-3\left(\frac{\gamma^2}{4\pi}\right)^{2/3}\big(k-\delta\big)^2\right)\,,\;\;k\in \Z
    \end{align*}
    with $C_{\gamma,\delta}$ being a normalizing constant corresponding to the sum over $k\in\Z$ of the exponential term.
\end{theorem}

\begin{remark}
    We remark that if $\mu/n^{1/4}$ is of constant order but does not admit a limit, or if $t_{\mu/n}-\floor{t_{\mu/n}}$ does not admit a limit, the asymptotic distribution of $h(\mathbf{T}_n)-\floor{t_{\mu/n}}+2$ is still close to that of $X_{\mu/n^{1/4},t_{\mu/n}-\floor{t_{\mu/n}}}$, although we do not obtain a distributional limit as stated in the last result.
\end{remark}

The proof of Theorem~\ref{thm:discrete clt} can be found in Section~\ref{sec:clt}.
We observe that, while the distribution of $X_{\gamma,\delta}$ strongly resembles that of a normal distribution, since it replaces the integral by a sum, it is different from the distribution of the integer part of a normal random variable.

 After addressing the first order behavior of the height of $\mathbf{T}_n$ in Theorem~\ref{thm:inter lln} and its fluctuations in Theorem~\ref{thm:inter clt} and Theorem~\ref{thm:discrete clt}, the next result covers the probabilities of large deviations from the mean.
\begin{theorem}[Large deviations of the height]\label{thm:ldp}
    Let $\mu=\mu_n$ be such that $1/\sqrt{n}\ll\mu\ll n$.
    Let $\mathbf{T}_n$ be a $\mu$-height-biased tree of size $n$.
    Then the sequence of random variables $(h(\bT_n)/(2\pi^2n/\mu)^{1/3})_{n \geq 1}$ satisfies a large deviation principle at speed $(\mu^2 n)^{1/3}$ and with good rate function given by
    \begin{align*}
        \Lambda^*(x)=\left\{\begin{array}{ll}
            \left(\frac{\pi}{2}\right)^{2/3}\frac{(x-1)^2(2x+1)}{x^2} & \textrm{if $x>0$} \\
            \infty & \textrm{otherwise.}
        \end{array}\right.
    \end{align*}
\end{theorem}

The proof of Theorem~\ref{thm:ldp} can be found in Section~\ref{sec:ldp}.
With this last result on the height in the intermediate regime, we next move on to other functionals of $\mathbf{T}_n$. Our first result in this direction concerns the width of $\mathbf{T}_n$, as defined in~\eqref{eq:def width}. This result covers \emph{almost} the entire intermediate regime, and extends beyond the intermediate regime well into the extreme regime.

\begin{theorem}[Width in the non-Brownian regime]\label{thm:width}
Let $\mu=\mu_n$ be such that $\mu\gg (\log n)^{3/2}/\sqrt{n}$.
Let $\mathbf{T}_n$ be a $\mu$-height-biased tree of size $n$.
Then for any $\eps>0$, we can find $K=K(\eps)>0$ large enough, such that
        \begin{align*}
            \liminf_{n \to \infty} \Prob{K^{-1}\leq\frac{w(\mathbf{T}_n)}{n\wedge\left((\mu n^2)^{1/3}\right)}\leq K}\geq 1-\eps\,.
        \end{align*}
\end{theorem}
The proof of Theorem~\ref{thm:width} can be found in Section~\ref{sec:inter_width_proof}. Note that in the current version we do not characterize the width behavior for the entire intermediate regime: namely, our result has a gap for the regime $1/\sqrt{n}\ll \mu= O\left((\log n)^{3/2}/\sqrt{n}\right)$. In the last theorem, we only provide a tightness result on the scaled width. We conjecture a \emph{law of large numbers} for the width in the entire non-Brownian regime:
\begin{conjecture}[LLN for the width in the non-Brownian regime]\label{conj:LLN_width}
    Let $\mu=\mu_n$ be such that $\mu\gg1/\sqrt{n}$.
    Let $\mathbf{T}_n$ be a $\mu$-height-biased tree of size $n$.
    There exists a constant $W>0$ such that as $n \to \infty$,
    \begin{align*}
        \frac{w(\mathbf{T}_n)}{n\wedge\left((\mu n^2)^{1/3}\right)} \toinp W.
    \end{align*}
\end{conjecture}

To prove Theorem \ref{thm:width}, we employ a probabilistic technique involving a multiscale decomposition of the associated \emph{random excursion} of the tree $\mathbf{T}_n$. More discussion on this technique and the proof of this result can be found in Section \ref{sec:inter_width_proof}. 

Finally, we state our  main result on the height of $\mathbf{T}_n$ in the extreme regime, also covering fluctuations when $\mu \gg n^{1/4}$. {Recall $\lambda_x$ and $t_x$ from \eqref{eq:func_lambdax_minvalue}. For any $x>0$ define $m_x=\max\{\floor{t_x},3\}$ and
\begin{align}
    \delta_n=\mu+n\log\left(\frac{1+\tan^2\left(\frac{\pi}{\lceil t_{\mu/n}\rceil}\right)}{1+\tan^2\left(\frac{\pi}{m_{\mu/n}}\right)}\right)\,. \label{eq:def_m_n_delta_n}
\end{align}
}

\begin{theorem}[Height fluctuations beyond the Gaussian sub-regime]\label{thm:high lln}
    Let $\mu=\mu_n$ be such that $\mu\gg n^{1/4}$ and such that $\mu/n$ converges to some $c\in[0,+\infty]$.
    Further, assume that $\delta_n$ admits a limit $\delta\in [-\infty,\infty]$.
    Let ~$\mathbf{T}_n$ be a $\mu$-height-biased tree of size $n$.
    Then, $h(\mathbf{T}_n)-m_{\mu/n}+2$ converges in distribution to a Bernoulli random variable with parameter
    \begin{align*}
        p_{c,\delta}=\left\{\begin{array}{ll}
            \displaystyle\frac{m_c\tan^2\left(\frac{\pi}{m_c+1}\right)e^{-\delta}}{(m_c+1)\tan^2\left(\frac{\pi}{m_c}\right)+m_c\tan^2\left(\frac{\pi}{m_c+1}\right)e^{-\delta}} & \textrm{if $c>0$} \\
            \displaystyle\frac{e^{-\delta}}{1+e^{-\delta}} & \textrm{otherwise}\,.
        \end{array}\right.
    \end{align*}
\end{theorem}
Observe that, since $m_c$ diverges to $\infty$ when $c$ converges to $0$, the two formulas in the definition of $p_{c,\delta}$ are coherent with each other, although need to be split since the formula is not well-defined for $m_c=\infty$.
The proof of Theorem~\ref{thm:high lln} can be found in Section~\ref{sec:bern}.
This theorem states that $h(\mathbf{T}_n)$ is asymptotically deterministic, except in specific regimes for $\mu$ where it transitions from one value to another.
This is to be expected, as {the formula for} the variance of $h(\mathbf{T}_n)$ {as found in Theorem~\ref{thm:inter clt}} converges to $0$ whenever $\mu\gg n^{1/4}$, and since $h(\mathbf{T}_n)$ takes integer values and thus cannot change from one value to another without transition. 

Furthermore, we observe that $p_{c,\delta}$ is equal to $0$ when $\delta=\infty$ and to $1$ when $\delta=-\infty$.
This can be understood by fixing the limiting value $m_{\mu/n}=\Tilde{m}$ and letting $\mu$ approach the threshold
\begin{align*}
        n\log\left(\frac{1+\tan^2\left(\frac{\pi}{\Tilde{m}+1}\right)}{1+\tan^2\left(\frac{\pi}{\Tilde{m}}\right)}\right)\,.
\end{align*}
Indeed, when $\delta=-\infty$, meaning that $\mu$ is far below the threshold, then the height is less biased and thus higher, leading to $h(T_n)=\Tilde{m}-1$.
Then, as $\mu$ approaches the threshold and $\delta$ becomes of constant order, it becomes possible for $h(T_n)$ to decrease by one level, from $\Tilde{m}-1$ to $\Tilde{m}-2$.
Finally, when $\mu$ has long passed the threshold and $\delta=+\infty$, the height-bias is too large and we have reduced the height by $1$, leading to $h(T_n)=\Tilde{m}-2$.

\begin{remark}[Dichotomy of the concentration of the height about its mean] We observe the following dichotomy on the concentration behavior of the height about its mean. When $\mu=O( 1/\sqrt{n})$, the height has random fluctuations at its highest order (Theorem \ref{thm:scaling limit}), but also when $\mu =\Omega( n)$ (Theorem \ref{thm:high lln}); while in the middle regime of $1/\sqrt{n}\ll \mu \ll n$, the height is more concentrated and satisfies a LLN (Theorem \ref{thm:inter lln}), where precise large-deviation estimates can be given (Theorem \ref{thm:ldp}).
\end{remark}

Using that the height of a tree with more than $2$ nodes cannot be less than $1$, it is interesting to extract from the previous result the regime in which $h(T_n)$ is minimal, and so the tree is simply the star $S_n$ with $n-1$ leaves, as stated in the following corollary.

\begin{corollary}[Observing a star]\label{cor:star}
    Let $\mu=\mu_n$ be such that $\mu-n\log2$ admits a limit $\delta\in\R\cup\{\pm\infty\}$.
    Then as $n \to \infty$
    \begin{align*}
        \Prob{h(\mathbf{T}_n)=1\big.}=\Prob{T_n=S_n\big.}\longrightarrow\frac{4}{4+e^{-\delta}}\,.
    \end{align*}
\end{corollary}
\begin{remark}
    Observe that this result does not require $\mu\gg n^{1/4}$, or even $\mu\gg1/\sqrt{n}$. Indeed, when $\mu$ is sublinear, i.e., $\mu\ll n$, then $\delta=-\infty$, and in this case the corollary simply tells us that the probability of obtaining the star graph is {asymptotically} $0$. 
\end{remark}

Let us quickly give the proof of the last corollary.

\begin{proof}[Proof of Corollary \ref{cor:star}]
    This is a direct consequence of Theorem~\ref{thm:high lln} by observing that
    \begin{align*}
        \log\left(\frac{1+\tan^2\left(\frac{\pi}{3}\right)}{1+\tan^2\left(\frac{\pi}{4}\right)}\right)=\log\left(\frac{1+3}{1+1}\right)=\log2\,.
    \end{align*}
    Indeed, if $m_{\mu/n}\neq3$ and $\lceil t_{\mu/n}\rceil\neq4$, then this means that $\mu$ is lower than the possible threshold considered here, leading to $\delta=-\infty$.
    Using the different results regarding the height of $h(T_n)$ (Theorem~\ref{thm:scaling limit},~\ref{thm:inter lln}, or~\ref{thm:high lln}), we see that $T_n=S_n$ with probability $0$ here, thus satisfying the claimed result.
    On the other hand, if $m_{\mu/n}=3$ and $\lceil t_{\mu/n}\rceil=4$, then we can directly apply Theorem~\ref{thm:high lln}.
    Finally, if $m_{\mu/n}=\lceil t_{\mu/n}\rceil=3$, then this means that $t_{\mu/n}\leq3$, in which case
    \begin{align*}
        \frac{\mu}{n}\geq\frac{2\pi}{9}\tan\frac{\pi}{3}=\frac{2\pi}{3^{3/2}}>1>\log2\,.
    \end{align*}
    It follows that $\delta=+\infty$ and then Theorem~\ref{thm:high lln} confirms that $T_n=S_n$ with high probability.
\end{proof}
Let us note that in the very special case when $\mu \gg n$, we have a positive answer to Conjecture \ref{conj:LLN_width} due to Corollary \ref{cor:star}. As $\delta =\infty$, we see a star with high probability, and thus $w(\mathbf{T}_n)/n \to 1$ in probability. This finishes all the statements about global functionals in the non-Brownian regime. In the next section, we state results on \emph{local} properties of $\mathbf{T}_n$.

\subsection{Results on local properties}

In this section, we state a dichotomy between two qualitatively different local behaviors around the root of the tree $\mathbf{T}_n$, depending on whether or not $\mu$ diverges to infinity. We name these regimes as,
\begin{itemize}
    \item \textbf{Local regime}, when $\mu=O( 1)$.
    \item \textbf{Condensation regime}, where $\mu\gg 1$.
\end{itemize}
For the next result, the relevant definitions for the topology can be found in Section~\ref{sec:prelims}.

\begin{theorem}[Local limit in the local regime]\label{thm:loclim}
Let $\mu=\mu_n$ be a sequence admitting a limit $\alpha \in [0,\infty)$.
Let $\mathbf{T}_n$ be a $\mu$-height-biased tree of size $n$. If $\alpha>0$, then the random rooted tree $(\mathbf{T}_{n},\varnothing)$ converges locally to the infinite \emph{Poisson tree} of parameter $\alpha$ as defined in \cite{abraham2020very}. Otherwise, if $\alpha=0$, equivalently $\mu=o(1)$, then the random rooted tree $(\mathbf{T}_{n},\varnothing)$ converges locally to Kesten's tree as defined in \cite{kesten1986subdiffusive}.
\end{theorem}
The proof of Theorem~\ref{thm:loclim} can be found in Section~\ref{sec:loc_results}.
We remark that the statement on the local limit when $\mu=\mu_n\equiv \alpha>0$ is a constant sequence was already identified in \cite[Theorem 4.3]{durhuus2023trees}. The above result additionally extends to the case where $\mu\ll 1$.

For any tree $T$, and any vertex $v \in T$, we denote by $\deg_T(v)$ the \emph{number of children} of the vertex $v$ in $T$, i.e., the number of neighbors of $v$ {in $T$} not lying on the ancestral path of $v$ to the root. When $\mu\gg1$, we observe a condensation phenomenon at the root, i.e., its degree diverges. 

\begin{theorem}[Condensation at the root]\label{thm:root deg}
    Let $\mu=\mu_n$ be such that $1\ll\mu\ll n^{1/4}$.
    Let $\mathbf{T}_n$ be a tree distributed according to~\eqref{eq:T dist}.
    Then the asymptotic behavior of the degree of the root $\varnothing$ in $\mathbf{T}_n$ is given by
    \begin{align*}
        \frac{\deg_{\mathbf{T}_n}(\varnothing)-2\mu}{\sqrt{6\mu}}\toindis\mathcal{N}(0,1)\,.
    \end{align*}
\end{theorem}
The proof of this theorem can be found in Section~\ref{sec:root deg}.
This result shows that, in the non-local regime and up to $\mu \ll n^{1/4}$, we observe a condensation phenomenon as the root degree starts diverging at scale $\mu$, as opposed to the local regime, where it is almost surely finite. We conjecture that the root degree satisfies a law of large numbers when $\mu=\Omega(n^{1/4})$.
\begin{conjecture}\label{conj:root degree}
 Assume $n\gg \mu =\Omega( n^{1/4})$. Then as $n \to \infty$, there is a constant $C>0$ such that
 \begin{align*}
     \frac{\deg_{\mathbf{T}_n}(\varnothing)}{\mu}\toinp C.
 \end{align*}
\end{conjecture}
The second order fluctuations of the root degree in the regime $n\gg \mu=\Omega( n^{1/4})$ appear to be more subtle, and we do not have a conjecture at this point.

\begin{remark}[Root degree when $\mu =\Omega(n)$]
    Observe that Corollary~\ref{cor:star} indirectly states that, whenever $\mu-n \log 2 \to \delta \in \R \cup \{\pm\infty\}$, then $\Prob{\deg_{\mathbf{T}_n}(\varnothing)=n-1}=4/(4+e^{-\delta})$.
    This suggests that there exists a regime for $\mu$ where the limit $(n-1)-\deg_{\mathbf{T}_n}(\varnothing)$ is non-trivial.
    Moreover, and in order to connect the regime $\mu=\Omega(n)$ with that of $\mu\ll n$ from Conjecture~\ref{conj:root degree}, we also expect the existence of a regime for $\mu$ where $(n-1)-\deg_{\mathbf{T}_n}(\varnothing)$ diverges but $\deg_{\mathbf{T}_n}(\varnothing)/n$ is non-trivial.
    At the moment, it is unclear to us where the exact change of behavior occurs and how they impact the two types of limit (subtracting $n$ and dividing by $n$).
\end{remark}

\begin{remark}[Comparison with the local results of \cite{abraham2020very}]
As a result of Theorems \ref{thm:loclim} and \ref{thm:root deg}, we observe the following qualitative picture in the transition of the local behavior around the root. As $\mu$ increases, locally around the root we begin with observing Kesten's tree, until it reaches a threshold when $\mu \to \alpha \in (0,\infty)$, when the behavior changes to that of the Poisson tree; beyond this threshold, we observe condensation at the root, in that its degree diverges. A similar picture was also observed in a different context, in the paper \cite{abraham2020very}. There, the authors consider Bienaymé trees conditoned on explicit generation sizes, thus penalizing the height \emph{indirectly}, unlike the direct exponential penalization of the measure \eqref{eq:T dist} that we consider.
This suggests that this qualitative transition landscape is in a sense \emph{universal}, and may emerge in other contexts.

\end{remark}
\begin{table}[H]
    \renewcommand{\arraystretch}{1.25}
    \centering

    \begin{tabular}{|c|c|c|c|c|c|c|}
        \hline
        \textbf{Regime} & $\mu=O( 1/\sqrt{n})$ & \multicolumn{4}{|c|}{$1/\sqrt{n}\ll\mu\ll n$} & $\mu=\Omega(n)$ \\
        \hline
        \textbf{Sub-regime} & $\mu=O( 1/\sqrt{n})$ & $1/\sqrt{n}\ll\mu\ll1$ & $\mu\asymp 1$ & $1\ll\mu=O(n^{1/4})$ & $n^{1/4}\ll\mu\ll n$ & $\mu=\Omega(n)$ \\
        \hline\hline
        \textbf{Height} & $\asymp_{\mathbb{P}} \sqrt{n}$ & \multicolumn{4}{|c|}{$\sim(2\pi^2n/\mu)^{1/3}$} & $\asymp_{\mathbb{P}} 1$ \\
        \hline
        \textbf{LDP speed for height} & NA
        & \multicolumn{4}{|c|}{$\asymp  (\mu^2n)^{1/3}$} & NA \\
        \hline
        \textbf{Width} & $\asymp_{\mathbb{P}} \sqrt{n}$& \multicolumn{4}{|c|}{$\asymp  (\mu n^2)^{1/3}$ (LLN conjectured)} & $\asymp  n$ \\
        \hline
        \textbf{Local behavior} & \multicolumn{2}{|c|}{Kesten's tree} & Poisson tree &  \multicolumn{3}{|c|}{Condensation} \\
        \hline
        \textbf{Root degree} & \multicolumn{2}{|c|}{$\asymp_{\mathbb{P}}  1$} & \multicolumn{3}{|c|}{$\asymp \mu$ (conjectured for $\mu=\Omega(n^{1/4})$)} & $\asymp  n$ \\
        \hline
        \textbf{Fluctuation of height} & $\asymp_{\mathbb{P}} \sqrt{n}$ & \multicolumn{3}{|c|}{$\asymp_{\mathbb{P}} (n/\mu^4)^{1/6}$} & \multicolumn{2}{|c|}{$O(1)$} \\
        \hline
    \end{tabular}
    \caption{
        A summary of our different results together with Conjectures \ref{conj:root degree} and \ref{conj:LLN_width} depending on the regime of $\mu$. 
    }
    \label{tab:results}
\end{table}

\subsection{Related work and discussion}\label{sec:lit}
As has been mentioned before, the Poisson tree which is the local limit of the model in hand for $\mu > 0$ constant has first been obtained by conditioning the geometric Bienaymé trees on their {generation sizes} in \cite{abraham2015introduction}. The decomposition of the local limit around a fixed finite tree found in \cite[Theorem 4.4]{durhuus2023trees} shows that, conditioned on the degree of the root being say $M$, one observes at least one and possibly more infinite subtrees emanating from those $M$ vertices at height 1, together with finite subtrees grafted on the rest. Each one of these infinite subtrees is a height-biased tree on its own and they have a \emph{joint} law. Tracing these infinite branches back to the stage where $n$ is large (and  $n \to \infty$ limit has not yet been taken), we find the large branches of size ${\asymp} n^{1/3}$. We expect a similar picture when $1 \ll \mu \ll n^{1/4}$: the degree of the root having mean $2\mu$, about half of which carrying large branches of size ${\asymp} ( n/\mu)^{1/3}$, each one of them being a height-biased tree with some positive \emph{constant} bias. Although we have {a number of heuristics that support this expectation}, we have not rigorously investigated it in the present work.

In \cite{bouchot2024scaling}, a very similar problem in case of one-dimensional simple random walk has been treated where the walk is penalized by its range $R$ i.e. the weight of a walk of size $n$ is proportional to $\exp(-\mu R)$. One again observes the stretched exponential in the asymptotics of the partition function. The limit laws for the range as $n\to \infty$ resemble closely to our results for the height and the windows of $\mu$ where different phases of the two models are observed are the same.

It is worth mentioning that the results from~\cite{bouchot2024scaling} require $\mu\gg(\log n)^{3/2}/\sqrt{n}$.
This threshold is quite natural when studying such models, as it arises when the leading term of the partition function, of order $e^{(\mu^2n)^{1/3}}$, becomes of the same order as any power of $n$ or $\mu$.
To circumvent this issue, we end up relying on very precise asymptotics for the number of trees of a given height.
In particular, we refer the reader to Theorem~\ref{thm:ht-bdd tree} where we provide the exact asymptotic behaviour of the number of trees of size $n$ and height $<m$:
\begin{align}\label{eq:Hnm}
    H_{n,m}=\Big|\big\{T\in\bbt_n:h(t)<m\big\}\Big|\,.
\end{align}
While the size of this set has been well-studied, with a closed form provided in~\cite{de1972average}, the best bound we could extract from the literature only covered the first order logarithm limit of $H_{n,m}$ (for example by applying the Mogul’ski\u{i} estimate~\cite{mogul1975small}), which only provided results when $\mu\gg(\log n)^{3/2}/\sqrt{n}$ and not all the way to $\mu\gg1/\sqrt{n}$.
\color{black}

The height-biased model where the weight function is power-like, i.e., where the weight of a tree $T$ of size $n$ is proportional to $h(T)^\alpha$ for $\alpha$ real, has been studied in \cite{durhuus2022trees}. There, the bias turns out to be not strong enough to affect the local limit even though the asymptotics of the partition function differ from the uniform case. The extension of these results to the case where $\alpha=\alpha_n$ depends on the size of the tree, remains open. 

Before concluding this section, we mention some existing models of biased/weighted random trees, even though they exhibit different characteristics from ours. A well-known class of trees are the recursive trees (RT) that are constructed by starting from a root of label 1 and joining new (labelled) vertices step by step. In \cite{leckey2020height}, the authors introduced \emph{depth-weighted random recursive trees} where the probability of joining a new leaf to an existing vertex is proportional to its depth in the tree. They also observe a change in the height of the tree when the weight function is exponential in depth, whereas the order of the height does not differ from the uniform case when the weight function is power-like, similar to phenomena we observe between \cite{durhuus2022trees} and \cite{durhuus2003probabilistic}. {More recently, \cite{dicaprio} studies recursive tree models where there is a preference for the new leaf to be joined towards \emph{more recently arrived} vertices, thus indirectly making the obtained trees \emph{taller} than when there is no such preference (i.e. taller than the uniform random recursive tree (URRT) \cite{pittel1994note}). Indeed, the authors exhibit both logarithmic and polynomial height regimes for such trees.}
In another recent work \cite{thevenin2023local} the authors studied a model of random trees that interpolates (in terms of the asymptotic form of the average path length, which is the sum of all distances to the
root) between uniformly random rooted labelled trees (a sub-class of {the so-called \emph{simply generated trees} \cite{janson2012simply}}) and recursive trees. There, an edge with the vertex closer to the root having label larger than the label of the other vertex it is adjacent to is called a {\em descent} and  the weight of a tree of size $n$ is proportional to a given parameter $q \in (0,1]$ to the power of the number of descents; when $q=1$, they recover the model of uniformly random labeled rooted trees.

\paragraph{Organization of the rest of the paper.} {In Section \ref{sec:prelims} we introduce some standard and useful notions related to viewing trees as excursions. Then in Section \ref{sec:global_result_proofs} we provide the proofs of Theorems \ref{thm:ldp} and \ref{thm:width_scaling_Brownian} for global properties in the Brownian regime. In Section \ref{sec:inter_height}, we provide all the results related to the height of the tree in the non-Brownian regimes, i.e., proofs of Theorems \ref{thm:inter lln}, \ref{thm:inter clt}, \ref{thm:discrete clt}, \ref{thm:ldp}, \ref{thm:width} and \ref{thm:high lln}. In Section \ref{sec:loc_results}, we provide the proofs of results on local behavior around the root, i.e., the proofs of Theorems \ref{thm:loclim} and \ref{thm:root deg}. In Appendix \ref{sec:func_poly_fac} we prove a technical result on functions of certain exponential form. We use this result first in Appendix \ref{sec:ht-bdd tree} to derive properties of the numbers $H_{n,m}$ as defined in \eqref{eq:Hnm}, and later in Appendices \ref{sec:on part fun}, \ref{sec:part fun 1} and \ref{sec:part fun 2}, where we provide the explicit analysis that derives the asymptotic behavior of the partition function as defined in \eqref{eq:part fun}.
}

\section{Encoding trees with walks}\label{sec:prelims}

In this section, we discuss \emph{encodings} of random trees, using random walk excursions. These universal concepts will aid us in proving the scaling limit result of Theorem \ref{thm:scaling limit}, as well as the result on the width in the non-Brownian regime, Theorem \ref{thm:width_scaling_Brownian}. 

For any integer $k>0$, let
\begin{align*}
    \cW^{(k)}:=\Big\{\big(w_t:0\leq t\leq k\big):w_0=0\textrm{ and }\big|w_t-w_{t+1}\big|=1{\textrm{ for }0\le t \le k-1}\Big\}
\end{align*}
be the set of \emph{walks} starting at $0$ of length $k$ and step size $1$.
Let us further define
\begin{align*}
    \cB^{(k)}_{x,y}:=\Big\{\big(x+w_t:0\leq t\leq k\big):w\in \cW^{(k)},x+w_k=y\textrm{ and }\big|w_t-w_{t+1}\big|=1{\textrm{ for }0\le t \le k-1}\Big\}
\end{align*}
for the set of \emph{bridges} of length $k$ and going from $x$ to $y$, and
\begin{align*}
    \cE^{(k)}:=\Big\{\big(w_t:0\leq t\leq k\big):w\in \cW^{(k)},w_k=-1, w_t\geq0\textrm{ for }0\leq t\leq k-1 \textrm{ and }\big|w_t-w_{t+1}\big|=1{\textrm{ for }0\le t \le k-1}\Big\}
\end{align*}
for the set of \emph{excursions} of length $k$.
It is worth noting that the set of bridges is empty whenever $k$ and $y-x$ do not have the same parity, which also means that the set of excursions is empty whenever $k$ is not odd.
In the case of bridges, when $x=0$ we simply denote it with $\cB^{(k)}_y$.
Since we will eventually be interested in scaling these objects and taking their limit, we equivalently see walks, bridges and excursions as their sequence of integer values or as continuous functions on $[0,k]$ where we linearly interpolate between the values at integer arguments.

\begin{definition}[Contours]\label{def:contour}
 For a plane tree $T$ with $n$ vertices, the contour walk $(C^{(T)}_t;0\leq t \leq 2n-1)$ of $T$ is an element of $\cE^{(2n-1)}$, that traces the contour of the tree in the depth-first manner from left to right, starting at $0$ and finishing at $2n-1$. We refer the reader to Figure \ref{fig:contour} for a visual representation.  
\end{definition}
\begin{figure}[htb]
    \centering
    \includegraphics[width=\textwidth]{contour.pdf}
    \caption{A representation of the bijection between a rooted plane tree with $n$ nodes and its contour --- an excursion of length $2n-1$. Each step of the excursion is represented by an arrow, both on the walk and on the tree, and corresponds to a single side of an edge (left or right). To highlight this relation, the nodes of the tree are shown at the bottom of the walk in the order in which they would be encountered when contouring around the tree.
    }
    
    \label{fig:contour}
\end{figure}

The set of all rooted plane trees of size $n$ is in bijective correspondence with the set $\cE^{(2n-1)}$, with a tree simply being mapped to its contour walk. For {the} random trees $\mathbf{T}_n$ or $\T_n$, (recall $\T_n$ is the uniform tree, i.e., the random tree with law \eqref{eq:T dist} corresponding to $\mu\equiv 0$) the contours {are} random variables taking values in $\cE^{(2n-1)}$. In particular, thanks to the bijection mentioned above, the contour $(C^{(\T_n)}_t;0\leq t\leq 2n-1)$ of the random tree $\T_n$, is a uniformly random element of $\cE^{(2n-1)}$.

\begin{remark}[Heights and maximums]\label{rem:ht_max_contour}
    It is worth mentioning that the height of a tree exactly corresponds to the maximal value of its corresponding contour walk seen as a continuous function.
\end{remark}

Sometimes we will translate events on bridges to events on excursions, and thus events on trees, which are a bit easier to understand. For this, we will need the \emph{cycle lemma}. The result is relatively standard, and we refer the reader to \cite[Lemma 15.3]{janson2012simply} for a proof. 
\begin{lemma}[Cycle lemma]\label{lem:cycle}
    Fix $n\geq 1$. There is a map $\phi:\cB^{(2n-1)}_{-1}\to \cE^{(2n-1)}$ from the set of all bridges starting $0$ and ending at $-1$, to the set of excursions, both of length $2n-1$, such that $\phi$ is a $(2n-1)$ to $1$ map.
\end{lemma}

Bridges that start at $0$ and end at $-1$ are called \emph{standard}. Sometimes it will be important to work with standard bridges. We will do this using the following result, that transforms general bridges starting at $0$ to standard bridges with a constraint.

\begin{lemma}[Cut bijection lemma]\label{lem:cut_bij}
    Let $n\geq1$ and $x\in\Z$.
    There is a function $\phi_x$ that bijectively sends the set $\cB^{(n)}_x$ of bridges from $0$ to $x$ onto
    \begin{align*}
        \tilde{\cB}^{(n)}_x:=\cB_{-\mathbb{1}_{\{n\textrm{ is odd}\}}}^{(n)}\cap\Big\{\big(w_t:0\leq t\leq n\big):\exists t,w_t=\floor{x/2}\Big\}
    \end{align*}
    of all bridges starting at $0$, ending at $-1$ when $n$ is odd or at $0$ when $n$ is even, and passing through $\floor{x/2}$. 
\end{lemma}

\begin{figure}[H]
    \centering
    \includegraphics[width=0.98\textwidth]{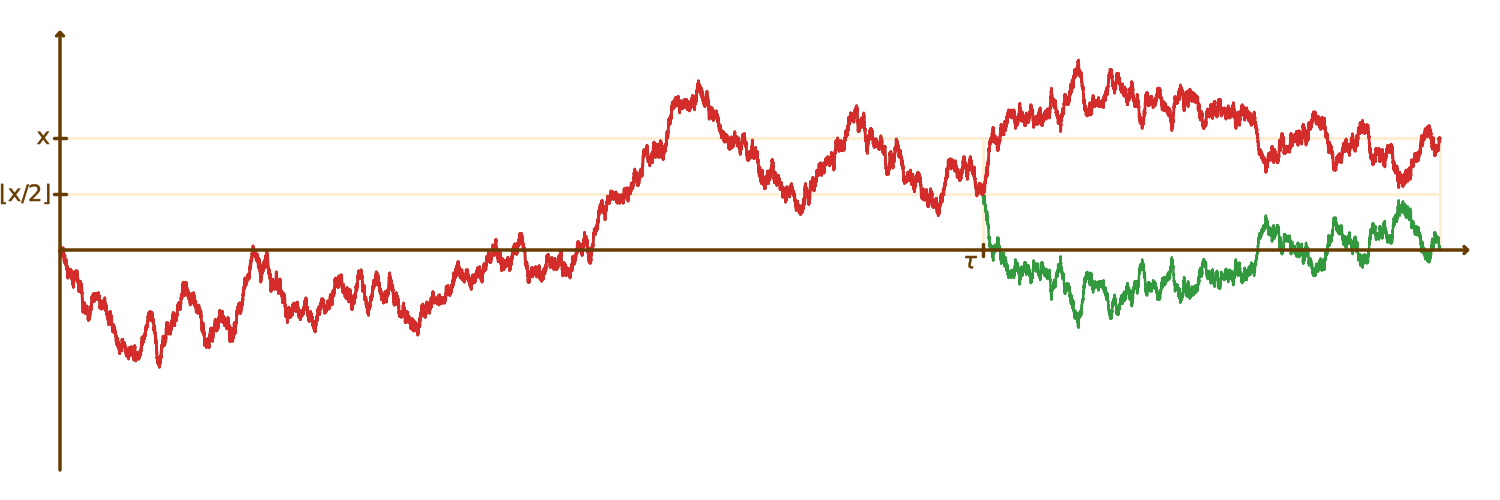}
    \caption{In the above picture for an element $w$ of $\cB^{(n)}_x$, displayed as the red path, we apply the map $\phi_{x}$ to reflect the segment of the path $w$ from the last passage $\tau$ at $\lfloor x/2 \rfloor$ to $n$ with respect to $y=\floor{x/2}$ to obtain the green segment. The image $\phi_{x}(w)$ is thus the concatenation of the red path up to time $\tau$ together with the green segment. }
    \label{fig:cutbij}
\end{figure}

\begin{proof}
We proceed by constructing such a function $\phi_x$. Let $w\in\cW^{(n)}$ be a walk such that there exists $t$ satisfying $w_t=\floor{x/2}$.
For such a $w$, let $\tau$ be its time of last passage from $\floor{x/2}$, i.e., 
\begin{align*}
    \tau:=\sup\left\{t\geq 0:w_t=\floor{\frac{x}{2}}\right\}\,.
\end{align*}
Then $\phi_x$ is the map such that $\phi_x(w)$ coincides with $w$ on $0 \leq t \leq \tau$ and is the reflection of $w$ around the horizontal axis $y=\floor{x/2}$ on $\tau \leq t \leq n$. In other words, seeing a walk as a sequence of $\pm1$ steps, $w$ and $\phi_x(w)$ share the same first $\tau$ steps and then have opposite steps.
We refer to Figure~\ref{fig:cutbij} for a representation of this function, where the original walk, in red and from the set $\cB^{(n)}_x$, is transformed into its green reflection which belongs to $\tilde{\cB}^{(n)}_x$.

As a first observation, for any $w$ such that there exists $t$ satisfying $w_t=\floor{x/2}$, we see that the last hitting time of $\floor{x/2}$ of $w$ and $\phi_x(w)$ is the same.
Indeed, $w_\tau=\floor{x/2}$ so we also have $\phi_x(w)_\tau=\floor{x/2}$ and since the two walks are the reflection of each other around $y=\floor{x/2}$ on $(\tau,n]$, each one intersects the axis if and only if the other one does.
This remark implies that $\phi_x\cdot\phi_x(w)=w$, meaning that $\phi_x$ is a symmetry.

Assume now that $w$ is such that $w_n=x$.
In that case, we necessarily have some $t$ such that $w_t=\floor{x/2}$, so $\phi_x(w)$ is well-defined.
Since $w$ and $\phi_x(w)$ coincide on $[0,\tau]$ and the steps of $w$ from $\tau$ to $n$ go from $\floor{x/2}$ to $x=\floor{x/2}+(x-\floor{x/2})$, the steps of $\phi_x(w)$ from $\tau$ to $n$ go from $\floor{x/2}$ to $\floor{x/2}-(x-\floor{x/2})=2\floor{x/2}-x$.
Recalling that $x$ has to have the the same parity as $n$, we see that $\phi_x(w)$ indeed ends at $0$ when $n$ is even and at $-1$ when $n$ is odd.
Conversely, if $w_n=-\mathbbm{1}_{\textrm{$n$ is odd}}$ and there exists $t$ such that $w_t=\floor{x/2}$, then $\phi_x(w)$ is a walk from $0$ to $2\floor{x/2}+\mathbbm{1}_{\textrm{$n$ is odd}}=x$.
This implies that $\phi_x$ sends $\cB^{(n)}_x$ into $\tilde{\cB}^{(n)}_x$ and vice-versa.
Using that $\phi_x\cdot\phi_x$ is the identity, this proves that $\phi_x$ bijectively maps to the sets.
\end{proof}

\section{Proof of global results in the Brownian regime}\label{sec:global_result_proofs}
In this section, we prove the scaling limit results in the Brownian regime, i.e., Theorems \ref{thm:scaling limit} and \ref{thm:width_scaling_Brownian}. We begin with Section \ref{sec:def_realtrees_GHP} discussing some relevant concepts, and then move on to the actual proofs.

\subsection{Gromov-Hausdorff-Prokhorov convergence and real trees}\label{sec:def_realtrees_GHP}
For more details on what follows, we refer the reader to \cite[Section 3]{gall2011scaling}. Informally, one can think of real trees as a connected collection of line segments in the plane with no loops. They can be constructed analogously to the discrete case as follows. Let $C([0,1], \R _+ )$ be the space of continuous functions from $[0,1]$ into $\R_+$ equipped with the topology of uniform convergence. Given a continuous function $g$ such that $g(0)=g(1)=0$, the construction of the corresponding real tree is standard: for every $s,t \in [0,1]$, set

$$
d_g(s,t) = g(s)+g(t)-2m_g(s,t)
$$ where $$m_g(s,t)= \inf _{z\in [s\wedge t,s\vee t]} g(z) \, .$$ Set $T_g = [0,1] /\!\!\sim$ where $s \sim t$ if and only if $d_g(s,t)=0$. Note that $d_g$ induces a distance on $T_g$, which we denote by $d_g$ as well and is in fact a pseudo-metric. The canonical projection $p_g : [0,1] \to T_g$ is continuous hence $T_g$ is compact. The compact metric space $(T_g,d_g)$ is then a real tree, which we view as a rooted real tree with the root $\rho=p_g(0)=p_g(1)$. 

In order to make sense of the convergence of rooted discrete trees towards real trees, we introduce the Gromov-Hausdorff-Prokhorov (GHP) distance on the set $\mathbb K$ of all isometry classes of {rooted} compact measured metric spaces, where a \emph{rooted} (also commonly referred to as \emph{pointed}) metric space is a metric space with a distinguished \emph{root}.

Let $(E,d)$ be a metric space with its Borel $\sigma-$algebra $\mathcal{B}(E)$. The Hausdorff distance $d_H$ between two non-empty subsets $X,Y$ of $E$ is defined as
$$
d_H(X,Y) = \inf \Big\{ \varepsilon>0 : X \subset U^\varepsilon(Y)~ \text{and} ~Y \subset U^\varepsilon(X) \Big\}\,,$$
where $U^\varepsilon(X) := \{x\in E: d(x,X)<\varepsilon \}$.
Let $(E_1,\rho_1,d_1)$ and $(E_2,\rho_2,d_2)$ be two {rooted} compact metric spaces with metrics $d_1$ and $d_2$ and with {root} $\rho_1$ and $\rho_2$, respectively. The Gromov-Hausdorff distance $d_{GH}$ between $E_1$ and $E_2$ is defined by
$$ d_{GH}(E_1,E_2)=\inf_{\phi_1,\phi_2} \Big\{d_H\big(\phi_1(E_1),\phi_2(E_2)\big) \vee d\big(\phi_1(\rho_1),\phi_2(\rho_2)\big)\Big\}$$ where the infimum is over all possible choices of metric spaces $({E},{d})$ and the isometric embeddings $\phi_1:E_1 \to E$ and $\phi_2:E_2 \to E$ of $E_1$ and $E_2$ into ${E}.$
Denoting by $P(E)$ the collection of all probability measures on $(E,\mathcal B (E))$, the Prokhorov distance $d_P: P(E) \times P(E) \to [0, \infty)$ between two Borel probability measures $\mu$ and $\nu$ on $E$ is given by
$$
d_P(\mu,\nu)= \inf \Big\{ \varepsilon>0 : \mu(A) \leq \nu \big(U^\varepsilon(A)\big) + \varepsilon ~\text{and} ~\nu(A) \leq \mu\big(U^\varepsilon(A)\big) + \varepsilon , ~\forall A \in \mathcal B (E) \Big\} \, .$$
Finally, the GHP distance between $(E_1,\rho_1,d_1,\mu_1)$ and $(E_2,\rho_2,d_2,\mu_2)$, where $\mu_i$ is a probability measure on $E_i$ for $i=1,2$, is given by
$$ d_{GHP} (E_1,E_2) = \inf \Big\{ d_{GH}\big(\phi_1(E_1) , \phi_2(E_2)\big) \vee d_P\big(\phi_{1*} \mu_1 , \phi_{2*} \mu_2\big) \Big\}$$
where the infimum is again over all possible choices of metric spaces $({E},{d})$ and the isometric embeddings $\phi_1:E_1 \to E$ and $\phi_2:E_2 \to E$ of $E_1$ and $E_2$ into ${E}$ and $\phi_*$ stands for the pushforward measure. For more details on GH and GHP distances and the topologies they induce, we refer the reader to \cite{burago2001course, miermont2013brownian, gall2011scaling}.

\begin{definition}[Brownian Continuum Random Tree]
    The Brownian Continuum Random Tree (BCRT in short) \cite{aldous1991continuum} is the random real tree $(T_{\mathbf{e}},\rho_{\mathbf{e}},d_{\mathbf{e}},\sigma_{\mathbf{e}})$, where $(\mathbf{e}(t):0\leq t\leq 1)$ is the standard Brownian excursion, and $\sigma_{\mathbf{e}}$ is the pushforward of the Lebesgue measure on $[0,1]$ under the projection map of the quotient $\sim$.
\end{definition}
For more discussion on this object and the Brownian excursion, we refer the reader to \cite{duquesne2002random}. The conclusion of Theorem \ref{thm:scaling limit} says that when $\mu=O(1/\sqrt{n})$, the scaling limit of the tree $\mathbf{T}_n$ is obtained by appropriately tilting the above defined BCRT.

Finally, let us state without proof two useful results that we will need in the proof of Theorem \ref{thm:scaling limit}. The first one is \cite[Corollary 3.7]{gall2011scaling}. {Recall that, for any function $f$, $\|f\|$ denotes its maximum value.}

\begin{proposition}\label{prop:gh_fromabove_supnorm}
    Let $\cT_g$ and $\cT_{g'}$ be two real trees encoded by excursions $g,g':[0,1]\to \R_+$. Then
    \begin{align*}
        d_{GH}(\cT_g,\cT_{g'})\leq 2\big\|g-g'\big\|.
    \end{align*}
\end{proposition}

The following proposition is \cite[Proposition 4.14]{duquesne2022scaling} adapted {for our purposes}. 
\begin{proposition} \label{GHP_supnorm}
Let $d_\cT$ and $d_{\cT'}$ be two pseudo-metrics  on $[0,1]$ with $\cT$ and $\cT'$ induced real trees equipped with the push-forward of the Lebesgue measure on $[0,1]$. Then
$$d_{GHP}(\cT,\cT') \leq \frac 3 2 \big\| d_\cT - d_{\cT'} \big\| \, .$$
\end{proposition}

\subsection{Scaling limit in the Brownian regime}

In this section, we provide the proof of theorem \ref{thm:scaling limit}. Recall the discussion on real trees and GHP convergence from Section \ref{sec:def_realtrees_GHP}. In what follows, whenever the dominated convergence theorem is evoked, we implicitly use Skorokhod's embedding theorem (e.g., see \cite{klenke2008probability}) to couple our random variables in a common probability space and pass to almost sure convergence.

\begin{proof}[Proof of Theorem \ref{thm:scaling limit}]

Recall the definition of the contour of a tree from Definition \ref{def:contour}. Let $\T_n$ be a \textit{uniformly} sampled tree, and let $\C^{(\T_n)}$ be the rescaled contour walk of the tree $\T _n$ obtained by applying the Brownian scaling on the original contour walk on $[0,2n-1]$, i.e.
$$\C^{(\T_n)}(t)=\frac{1}{\sqrt{n}}C^{(\T_n)}\big((2n-1)t\big)$$
for $t \in [0,1]$, where $C^{(\T_n)}$ is the contour of $\T_n$.
Let us also denote the set of such rescaled contour walks of length $2n-1$ by $\cC^{(2n-1)}$.
We note that $\cC^{(2n-1)}$ is a finite subset of $C([0,1],\R_+)$ and further, since the association of trees to their respective contours is bijective, we note that the law of $\C^{(\T_n)}$ is uniform on $\cC^{(2n-1)}$.

Consider the contour walk $C^{(\bT_n)}$ of the $\mu$-height-biased tree $\bT_n$. Note that by the definition of the distribution \eqref{eq:T dist} and Remark \ref{rem:ht_max_contour}, the law of $C^{(\bT_n)}$, has a Radon-Nikodým derivative with respect to the law of $C^{(\T_n)}$, that is proportional to the function $\exp(-\mu \|f\|)$, on the space of all discrete contours of length $2n-1$. As a consequence, the contour $\C^{(\bT_n)}$ obtained after applying Brownian scaling on $C^{(\bT_n)}$, i.e.,
$$\C^{(\bT_n)} (t)=\frac{1}{\sqrt{n}}C^{(\bT_n)} \big( (2n-1)t) \big)$$
for $t \in [0,1]$, has a Radon-Nikodým derivative with respect to the rescaled uniform contour $\C^{(\T_n)}$, that is proportional to $\exp{(-\mu \sqrt{n} \|f\|)}$. Since the law of $\C^{(\T_n)}(t)$ is uniform on $\cC^{(2n-1)}$, this is equivalent to simply saying that
\begin{align*}
    \Prob{\C^{(\bT_n)} = \C} \propto \exp\Big(-\mu \sqrt{n}\big\|\C\big\|\Big)
\end{align*}
over all $\C \in \cC ^{(2n-1)}$.
{These observations at hand, let us first} show that $\C^{(\bT_n)}$ converges in distribution to $\mathbf{e}_{\alpha}$ in the supremum norm topology, {$\mathbf{e}_{\alpha}$ is the tilted Brownian excursion as defined right before Theorem \ref{thm:scaling limit}}. {Endowing $C([0,1],\R_+)$ with the supremum topology,} by the Portmanteau theorem (see \cite{klenke2008probability}), it is sufficient to show that for any continuous bounded function
$\Lambda: C([0,1],\R_+) \to \R$,
$$\E\left[\Lambda\left(\C^{(\bT_n)}\right)\right]= \frac{\Exp{\Lambda(\C^{(\T_n)}) e^{-\mu \sqrt{n}  \left\|\C^{(\T_n)}\right\|}}}{\Exp{e^{-\mu \sqrt{n} \left\|\C^{(\T_n)}\right\|}}} \longrightarrow \Exp{\Lambda(\mathbf{e}_{\alpha})} \, .$$ 

It is standard that, in the uniform case, we have convergence of contours $\C^{(\T_n)}$ towards the standard Brownian excursion $(\mathbf{e}(t))_{t \in [0,1]}$ on $[0,1]$, see e.g., \cite{aldous1991continuum}. Since $\| \C^{(\T_n)} \|$ admits sub-Gaussian tails by \cite[Theorem 1.2]{LAB_LD_SJ_13}, we can deduce that  $\exp(-\mu \sqrt{n} \| \C^{(\T_n)}\|)$ is uniformly integrable.
Together with the boundedness of $\Lambda$ this implies that

$$
\frac{\E\left[\Lambda(\C^{(\T_n)}) e^{-\mu \sqrt{n} \left\| \C^{(\T_n)}\right\|}  \right]}{\Exp{e^{-\mu \sqrt{n} \left\| \C^{(\T_n)}\right\|} }} \longrightarrow \frac{\E \left[ \Lambda (\mathbf{e}) e^{-\alpha \|\mathbf{e}\|}\right]}{\Exp{e^{-\alpha \|\mathbf{e}\|}}}
$$
which itself implies that
$$
\E\left[\Lambda\left(\C^{(\bT_n)}\right) \right]\longrightarrow \Exp{\big.\Lambda(\mathbf{e}_{\alpha})} \, 
$$
where the law $\rho_\alpha$ of $\mathbf{e}_\alpha$ satisfies     
\begin{align*}
 \frac{d\rho_{\alpha}}{d\rho_{0}}(f)=\frac{e^{-\alpha \|f\|}}{\int e^{-\alpha \| g \| }d\rho_0 (g)} \, ,
 \end{align*}
{with $\rho_0$ denoting the law of $\mathbf{e}$. In particular, $\C^{(\bT_n)}$ converges in the topology of uniform convergence to $\mathbf{e}_\alpha$. Together with \cref{prop:gh_fromabove_supnorm}, the convergence of $\C^{(\bT_n)}$ to $\mathbf{e}_\alpha$  implies the convergence in distribution of the corresponding real tree $\ccT^\mu _n$ to the continuum tree $\ccT_\alpha$ encoded by $\mathbf{e}_\alpha$, in the Gromov-Hausdorff sense. Let us now discuss how to get the improvement to the Gromov-Hausdorff-Prokhorov (GHP) sense of convergence.}

For this improvement, let us first obtain some bounds in the deterministic setting. Fix some tree $T_n$ of size $n$ with the standard graph distance $d_{gr}$ and consider the corresponding real tree $\cT_n$ generated by its contour function with the distance function $d_{\cT_n}$ assigning unit length to every edge. The discrete tree can be embedded in the real one isometrically: see for example \cite[Example 4.9]{duquesne2022scaling} for this standard construction. Further, denote by $a\cT_n$ the real tree with distances scaled by $a$. Below, we will compute the Prokhorov distance between certain measures defined on $a\cT_n$. When we say the push-forward of some measure originally defined on the discrete tree $T_n$, it is by the composition of these two maps to the real tree $a\cT_n$.

Given $T_n$, recall that we write $\deg_{T_n}(v)$ for the degree of the vertex $v \in T_n$. The degree-biased measure $\sigma_{T_n} ^*$ on $T_n$ is the measure defined on the vertex set $V(T_n)$ of $T_n$ satisfying
$$
\sigma_{T_n} ^*(v) = \frac{\deg_{T_n} (v)}{\sum \limits_{v \in T_n } \deg_{T_n} (v)} = \frac{\deg_{T_n} (v)}{2(n-1)} \, .
$$

First, we bound the Prokhorov distance $d_{\mathrm{Prok}}(\mathrm{Leb} _{\cT_n},\sigma_{\cT_n} ^*)$ between the push-forward of the degree-biased measure $\sigma _{\cT_n} ^* $ viewed as a measure on $a\cT_n$ and the push-forward of Lebesgue measure $\mathrm{Leb} _{\cT_n}$ from $[0,1]$ to $a\cT_n$. For this purpose, first note that $2(n-1)$ is two times the number of edges of $T_n$. As each edge in the tree $T_n$ corresponds to one up and one down step in the corresponding contour $\C^{(T_n)}$, we will call the set of points in the real tree $\cT_n$ (equivalently in its rescaled version) induced by those two steps an \emph{edge} as well. Hence the measure $\mathrm{Leb} _{\cT_n}^*$ assigns mass $1/(n-1)$ to each edge of $a\cT_n$.

For a metric space $(E',d')$ and $U \subset E'$, {recall} the $\varepsilon$ neighborhood of $U$ {is} defined {as} $U^{\varepsilon}=\{x:d'(x,U)<\varepsilon\}$. It follows that for all $U\subset a\cT_n$, $\mathrm{Leb} _{\cT_n} (U)\leq \sigma_{T_n} ^*(U^{\frac a 2})+ \frac{a}{2}$ with the same inequality valid for $\mathrm{Leb}_{\cT_n} ^* \leftrightarrow \sigma_{\cT_n} ^*$. Hence 

\bb \label{prok1}
d_{\mathrm{Prok}}(\mathrm{Leb}_{\cT_n} ^*,\sigma_{\cT_n} ^*) \leq \frac{a}{2} \, .
\ee
In a second step, we compare $\sigma_{\cT_n} ^*$ with the push-forward of the uniform measure on $T_n$ to $a\cT_n$ which we denote by $\sigma_{\cT_n}$. Note that both the uniform and size-biased measures defined on $T_n$ are supported on the vertices and denote the image of the vertex set $V(T_n)$ of $T_n$ in $a\cT_n$ under the aforementioned canonical embedding by $V(a\cT_n)$. Fix any subset $S \subset V(a\cT_n)$ and set $S_+ := \{v \in V(a\cT_n):d(v,S) \leq a \}$. Writing $e(S)$ for the set of edges having at least one end-point in $S$, one has $e(S) \geq |S|-1$, which can easily be seen by counting the edge that connects every $v\in S$ to its ancestor, except possibly the root. Defining the degree $\deg_{\cT_n}(v)$ of a vertex $v \in V(a\cT)$ as the degree of its pre-image in the discrete tree $T_n$ we have
\begin{align*}
    \sigma_{\cT_n} ^*(S_+) = \sum \limits_{v \in S_+} \sigma_{\cT_n} ^* (v) & = \frac{\sum \limits_{v \in S_+} \deg_{\cT_n} (v)}{2(n-1)} \geq \frac{2|e(S)|}{2(n-1)} \geq \frac{|S|-1}{n} = \sigma_{\cT_n}(S) - \frac{1}{n}
 \end{align*}
Note that $S_+$ contains the same vertices as the $(1+\delta)a$-neighborhood of $S$ for some $\delta<1$ small. Hence the last inequality implies
$$
\sigma_{\cT_n} (S)\leq \sigma_{\cT_n} ^*\big(S^{(1+\delta)a}\big) + \frac 1 n \, .
$$
Then for $ (1+\delta)a\geq1/n$, we have
$$\sigma_{\cT_n} (S)\leq \sigma_{\cT_n} ^*\big(S^{(1+\delta)a}\big)+(1+\delta)a$$
and for $(1+\delta)a < 1/n$, we have
$$\sigma_{\cT_n}^*\big(S^{(1+\delta)a}\big) \leq \sigma^*\big(S^{1/n}\big)~\implies~\sigma_{\cT_n}^*(S) \leq \sigma_{\cT_n}^*\big(S^{1/n}\big) + \frac 1 n\,.$$
Furthermore, $\sigma_{\cT_n}^*(S)\leq \sigma_{\cT_n}(S^{(1+\delta)a})$ trivially holds. Hence we conclude that
\bb \label{prok2}
d_{\mathrm{Prok}}\big(\sigma_{\cT_n}^*,\sigma_{\cT_n}\big) \leq \max\left((1+\delta)a,\frac 1 n \right) \, .
\ee
It now follows by an application of triangle inequality to \eqref{prok1} and \eqref{prok2} that 
\bb \label{prok3}
d_{\mathrm{Prok}}\big(\mathrm{Leb}_{\cT_n} ^*, \sigma_{\cT_n}\big) \leq \frac a 2 + \max\left((1+\delta)a,\frac 1 n \right) \, .
\ee
Finally, using \cref{GHP_supnorm}, we have for any sequence $\cT_n \to \cT$, 
\bb \label{prok4}
d_{\mathrm{Prok}}\big(\mathrm{Leb}_{\cT_n} ^*, \mathrm{Leb}_{\cT} ^*\big) \leq d_{GHP}\big(\cT_n, \cT\big) \leq \frac 3 2 \left\| d_{C^{(\cT_n)}}-d_{C^{(\cT)}} \right\| \leq 6\left\| C^{(\cT_n)} -C^{(\cT)} \right\| \, . 
\ee
Together with the GH convergence, \eqref{prok3} with $a$ chosen to be $1/\sqrt{n}$ and \eqref{prok4} imply that
\begin{align*}
\left(\ccT_n ^\mu, \frac{d_{\ccT_n ^\mu} }{\sqrt{n}}, \sigma_{\ccT_n ^\mu}\right) {\toindis} (\ccT_\alpha, d_\alpha, \sigma_\alpha)
\end{align*}
as $n\to \infty$ in the GHP {topology}.
Finally, using the canonical embedding of the discrete tree $T_n$ into its corresponding real tree $\cT_n$ one has
\begin{align*}
d_{GHP}\left( \left(T_n, \frac{d_{gr}}{\sqrt{n}}, \sigma_n \right), \left(\cT_n, \frac{d_{\cT_n}}{\sqrt{n}}, \sigma_{\cT_n} \right) \right) \leq \frac{1}{\sqrt{n}} \, .
\end{align*}
Taking $T_n=\bT_n$ and $\cT_n=\ccT_n ^\mu$ in the above display, the last two displays together imply the desired convergence in the statement of the theorem and conclude the proof.
\end{proof}

\subsection{Scaling limit of the width in the Brownian regime}\label{sec:proof_width_brownian}

In this section, we prove Theorem \ref{thm:width_scaling_Brownian}. {Recall for any function $f$ we denote its supremum value by $\|f\|$.} The key result is that if we consider the uniform plane tree $\T_n$ on $n$ vertices, and denote the corresponding height by $H_n$ and width by $W_n$, then one has {the joint convergence} (see e.g., \cite[Corollary 2.5.1]{duquesne2002random}) 
\begin{align*}
    \left(\frac{H_n}{\sqrt{n}},\frac{W_n}{\sqrt{n}} \right) \dlim \left(\|\mathbf{e}\|, \frac{1}{2}\|\ell\|\right), \numberthis \label{eq:joint_ht_wt}
\end{align*}
where $(\ell(x):x\geq 0)$ is the \textit{local-time process} of the standard Brownian excursion $(\mathbf{e}(t):0\leq t\leq 1)$, defined by
\begin{align*}
    \int_0^a \ell(x)dx=\int_0^1\indE{[0,a]}\big(\mathbf{e}(s)\big)ds, \numberthis \label{eq:loc_time_def}
\end{align*}
for any $a>0$. 

\begin{proof}[Proof of Theorem \ref{thm:width_scaling_Brownian}]

Assume that $\mu \sqrt{n} \to \alpha \in \R$, and note that by the definition of our measure \eqref{eq:T dist}, for any bounded continuous function $F:\R \to \R$, one has
\begin{align*}
    \Exp{F\left(\frac{w(\mathbf{T}_n)}{\sqrt{n}}\right)}=\frac{\Exp{F\left(\frac{W_n}{\sqrt{n}}\right)\exp\left(-\mu\sqrt{n}\frac{H_n}{\sqrt{n}}\right)}}{\Exp{\exp\left(-\mu\sqrt{n}\frac{H_n}{\sqrt{n}}\right)}}.
\end{align*}
Observe that as $n \to \infty$, if $\mu\geq 0$, using \eqref{eq:joint_ht_wt} we immediately obtain as a consequence of dominated convergence theorem that
\begin{align*}
    \Exp{F\left(\frac{w(\mathbf{T}_n)}{\sqrt{n}}\right)}\to \frac{\Exp{F\left(\frac{1}{2}\|\ell\|\right)\exp(-\alpha\|\mathbf{e}\|)}}{\Exp{\exp(-\alpha\|\mathbf{e}\|)}}.\numberthis \label{eq:brownian_width_lim_1}
\end{align*}
A uniform integrability argument similar to the one used in the proof of Theorem \ref{thm:scaling limit} gives also the same conclusion for $\mu<0$ satisfying $\mu \sqrt{n} \to \alpha \in [0,\infty)$.

Recall the well known identity of Jeulin, (see e.g., \cite[Eq. (2.8)]{chassaing2000height})
\begin{align*}
    \left(2\|\mathbf{e}\|,\frac{1}{2}\|\ell\| \right)\stackrel{d}{=} \left(\int_0^1\frac{ds}{\Tilde{\mathbf{e}}(s)},\|\Tilde{\mathbf{e}}\|  \right), \numberthis \label{eq:Jeulin}
\end{align*}
where $(\Tilde{\mathbf{e}}(s):0\leq 1\leq 1)$ is distributed as another standard Brownian excursion.

Using the above, {recalling the definition of $B_\alpha$ from the statement of Theorem \ref{thm:width_scaling_Brownian},} the RHS of \eqref{eq:brownian_width_lim_1} can be also written as,
\begin{align*}
    \frac{\Exp{F\left(\|\Tilde{\mathbf{e}}\|\right)\exp\left(-\frac{\alpha}{2}\int_0^1\frac{ds}{\Tilde{\mathbf{e}}(s)}\right)}}{\Exp{\exp\left(-\frac{\alpha}{2}\int_0^1\frac{ds}{\Tilde{\mathbf{e}}(s)}\right)}}=\Exp{\Big.F\big(\|B_\alpha\|\big)}\,,
\end{align*}
as claimed, {where {we} recall {the definition of the process} $B_\alpha$ with law $\xi_\alpha$ from \eqref{eq:law_balpha_xialpha}}.
\end{proof}

\begin{remark}[Local times of max tilted excursions]
    It is in fact possible to interpret the RHS of \eqref{eq:brownian_width_lim_1} as the maximum `local time' of a random excursion on $[0,1]$. Namely, consider the random excursion $(\mathbf{e}_{\alpha}(t):0\leq t \leq 1)$ with law $\rho_{\alpha}$ from Theorem \ref{thm:scaling limit}, and recall that it has a Radon-Nikodým derivative with respect to the standard Brownian excursion given by
    \begin{align*}
      f \mapsto  \frac{e^{-\alpha\|f\|}}{\int e^{-\alpha \|g\|}d\rho_0(g)},
    \end{align*}
    for any $f \in C[0,1]$, where $\rho_0$ denotes the law of a standard Brownian excursion $(\mathbf{e}(t):0\leq t\leq 1)$. Consider the local time process $(\ell(x);x\geq 0)$ of $\mathbf{e}(t)$ {as defined in \eqref{eq:loc_time_def}}. We aim to define an analogous process $(\ell_{\alpha}(x):x\geq 0)$ for the random excursion $\mathbf{e}_{\alpha}$. Observe that for any $0<x_1<x_2<\dots <x_m$, one can define the finite dimensional distributions as,
    \begin{align*}
        \Prob{\Big.\ell_{\alpha}(x_i)\in B_i,1\leq i \leq m}:=\frac{\Exp{\ind{\ell(x_i)\in B_i, 1\leq i \leq m}e^{-\alpha \|\mathbf{e}\|}}}{\Exp{e^{-\alpha \|\mathbf{e}\|}}},
    \end{align*}
where the $B_i$ are arbitrary Borel subsets of $\R$. The Kolmogorov consistency theorem, and the existence of the process $(\ell(x):x\geq 0)$ for the standard Brownian excursion $(\mathbf{e}(t):0\leq t\leq 1)$ then guarantees the existence of a \emph{unique} process $(\ell_{\alpha}(x):x\geq 0)$ with finite dimensional distributions as given in the last display. Furthermore, the identity $\int_0^a \ell(x)dx=\int_0^1\indE{[0,a]}(\mathbf{e}(s))ds$ implies a similar identity,
\begin{align*}
    \int_0^a\ell_{\alpha}(x)dx=\int_0^1\indE{[0,a]}\big(\mathbf{e}_{\alpha}(s)\big)ds,
\end{align*}
which means we can indeed interpret the RHS of \eqref{eq:brownian_width_lim_1} as the maximum local time of the random excursion $(\mathbf{e}_{\alpha}(t):0\leq t \leq 1)$.
\end{remark}

\section{The partition function, and asymptotics of the height in the non-Brownian regime}\label{sec:inter_height}

In this section we prove all the results regarding the behaviour of the height in the non-Brownian regime, i.e., $\mu\gg 1/\sqrt{n}$. These comprise the Theorems~\ref{thm:inter lln},~\ref{thm:inter clt},~\ref{thm:discrete clt},~\ref{thm:ldp}, and~\ref{thm:high lln}.
The main ingredient for all these proofs is the following detailed result characterizing the asymptotic behaviour of the partition function as defined in~\eqref{eq:part fun}. {Recall the function $\lambda_x(t)=xt+\log(1+\tan^2(\pi/t))$ and its unique minimizer $t_x$ from \eqref{eq:func_lambdax_minvalue}.}

\begin{theorem}\label{thm:part-fun}
    Let $\mu=\mu_n$ be such that $\mu\gg1/\sqrt{n}$. Then the partition function $Z_n$ as defined in~\eqref{eq:part fun} behaves according to the following regimes.
    \begin{itemize}
        \item[\numitem{1}] If $1/\sqrt{n}\ll\mu\ll n^{1/4}$, we have
        \begin{align*}
            Z_n=\big(1+o(1)\big)4^ne^\mu(e^\mu-1)\frac{\pi^{5/6}}{2^{1/3}3^{1/2}}\frac{\mu^{1/3}}{n^{5/6}}\exp\left(-3\left(\frac{\pi^2\mu^2n}{4}\right)^{1/3}\right)\,.
        \end{align*}
        \item[\numitem{2}] If $(\log^{3/2}n)/\sqrt{n}\ll\mu\ll n$, we have
        \begin{align*}
            Z_n=\big(1+o(1)\big)4^ne^{\mu}(e^\mu-1)\frac{\mu}{2n}e^{-n\lambda_{\mu/n}(t_{\mu/n})}\sum_{m\in\Z}\exp\left(-\left(\frac{3}{2}+o(1)\right)\left(\frac{\mu^4}{2\pi^2n}\right)^{1/3}(m-t_{\mu/n})^2\right)\,.
        \end{align*}
        \item[\numitem{3}] If $n^{1/4}\ll\mu$, we have
        \begin{align*}
            Z_n=\big(1+o(1)\big)4^ne^{2\mu}\sum_{m=\max\{\lfloor t_{\mu/n}\rfloor,3\}}^{\lceil t_{\mu/n}\rceil}\frac{\tan^2\left(\frac{\pi}{m}\right)}{m}e^{-n\lambda_{\mu/n}(m)}\,.
        \end{align*}
    \end{itemize}
    
\end{theorem}

\begin{remark}
    Some of the regimes from Theorem~\ref{thm:part-fun} overlap, but it can be easily verified that the formulas are coherent with each other.
    More precisely, the first two agree thanks to properties of Riemann's sum and the last two are directly compared in the proof of Proposition~\ref{prop:W mu-n}.
\end{remark}

\begin{proof}
    All of these results are covered in the different sections of the appendix.
    First of all, Lemma~\ref{lem:z and w} states that
    \begin{align*}
        Z_n=\big(1+o(1)\big)4^ne^\mu(e^\mu-1)W_n
    \end{align*}
    where $W_n$ is defined in~\eqref{eq:def W}.
    Now, the three regimes directly follow from Theorem~\ref{thm:W asympt}, Proposition~\ref{prop:W as a sum}, and Proposition~\ref{prop:W mu-n} respectively.
\end{proof}

We now use this theorem to prove the asymptotic behaviour of the height of the trees as covered in Theorem~\ref{thm:inter lln},~\ref{thm:inter clt},~\ref{thm:discrete clt},~\ref{thm:ldp}, and~\ref{thm:high lln}, respectively in Section~\ref{sec:lln},~\ref{sec:clt},~\ref{sec:clt},~\ref{sec:clt}, and~\ref{sec:bern}.
Before doing, it is worth highlighting the results from Proposition~\ref{prop:asympt lx} regarding $\lambda_x$ and $t_x$ which are repeatedly used throughout this section.
In particular, whenever $\mu\ll n$, the proposition states that
\begin{align*}
    t_{\mu/n}=\left(\frac{2\pi^2n}{\mu}\right)^{1/3}+o(1)
\end{align*}
and that
\begin{align*}
    n\lambda_{\mu/n}(t_{\mu/n})=3\left(\frac{\pi^2\mu^2n}{n}\right)^{1/3}+O\left(\left(\frac{\mu^4}{n}\right)^{1/3}\right)\,.
\end{align*}

\subsection{First order behaviour}\label{sec:lln}

This section focuses on proving Theorem~\ref{thm:inter lln} and provides a first and straightforward application of Theorem~\ref{thm:part-fun}.

\begin{proof}[Proof of Theorem~\ref{thm:inter lln}]
    To prove this result, we simply apply a Chernoff bound (or, exponential Markov inequality) with the formula for the partition function from Theorem~\ref{thm:part-fun}.
    In fact, we do not exactly prove the convergence stated in the theorem, but rather that
    \begin{align*}
        \frac{h(\bT_n)}{t_{\mu/n}}\longrightarrow1\,,
    \end{align*}
    where $t_{\mu/n}$ is the minimum of $\lambda_{\mu/n}$ as defined in~\eqref{eq:func_lambdax_minvalue}.
    Note that $t_{\mu/n}$ behaves as $(2\pi^2n/\mu)^{1/3}$ thanks to Proposition~\ref{prop:asympt lx}~\numitem{1}.
    We start by fixing $\epsilon>0$ and proving that
    \begin{align*}
        \Prob{\frac{h(\bT_n)}{t_{\mu/n}}\geq1+\epsilon}\longrightarrow0
    \end{align*}
    for all considered regimes of $\mu$.
    Observe that, for any (fixed) $\gamma\in(0,1)$, by a Chernoff bound we have
    \begin{align*}
        \Prob{\frac{h(\bT_n)}{t_{\mu/n}}\geq1+\epsilon}\leq\E\left[e^{\gamma\mu h(T_n)}\right]e^{-\gamma(1+\epsilon)\mu t_{\mu/n}}=\frac{e^{-\gamma(1+\epsilon)\mu t_{\mu/n}}Z_n^{(1-\gamma)\mu}}{Z_n^\mu}\,.
    \end{align*}
    Thus, if $1/\sqrt{n}\ll\mu\ll n^{1/4}$, Theorem~\ref{thm:part-fun}~\numitem{1} tells us that
    \begin{align*}
        \Prob{\frac{h(\bT_n)}{t_{\mu/n}}\geq1+\epsilon}\sim e^{-\gamma(1+\epsilon)\mu t_{\mu/n}}e^{-\gamma\mu}\frac{e^{(1-\gamma)\mu}-1}{e^\mu-1}\left(1-\gamma\right)^{1/3}\exp\left(3\left(\frac{\pi^2\mu^2n}{4}\right)^{1/3}-3\left(\frac{\pi^2(1-\gamma)^2\mu^2n)}{4}\right)^{1/3}\right)\,.
    \end{align*}
    Further use Proposition~\ref{prop:asympt lx}~\numitem{1} to see that
    \begin{align*}
        \mu t_{\mu/n}=2\left(\frac{\pi^2\mu^2n}{4}\right)^{1/3}+o(1)\,,
    \end{align*}
    so that the previous bound can transform into
    \begin{align*}
        \Prob{\frac{h(\bT_n)}{t_n}\geq1+\epsilon}\leq\big(1+o(1)\big)\exp\left(-\Big[3(1-\gamma)^{2/3}+2\gamma(1+\epsilon)-3\Big]\left(\frac{\pi^2\mu^2n}{4}\right)^{1/3}\right)\,.\numberthis \label{eq:LLN_UB}
    \end{align*}
    Using that $\epsilon$ is fixed and that, for $\gamma$ small enough, we have
    \begin{align*}
        \Big[3(1-\gamma)^{2/3}+2\gamma(1+\epsilon)-3\Big]\sim2\epsilon\gamma>0\,,
    \end{align*}
    we can find some $\gamma>0$ such that the right-hand side converges to $0$, thus proving the upper bound when $1/\sqrt{n}\ll\mu\ll n^{1/4}$.
    Using the same method when $\mu\asymp  n^{1/4}$ but with Theorem~\ref{thm:part-fun}~\numitem{2} and observing that the sum appearing in the formula of $Z_n$ is of order $1$, we again obtain the bound \eqref{eq:LLN_UB}.
    This leads to the same conclusion as for the case $\mu\ll n^{1/4}$.
    Finally, for the case $\mu\gg n^{1/4}$, observe again that the same method but applied with Theorem~\ref{thm:part-fun}~\numitem{3} works, thanks to the fact that $(\mu^4/n)^{1/3}\ll(\mu^2n)^{1/3}$.

    To conclude the proof of the theorem, note that a Chernoff bound applied to the lower tail of the height also leads to
    \begin{align*}
        \Prob{\frac{h(\bT_n)}{t_{\mu/n}}\leq1-\epsilon}\leq\big(A+o(1)\big)e^{2\gamma\mu}\exp\left(-\Big[3(1+\gamma)^{2/3}-2\gamma(1-\epsilon)-3\big]\left(\frac{\pi^2\mu^2n}{4}\right)^{1/3}\right)\,,
    \end{align*}
    where $A>0$ is some positive constant and the convergence to $0$ follows from the same argument as before and since $\mu\ll(\mu^2n)^{1/3}$.
\end{proof}

\subsection{Fluctuations and deviations from the mean}\label{sec:clt}

In this section, we extend the results from Section~\ref{sec:lln} to the second order behaviour and the large deviation of the height of $T_n$.
We heavily rely on the asymptotic behaviour of the partition function from Theorem~\ref{thm:part-fun}.
\subsubsection{Second order fluctuations}
In this section we provide the proofs of Theorems \ref{thm:inter clt} and \ref{thm:discrete clt}.
\begin{proof}[Proof of Theorem~\ref{thm:inter clt}]
Let us denote
\begin{align*}
    \sigma=\sqrt{\frac{1}{3}\left(\frac{2\pi^2n}{\mu^4}\right)^{1/3}}.
\end{align*}
For any $t\in\R$, consider the moment-generating function of the shifted and scaled height
    \begin{align*}
        f(t)=f_n(t)=\E\left[e^{t\frac{h(\bT_n)-(2\pi^2n/\mu)^{1/3}}{\sigma}}\right]&=e^{-t\frac{(2\pi^2n / \mu)^{1/3}}{\sigma}}\frac{Z_n^{\mu-t/\sigma}}{Z_n^\mu}\,.
    \end{align*}
    Now, applying Theorem~\ref{thm:part-fun}~\numitem{1} on the numerator and the denominator of the ratio of the partition functions above, we obtain that
    \begin{align*}
        f(t)\sim e^{-t/\sigma}\frac{e^{\mu-t/\sigma}-1}{e^\mu-1}\left(1-\frac{t}{\sigma\mu}\right)^{1/3}\exp\left(-t\frac{(2\pi^2n/\mu)^{1/3}}{\sigma}-3\left(\frac{\pi^2\mu^2n}{4}\right)^{1/3}\left[\left(1-\frac{t}{\sigma\mu}\right)^{2/3}-1\right]\right)\,.
    \end{align*}
    Observe first that $\sigma\asymp (n/\mu^4)^{1/6}$ so that $\sigma\gg1$ but also $\sigma\mu\asymp (\mu^2n)^{1/6}\gg1$, implying that all the terms outside of the final exponential are of order $1$.
    Combining this with a Taylor expansion of the term $(1-t/(\sigma\mu))^{2/3}$ leads to
    \begin{align*}
        f(t)&\sim\exp\left(-t\frac{(2\pi^2n/\mu)^{1/3}}{\sigma}-3\left(\frac{\pi^2\mu^2n}{4}\right)^{1/3}\left[-\frac{2}{3}\frac{t}{\sigma\mu}-\big(1+o(1)\big)\frac{1}{9}\frac{t^2}{\sigma^2\mu^2}\right]\right)\\
        &=\exp\left(\big(1+o(1)\big)\frac{1}{2}\left[\frac{1}{3}\left(\frac{2\pi^2n}{\mu^4}\right)^{1/3}\right]\frac{t^2}{\sigma^2}\right)\,.
    \end{align*}
    To conclude the proof, simply note that $\sigma$ is exactly defined so that $f(t)\sim e^{t^2/2}$, as required for the convergence in distribution to the standard normal.
\end{proof}

\begin{proof}[Proof of Theorem~\ref{thm:discrete clt}]
   {Recall $\lambda_x(t)=xt+\log(1+\tan^2(\pi/t))$ and its unique minimizer $t_x$ from \eqref{eq:func_lambdax_minvalue}.} For this proof, we can directly use Theorem~\ref{thm:part-fun}~\numitem{2} with the moment generating function of $h(\bT_n)$ to obtain
    \begin{align*}
        \E\left[e^{sh(\bT_n)}\right]&=\frac{Z_n^{\mu-s}}{Z_n^\mu}\sim e^{-2s}\left(1-\frac{s}{\mu}\right)e^{n(\lambda_{\mu/n}(t_{\mu/n})-\lambda_{(\mu-s)/n}(t_{(\mu-s)/n}))}\frac{S_{(\mu-s)/n}}{S_{\mu/n}}\,,
    \end{align*}
    where {for any $x>0$},
    \begin{align*}
        S_x=\sum_{m\in\Z}\exp\left(-\left(\frac{3}{2}+o(1)\right)\left(\frac{n^3x^4}{2\pi^2}\right)^{1/3}(m-t_x)^2\right)\,,
    \end{align*}
    corresponding to the sum in the asymptotics from Theorem~\ref{thm:part-fun}~\numitem{2}.
    It is worth observing that, since we consider $x\in\{\mu/n,(\mu-s)/n\}$, for both value it holds that $x\sim\gamma n^{-3/4}$.
    Moreover, Proposition~\ref{prop:asympt lx}~\numitem{1} tells us that
    \begin{align*}
        t_{(\mu-s)/n}=\left(\frac{2\pi^2n}{\mu-s}\right)^{1/3}+o(1)=\left(\frac{2\pi^2n}{\mu}\right)^{1/3}\left(1+\frac{1}{3}\frac{s}{\mu}+o\left(\frac{1}{\mu}\right)\right)+o(1)=t_{\mu/n}+\left(\frac{2\pi^2}{\gamma^4}\right)^{1/3}\frac{s}{3}+o(1)\,,
    \end{align*}
    where we used that $\mu\sim\gamma n^{1/4}$, that these formulae work for arbitrary $s$, so also for $s=0$, 
    and that $t_{\mu/n}=(2\pi^2n/\mu)^{1/3}+o(1)$, which is itself a direct consequence of Proposition~\ref{prop:asympt lx}~\numitem{1}.
    In a similar manner, using Proposition~\ref{prop:asympt lx}~\numitem{2}, we obtain that
    \begin{align*}
        n\lambda_{(\mu-s)/n}(t_{(\mu-s)/n})&=3\left(\frac{\pi^2(\mu-s)^2n}{4}\right)^{1/3}+\frac{1}{6}\left(\frac{(\mu-s)^4}{(2\pi)^4n}\right)^{1/3}+o(1)\\
        &=3\left(\frac{\pi^2\mu^2n}{4}\right)^{1/3}\left[1-\frac{2}{3}\frac{s}{\mu}-\frac{1}{9}\frac{s^2}{\mu^2}+o\left(\frac{1}{\mu^2}\right)\right]+\frac{1}{6}\left(\frac{\gamma}{2\pi}\right)^{4/3}+o(1)\\
        &=n\lambda_{\mu/n}(t_{\mu/n})-s\left(\frac{2\pi^2n}{\mu}\right)^{1/3}-\left(\frac{\pi}{2\gamma^2}\right)^{2/3}\frac{s^2}{3}+o(1)\,,
    \end{align*}
    where we used again that the formula holds for $s=0$.
    Plugging the previous results in the original formula, along with the fact $s/\mu\ll1$, we obtain that
    \begin{align*}
        \E\left[e^{sh(\bT_n)}\right]&\sim\exp\left(-2s+s\left(\frac{2\pi^2n}{\mu}\right)^{1/3}+\left(\frac{\pi}{2\gamma^2}\right)^{2/3}\frac{s^2}{3}\right)\frac{S_{(\mu-s)/n}}{S_{\mu/n}}\,,
    \end{align*}
    where
    \begin{align*}
        S_{\mu/n}\sim\sum_{m\in\Z}\exp\left(-\frac{3}{2}\left(\frac{\gamma^4}{2\pi^2}\right)^{1/3}\left(m-\left(\frac{2\pi^2n}{\mu}\right)^{1/3}\right)^2\right)
    \end{align*}
    and
    \begin{align*}
        S_{(\mu-s)/n}\sim\sum_{m\in\Z}\exp\left(-\frac{3}{2}\left(\frac{\gamma^4}{2\pi^2}\right)^{1/3}\left(m-\left(\frac{2\pi^2n}{\mu}\right)^{1/3}-\left(\frac{2\pi^2}{\gamma^4}\right)^{1/3}\frac{s}{3}\right)^2\right)\,.
    \end{align*}
    To conclude this proof, observe that the last term within the exponential simplifies as
    \begin{align*}
        \frac{3}{2}\left(\frac{\gamma^4}{2\pi^2}\right)^{1/3}\left(m-\left(\frac{2\pi^2n}{\mu}\right)^{1/3}-\left(\frac{2\pi^2}{\gamma^4}\right)^{1/3}\frac{s}{3}\right)^2&=\frac{3}{2}\left(\frac{\gamma^4}{2\pi^2}\right)^{1/3}\left(m-\left(\frac{2\pi^2n}{\mu}\right)^{1/3}\right)^2\\
        &\hspace{0.5cm}-ms+s\left(\frac{2\pi^2n}{\mu}\right)^{1/3}+\left(\frac{\pi}{2\gamma^2}\right)^{2/3}\frac{s^2}{3}\,,
    \end{align*}
    where we see that the first term exactly corresponds to the term in the exponential of the sum of $S_{\mu/m}$ and the last two terms are also the two terms in the exponential of the moment generating function.
    Doing all the subsequent simplifications lead to
    \begin{align*}
        \E\left[e^{sh(\bT_n)}\right]&\sim\frac{\displaystyle\sum_{m\in Z}e^{(m-2)s}\exp\left(-\frac{3}{2}\left(\frac{\gamma^4}{2\pi^2}\right)^{1/3}\left(m-\left(\frac{2\pi^2n}{\mu}\right)^{1/3}\right)^2\right)}{\displaystyle\sum_{m\in\Z}\exp\left(-\frac{3}{2}\left(\frac{\gamma^4}{2\pi^2}\right)^{1/3}\left(m-\left(\frac{2\pi^2n}{\mu}\right)^{1/3}\right)^2\right)}\,.
    \end{align*}
    The claimed convergence now follows by shifting the sum $m\in\Z$ by $m_{\mu/n}=\lfloor t_{\mu/n}\rfloor$ and recalling that Proposition~\ref{prop:asympt lx}~\numitem{1} states that $t_{\mu/n}=(2\pi^2n/\mu)^{1/3}+o(1)$.
\end{proof}

\subsubsection{Large deviations}\label{sec:ldp}

In this section we provide the proof of Theorem \ref{thm:ldp}. Let us begin by stating the Gärtner-Ellis theorem \cite[Theorem 2.3.6]{dembo2009large}, as this is the key ingredient of our proof. We don't define the large deviations principle (LDP), or the notion of good rate functions here, and refer the reader for these standard concepts to \cite[Chapter 1]{dembo2009large}.
\begin{theorem}[Gärtner-Ellis theorem]\label{thm:gartner_ellis}
    Consider a sequence of random variables $(X_n)_{n \geq 1}$ and a sequence $(a_n)_{n \geq 1}$ with $a_n \gg 1$. Assume for every $\lambda \in \R$ the limit
    \begin{align*}
    \lim_{n \to \infty}\frac{1}{a_n}\log \Exp{\Big.\exp\big(\lambda a_nX_n\big)}= \Lambda(\lambda),  \numberthis \label{eq:lim_log_mom_gen_fn}
\end{align*}
    exists point-wise as a convex function of $\lambda$, possibly assuming an extended real value for some $\lambda$, with the origin lying inside the interior of the domain $\cD_{\Lambda}:=\{\lambda \in \R:\Lambda(\lambda)<\infty\}.$ Let $\Lambda$ further satisfy the following:
    \begin{itemize}
    \item Interior of $\cD_{\Lambda}$ is non empty.
    \item $\Lambda$ is differentiable at all points of the interior of $\cD_{\Lambda}$.
    \item $\Lambda(\cdot)$ is steep, namely, $\lim_{n \to \infty} |\Lambda'(\Lambda_n)|=\infty$ whenever $\lambda_n$ is a sequence of interior points of $\cD_{\Lambda}$ converging to some boundary point of the interior of $\cD_{\Lambda}$.
\end{itemize}
Then the sequence $(X_n)_{n \geq 1}$ satisfies a LDP with speed $a_n$ and with good rate function $\Lambda^*$, the \emph{Fenchel-Legendre} transform of $\Lambda$, defined as,
\begin{align*}
    \Lambda^*(x)=\sup_{\lambda \in \R}\Big\{\lambda x -\Lambda(\lambda)\Big\}\,.
\end{align*}
\end{theorem}

We are now ready to provide the proof of Theorem \ref{thm:ldp}.
\begin{proof}[Proof of Theorem~\ref{thm:ldp}]
    
    Using Proposition~\ref{prop:asympt lx}~\numitem{1}, 
    we see that $t_{\mu/n}$, the minimum of $\lambda_{\mu/n}$ as defined in \eqref{eq:func_lambdax_minvalue}, is of order $(2\pi^2n/\mu)^{1/3}$.
    Thus, it is enough to prove that $h(\mathbf{T}_n)/t_{\mu/n}$ satisfies the claimed large deviation principle with speed $(\mu^2 n)^{1/3}$.
    Thus, the main technical part of the proof reduces to proving that
    \begin{align*}
        \Lambda(\gamma)=\lim_{n\rightarrow\infty}\frac{1}{(\mu^2n)^{1/3}}\log\E\left[e^{(\mu^2n)^{1/3}\gamma h(\mathbf{T}_n)/t_{\mu/n}}\right]
    \end{align*}
    exists and has the desired properties of Theorem \ref{thm:gartner_ellis}.
    To do so, we start by characterizing the above limit when replacing $(\mu^2n)^{1/3}$ with $\mu t_{\mu/n}$, which has the same order, thanks to Proposition~\ref{prop:asympt lx}~\numitem{1}.

    Assume first that $1/\sqrt{n}\ll\mu\ll n^{1/4}$.
    Then, using the asymptotic behaviour of the partition function as provided in Theorem~\ref{thm:part-fun}~\numitem{1}, we see that, for any $\gamma<1$
    \begin{align*}
        \E\left[e^{\mu\gamma h(T_n)}\right]=\frac{Z_n^{\mu(1-\gamma)}}{Z_n^\mu}\sim\frac{e^{\mu(1-\gamma)}(e^{\mu(1-\gamma)}-1)}{e^\mu(e^\mu-1)}(1-\gamma)^{2/3}\exp\left(-3\left(\frac{\pi^2\mu^2n}{4}\right)^{1/3}\Big((1-\gamma)^{1/3}-1\Big)\right)\,.
    \end{align*}
    Observing that $\mu\ll(\mu^2n)^{1/3}$ here, it follows that
    \begin{align*}
        \frac{1}{\mu t_{\mu/n}}\log\E\left[e^{\mu\gamma h(\mathbf{T}_n)}\right]&=-\frac{3}{\mu t_{\mu/n}}\left(\frac{\pi^2\mu^2n}{4}\right)^{1/3}\Big((1-\gamma)^{2/3}-1\Big)+o(1)\\
        &=\frac{3}{2}\Big(1-(1-\gamma)^{2/3}\Big)+o(1)\,.
    \end{align*}
    Note that this formula also works when $\gamma=1$ by replacing the numerator with $C_{n-1}\asymp 4^n/n^{3/2}$.

    Assume now that $n^{1/4}\ll\mu\ll n$.
    Then, using the asymptotic behaviour of the partition function as provided in Theorem~\ref{thm:part-fun}~\numitem{2}, we see that, for any $\gamma<1$
    \begin{align*}
        \E\left[e^{\mu\gamma h(\mathbf{T}_n)}\right]=\frac{Z_n^{\mu(1-\gamma)}}{Z_n^\mu}\sim\frac{e^{2\mu(1-\gamma)}}{e^{2\mu}}(1-\gamma)e^{-n\lambda_{\mu(1-\gamma)/n}(t_{\mu(1-\gamma)/n})+n\lambda_{\mu/n}(t_{\mu/n})}R\,,
    \end{align*}
    where $R$ is a fraction between two sums, both of them of order $e^{O((\mu^4/n)^{1/3})}=e^{o(\mu^2n)^{1/3}}$.
    Furthermore, thanks to Proposition~\ref{prop:asympt lx}~\numitem{2}, we know that
    \begin{align*}
        -n\lambda_{\mu(1-\gamma)/n}(t_{\mu(1-\gamma)/n})+n\lambda_{\mu/n}(t_{\mu/n})&=-3\left(\frac{\pi^2\mu^2n}{4}\right)^{1/3}\Big((1-\gamma)^{2/3}-1\Big)+O\left(\left(\frac{\mu^4}{n}\right)^{1/3}\right)\,.
    \end{align*}
    Using again that $(\mu^4/n)^{1/3}\ll(\mu^2n)^{1/3}\asymp \mu t_{\mu/n}$, we obtain that
    \begin{align*}
        \frac{1}{\mu t_{\mu/n}}\log\E\left[e^{\mu\gamma h(T_n)}\right]&=\frac{3}{2}\Big(1-(1-\gamma)^{2/3}\Big)+o(1)\,,\numberthis \label{eq:def_Gamma}
    \end{align*}
    and can further check that this formula still holds when $\gamma=1$.

    Defining
    \begin{align*}
        \Gamma(\gamma)=\frac{3}{2}\left(1-(1-\gamma)^{2/3}\right),
    \end{align*}
    we note thus that the limit $\Gamma(\gamma)$ exists for any $\gamma\leq1$, is convex, differentiable on $(-\infty,1)$, and its derivative diverges to $\infty$ at $1$. It is worth observing that, since $\Gamma$ is convex, we can also deduce that $\Gamma(\gamma)=\infty$ whenever $\gamma>1$, without having to compute it (which would otherwise require an understanding of the partition function when $\mu$ is negative).    
    All of these observations imply that $\Gamma$ satisfies the properties in Theorem \ref{thm:gartner_ellis}. Thus, by Theorem \ref{thm:gartner_ellis}, it follows that $(h(\mathbf{T}_n)/t_{\mu/n})_{n \geq 1}$ satisfies an LDP with respect to the good rate function that is the Fenchel-Legendre transform of $\Gamma(\gamma)=(3/2)(1-(1-\gamma)^{2/3})$ and with speed $\mu t_{\mu/n}$.

Observe that defining $\Lambda(\gamma):=\lim_{n \to \infty}\frac{1}{(\mu^2n)^{1/3}}\log\E\left[e^{(\mu^2n)\gamma h(\bT_n)/t_{\mu/n}}\right]$, we have the relation
    \begin{align*}
        \Lambda(\gamma)
        &\sim\frac{(2\pi^2)^{1/3}}{\mu t_{\mu/n}}\log\E\left[\exp\left(\mu\frac{\gamma+o(1)}{(2\pi^2)^{1/3}}h(\bT_n)\right)\right]\sim(2\pi^2)^{1/3}\Gamma\left(\frac{\gamma}{(2\pi^2)^{1/3}}\right)\,. \numberthis \label{eq:rel_lambda_gamma}
    \end{align*}
    
    We now compute the rate function $\Lambda^*(\cdot)$. For any $x\in\R$, we have
    \begin{align*}
        x=\Lambda'(\gamma)~\Longleftrightarrow&~x=\Gamma'\left(\frac{\gamma}{(2\pi^2)^{1/3}}\right)=\left(1-\frac{\gamma}{(2\pi^2)^{1/3}}\right)^{-1/3}~\Longleftrightarrow~\gamma=(2\pi^2)^{1/3}\left(1-\frac{1}{x^3}\right)\,.
    \end{align*}
    This value of $\gamma$ is defined for any $x>0$ and is less than $(2\pi^2)^{1/3}$, thus $\gamma \in \cD_\Lambda$.
    In that case, we see that the supremum in the definition of the Fenchel-Legendre transform of $\Lambda$ is reached when $\gamma$ is the above value, leading to
    \begin{align*}
        \Lambda^*(x)&=x(2\pi^2)^{1/3}\left(1-\frac{1}{x^3}\right)-\Lambda\left((2\pi^2)^{1/3}\left(1-\frac{1}{x^3}\right)\right)\\
        &=(2\pi^2)^{1/3}\left[x-\frac{1}{x^2}-\frac{3}{2}\left(1-\left(\frac{1}{x^3}\right)^{2/3}\right)\right]\\
        &=(2\pi^2)^{1/3}\left[x+\frac{1}{2x^2}-\frac{3}{2}\right]\,.
    \end{align*}
    This agrees with the rate function appearing in the statement of Theorem \ref{thm:ldp} when $x>0$, while for $x\leq 0$, the supremum in the definition of Fenchel-Legendre transform is infinite, as follows from \eqref{eq:def_Gamma} and the relation \eqref{eq:rel_lambda_gamma} between $\Lambda$ and $\Gamma$, agreeing with the claim of Theorem \ref{thm:ldp}.
\end{proof}

\subsection{Bernoulli fluctuations}\label{sec:bern}

In this section we provide the proof Theorem~\ref{thm:high lln} by splitting the result into two propositions, covering the regimes $\mu\ll n$ and $\mu=\Omega (  n)$.
We then combine them together to prove Theorem~\ref{thm:high lln}. {Recall the function $\lambda_x$, its minimizer $t_x$ as defined in \eqref{eq:func_lambdax_minvalue}, and the definition of $\delta_n$ and $m_x=\max\{\floor{t_x},3\}$ from \eqref{eq:def_m_n_delta_n}. }

\begin{proposition}\label{prop:inter ber}
    Let $\mu=\mu_n$ be such that $n^{1/4}\ll\mu\ll n$, and assume that $\delta_n$ admits a limit $\delta$ (possibly infinite).
    Let $\mathbf{T}_n$ be a $\mu$-height-biased tree of size $n$.
    Then, $h(\bT_n)-m_{\mu/n}+2$ converges in distribution to a Bernoulli random variable with parameter
    \begin{align*}
        p_\delta=\frac{e^{-\delta}}{1+e^{-\delta}}\,.
    \end{align*}
\end{proposition}

\begin{proof}
    Fix some $\gamma\in\R$.
    Recall that Theorem~\ref{thm:part-fun}~\numitem{3} tells us that, in the regime considered here, we have
    \begin{align*}
        Z_n=\big(1+o(1)\big)4^ne^{2\mu}\sum_{m=\max\{\lfloor t_{\mu/n}\rfloor,3\}}^{\lceil t_{\mu/n}\rceil}\frac{\tan^2\left(\frac{\pi}{m}\right)}{m}e^{-n\lambda_{\mu/n}(m)}\,,
    \end{align*}
    Further using that $t_{\mu/n}$ and $t_{(\mu+\gamma)/n}$ diverge here thanks to
    Proposition~\ref{prop:asympt lx}~\numitem{1}, we can remove the max from the lower bound of the sum and identify $m_{\mu/n}$ with $\lfloor t_{\mu/n}\rfloor$.
    Finally, observing that $\lambda_{x+y}(m)=ym+\lambda_x(m)$, we see that
    \begin{align*}
        \E\left[e^{\gamma h(\bT_n)}\right]=\frac{Z_n^{\mu-\gamma}}{Z_n^\mu}\sim e^{-2\gamma}\frac{\sum_{m=\lfloor t_{(\mu-\gamma)/n}\rfloor}^{\lceil t_{(\mu-\gamma)/n}\rceil}\frac{\tan^2\left(\frac{\pi}{m}\right)}{m}e^{\gamma m-n\lambda_{\mu/n}(m)}}{\sum_{\lfloor t_{\mu/n}\rfloor}^{\lceil t_{\mu/n}\rceil}\frac{\tan^2\left(\frac{\pi}{m}\right)}{m}e^{-n\lambda_{\mu/n}(m)}}\,. \numberthis \label{eq:mulln_mom_gen_func}
    \end{align*}
{To prove the Proposition, it is enough to show that 
\begin{align*}
    \Exp{e^{\gamma (h(\bT_n)-m_{\mu/n}+2)}}\to 1-p_\delta+p_\delta e^\gamma. \numberthis\label{eq:limiting_bernoulli_mgf}
\end{align*}}

   {Note that the sequence $t_{\mu/n}-\lfloor t_{\mu/n} \rfloor$ is bounded in $[0,1]$, and thus admits subsequential limits. Below we separate our analysis into three possible cases of these limits, and exhibit that the result claimed by the Proposition is coherent across these cases. Further, to simplify notation, we hide the subsequence notation and assume the sequence itself has the assumed behavior.}

    \noindent \textbf{Case: $t_{\mu/n}-\lfloor t_{\mu/n} \rfloor \to 0$.} Using equation \eqref{eq:mulln_mom_gen_func}, if $t_{\mu/n}-\lfloor t_{\mu/n}\rfloor$ converges to $0$, Proposition~\ref{prop:asympt lx}~\numitem{3} tells us that
    \begin{align*}
        n\lambda_{\mu/n}(\lfloor t_{\mu/n}\rfloor\pm1)-n\lambda_{\mu/n}(\lfloor t_{\mu/n}\rfloor)\geq Cn\left(\frac{\mu}{n}\right)^{2/3}\frac{1}{\left(\frac{n}{\mu}\right)^{2/3}}=C\left(\frac{\mu^4}{n}\right)^{1/3}\gg1\,,
    \end{align*}
    with a similar inequality when applied to $\lfloor t_{(\mu-\gamma)/n}\rfloor$, implying that we can ignore the term different from $\lfloor t_{\mu/n}\rfloor$ in both sums.
    In that case, we directly see that
    \begin{align*}
        \E\left[e^{\gamma h(\bT_n)}\right]\sim e^{\gamma\left(\lfloor t_{\mu/n}\rfloor-2\right)}\,.
    \end{align*}
    To see that this implies \eqref{eq:limiting_bernoulli_mgf}, using $t_{\mu/n}$ diverges, it suffices to show that $\delta_n\rightarrow\infty$.
    But this is straightforward since, either $t_{\mu/n}$ is an integer, and then $\lfloor t_{\mu/n}\rfloor=\lceil t_{\mu/n}\rceil$ so that $\delta_n=\mu\gg1$, or $t_{\mu/n}$ is not an integer, and then $\lceil t_{\mu/n}\rceil=\lfloor t_{\mu/n}\rfloor+1$, leading to
    \begin{align*}
        \delta_n&=n\lambda_{\mu/n}(\lfloor t_{\mu/n}\rfloor+1)-n\lambda_{\mu/n}(\lfloor t_{\mu/n}\rfloor)\geq C\left(\frac{\mu^4}{n}\right)^{1/3}\gg1\,.
    \end{align*}
    
    \noindent \textbf{Case: $t_{\mu/n}-\lfloor t_{\mu/n} \rfloor \to 1$.} Note that $t_{\mu/n}-\lfloor t_{\mu/n}\rfloor$ converges to $1$ is equivalent to saying that $t_{\mu/n}-\lceil t_{\mu/n}\rceil$ converges to $0$, and similar arguments as in the previous case apply, leading to
    \begin{align*}
        \E\left[e^{\gamma h(\bT_n)}\right]\sim e^{\gamma\left(\lceil t_{\mu/n}\rceil-2\right)}\,.
    \end{align*}
    To see that this is desired result, observe again that if $t_{\mu/n}$ is an integer the result follows, and otherwise $-\delta_n\gg1$, leading to $\delta=-\infty$, and the last display again agrees with \eqref{eq:limiting_bernoulli_mgf}.

    \noindent \textbf{Case: $\lim_{n \to \infty} t_{\mu/n}-\lfloor t_{\mu/n} \rfloor \in (0,1)$.} Assume finally that $t_{\mu/n}-\lfloor t_{\mu/n}\rfloor$ converges to some number in $(0,1)$.
    In that regime, Proposition~\ref{prop:asympt lx}~\numitem{1} implies that
    \begin{align*}
        t_{(\mu-\gamma)/n}=\left(\frac{2\pi^2n}{\mu-\gamma}\right)^{1/3}+o(1)=\left(\frac{2\pi^2n}{\mu}\right)^{1/3}+O\left(\left(\frac{n}{\mu^4}\right)^{1/3}\right)+o(1)=t_{\mu/n}+o(1)\,,
    \end{align*}
    which means that $\lfloor t_{\mu/n}\rfloor=\lfloor t_{(\mu-\gamma)/n}\rfloor$ for $n$ large enough.
    It follows that
    \begin{align*}
        \E\left[e^{\gamma h(\bT_n)}\right]&\sim e^{-2\gamma}\frac{\sum_{m=\lfloor t_{\mu/n}\rfloor}^{\lceil t_{\mu/n}\rceil}\frac{\tan^2\left(\frac{\pi}{m}\right)}{m}e^{\gamma m-n\lambda_{\mu/n}(m)}}{\sum_{m=\lfloor t_{\mu/n}\rfloor}^{\lceil t_{\mu/n}\rceil}\frac{\tan^2\left(\frac{\pi}{m}\right)}{m}e^{-n\lambda_{\mu/n}(m)}}\\
        &\sim e^{\gamma\left(\lfloor t_{\mu/n}\rfloor-2\right)}\frac{1+e^{\gamma}e^{n\lambda_{\mu/n}(\lfloor t_{\mu/n}\rfloor)-n\lambda_{\mu/n}(\lfloor t_{\mu/n}\rfloor+1)}}{1+e^{n\lambda_{\mu/n}(\lfloor t_{\mu/n}\rfloor)-n\lambda_{\mu/n}(\lfloor t_{\mu/n}\rfloor+1)}}\,,
    \end{align*}
    where we used that $\lfloor t_n\rfloor\sim\lceil t_n\rceil\gg1$ to take the fraction out of the sum.
    But now, we observe that the definition of $\delta_n$ implies that
    \begin{align*}
        n\lambda_{\mu/n}(\lfloor t_{\mu/n}\rfloor)-n\lambda_{\mu/n}(\lfloor t_{\mu/n}\rfloor+1)&=-\delta_n\,,
    \end{align*}
    leading to
    \begin{align*}
        \E\left[e^{\gamma h(\bT_n)}\right]&\sim e^{\gamma\left(\lfloor t_{\mu/n}\rfloor-2\right)}\frac{1+e^{\gamma}e^{-\delta_n}}{1+e^{-\delta_n}}\,,
    \end{align*}
    implying \eqref{eq:limiting_bernoulli_mgf}.
\end{proof}

\begin{proposition}\label{prop:high lln}
    Let $\mu=\mu_n$ be such that $c_n=\mu/n$ converges to some $c\in(0,\infty]$ and assume that $\delta_n$ admits a limit $\delta$ (possibly infinite).
    Let $\mathbf{T}_n$ be a $\mu$-height-biased tree of size $n$.
    Then, $h(\bT_n)-m_{\mu/n}+2$ converges in distribution to a Bernoulli random variable with parameter
    \begin{align*}
        p_{c,\delta}=\frac{m_c\tan^2\left(\frac{\pi}{m_c+1}\right)e^{-\delta}}{(m_c+1)\tan^2\left(\frac{\pi}{m_c}\right)+m_c\tan^2\left(\frac{\pi}{m_c+1}\right)e^{-\delta}}\,.
    \end{align*}
\end{proposition}

\begin{proof}
    The proof is very similar to that of Proposition~\ref{prop:inter ber}, except that Proposition~\ref{prop:asympt lx} does not apply here.
    We start by applying Theorem~\ref{thm:part-fun} to see that
    \begin{align*}
        Z_n\sim4^ne^{2\mu}\sum_{m=m_{\mu/n}}^{\lceil t_{\mu/n}\rceil}\frac{\tan^2\left(\frac{\pi}{m}\right)}{m}e^{-n\lambda_{\mu/n}(m)}\,,
    \end{align*}
    so that, using again $\lambda_{x+y}(m)=ym+\lambda_x(m)$, we obtain
    \begin{align*}
        \E\left[e^{\gamma h(\bT_n)}\right]=\frac{Z_n^{\mu-\gamma}}{Z_n^\mu}\sim e^{-2\gamma}\frac{\sum_{m=m_{(\mu-\gamma)/n}}^{\lceil t_{(\mu-\gamma)/n}\rceil}\frac{\tan^2\left(\frac{\pi}{m}\right)}{m}e^{\gamma m-n\lambda_{\mu/n}(m)}}{\sum_{m=m_{\mu/n}}^{\lceil t_{\mu/n}\rceil}\frac{\tan^2\left(\frac{\pi}{m}\right)}{m}e^{-n\lambda_{\mu/n}(m)}}\,.
    \end{align*}
    Observe now that, whenever $t_{\mu/n}\leq 3$, then the two sums can be reduced to $m=3$ since the possible other term (corresponding to the case $\lceil t_{(\mu-\gamma)n}\rceil=3$, {implying} $t_{(\mu-\gamma)/n}\sim3$) is negligible.
    This leads to the desired result, since in that case $\delta_n=\mu$ and so $\delta=\infty$ with $p_{c,\delta}=0$.
    We assume from now on that $t_{\mu/n}>3$ so that $m_{\mu/n}=\lfloor t_{\mu/n}\rfloor$ and $m_{(\mu-\gamma)/n}=\lfloor t_{(\mu-\gamma)/n}\rfloor$ (since $x\mapsto t_x$ is decreasing).
    
    Assume here that $t_{\mu/n}$ converges to some integer $\Tilde{m}_n$, so that $t_{(\mu-\gamma)/n}$ does as well.
    Using that $\mu/n\sim c$ and the fact that $\lambda_c$ is strictly convex, the previous two sums can be reduced to $m=\Tilde{m}_n$  since the other terms decrease exponentially faster.
    In that case, we see that
    \begin{align*}
        \E\left[e^{\gamma h(\bT_n)}\right]=\frac{Z_n^{\mu-\gamma}}{Z_n^\mu}\sim e^{\gamma(\Tilde{m}_n-2)}\,.
    \end{align*}
    If $t_{\mu/n}$ is an integer, then $\delta_n=\mu$ and so $\delta=+\infty$.
    Otherwise, we see that
    \begin{align*}
        \delta_n=n\lambda_{\mu/n}(\lceil t_{\mu/n}\rceil)-n\lambda_{\mu/n}(\lfloor t_{\mu/n}\rfloor)\,,
    \end{align*}
    which is positive if $t_{\mu/n}$ converges to $\lfloor t_{\mu/n}\rfloor$ and negative if it converges to $\lceil t_{\mu/n}\rceil$.
    Moreover, using that $\mu/n\sim c$ and that $\lambda_c$ is strictly convex, we see that $|\delta_n|\gg1$.
    Thus, either $\delta=+\infty$, in which case $\Tilde{m}_n=\lfloor t_{\mu/n}\rfloor$ and so $h(\bT_n)\rightarrow\lfloor t_{\mu/n}\rfloor-2$ and $p_{c,\delta}=0$, or $\delta=-\infty$, in which case $\Tilde{m}_n=\lceil t_{\mu/n}\rceil$ and so $h(\bT_n)\rightarrow\lfloor t_{\mu/n}\rfloor-1$ and $p_{c,\delta}=1$.
    In all cases, the desired convergence is verified.

    We assume now that $t_{\mu/n}-\lfloor t_{\mu/n}\rfloor$ converges to some $r\in(0,1)$.
    Observe that, for $\gamma$ small enough and $n$ large enough, this implies that $\lfloor t_{(\mu-\gamma)/n}\rfloor=\lfloor t_{\mu/n}\rfloor=m_c$.
    Recalling that we assume that $t_{\mu/n}$ and thus also $t_{(\mu-\gamma)/n}$ are both larger than $3$, we see that
    \begin{align*}
        \E\left[e^{\gamma h(\bT_n)}\right]&\sim e^{-2\gamma}\frac{\sum_{m=m_c}^{m_c+1}\frac{\tan^2\left(\frac{\pi}{m}\right)}{m}e^{\gamma m-n\lambda_{\mu/n}(m)}}{\sum_{m=m_c}^{m_c+1}\frac{\tan^2\left(\frac{\pi}{m}\right)}{m}e^{-n\lambda_{\mu/n}(m)}}\\
        &=e^{\gamma(m_c-2)}\frac{(m_c+1)\tan^2\left(\frac{\pi}{m_c}\right)+m\tan^2\left(\frac{\pi}{m_c+1}\right)e^{\gamma+n(\lambda_{\mu/n}(m_c)-\lambda_{\mu/n}(m_c+1))}}{(m_c+1)\tan^2\left(\frac{\pi}{m_c}\right)+m\tan^2\left(\frac{\pi}{m_c+1}\right)e^{n(\lambda_{\mu/n}(m_c)-\lambda_{\mu/n}(m_c+1))}}\,.
    \end{align*}
    To conclude the proof, observe that
    \begin{align*}
        n\big(\lambda_{\mu/n}(m_c)-\lambda_{\mu/n}(m_c+1)\big)&=-\mu+n\log\left(\frac{1+\tan^2\left(\frac{\pi}{m_c}\right)}{1+\tan^2\left(\frac{\pi}{m_c+1}\right)}\right)=-\delta_n\,,
    \end{align*}
    leading to the desired limiting distribution.
\end{proof}

With our two main results stated and proven, we can now straightforwardly combine them to obtain Theorem~\ref{thm:high lln}.

\begin{proof}[Proof of Theorem~\ref{thm:high lln}]
    The proof of Theorem~\ref{thm:high lln} directly follows from Proposition~\ref{prop:inter ber} and~\ref{prop:high lln} by observing that the case $\mu\ll n$ is equivalent to the case $c=0$.
\end{proof}

\section{Width in the non-Brownian regime}\label{sec:inter_width_proof}

In this section, we provide the proof of Theorem \ref{thm:width}. We employ a multiscale decomposition technique using random excursions that encode our height biased trees. 

Since we are interested in bounding the width of a tree by using its contour representation, it is worth pointing out the relation between the width and the corresponding walk.
For a discrete contour $C=(C_t:0\leq t\leq 2n-1)$ (seen as a continuous function on $[0,2n-1]$), we let
\begin{align*}
    w(C)=\frac{1}{2}\max\Big\{\big|\{t:C_t=\ell\}\big|:\ell\in\R\Big\}\,.
\end{align*}
Note that for any discrete walk the maximum above is finite. By observing that every edge of the tree exactly corresponds to one up-step and one down-step of the contour, we see that $w(C)$ exactly counts the number of edges at a certain depth in the tree.
By further using that each such edge leads to a unique node at the given depth, we see that $w(C^{(T)})=w(T)$.
It is worth noting that this equality does not apply in the case where the tree is a single node, but this definition will not be applied to this case.
Finally, we observe that $w$ is naturally extended to bridges and walks and that, for any bridge $B^{(n)}_{x,y}$ going from $x$ to $y$ and obtained as the concatenation of $B^{(k)}_{x,z}$ and $B^{(n-k)}_{z,y}$ respectively going from $x$ to $z$ and from $z$ to $y$, we have
\begin{align*}
    w\big(B^{(n)}_{x,y}\big)\leq w\big(B^{(k)}_{x,z}\big)+w\big(B^{(n-k)}_{z,y}\big)\,.
\end{align*}
This inequality will prove to be useful in order to decompose our random excursion in a collection of bridges.

\subsection{Relating height-biased and uniform trees}

\begin{proposition}\label{prop:contour_random_narrow_exc_pos}
    Let $\mu=\mu_n$ such that $1/\sqrt{n}\ll\mu\ll n$.
    Then for any $\epsilon>0$, there exists $\delta>0$ such that, for any sequence of events $(A_n)_{n\geq1}$, we have
    \begin{align*}
        \limsup_{n\rightarrow\infty}\frac{1}{(\mu^2n)^{1/3}}\log\left(\frac{\Prob{\bT_n\in A_n,\left|\frac{h(\bT_n)}{(2\pi^2n/\mu)^{1/3}}-1\right|\leq\delta}}{\CProb{\T_n\in A_n}{\left|\frac{h(\T_n)}{(2\pi^2n/\mu)^{1/3}}-1\right|\leq\delta}}\right)\leq\epsilon\,.
    \end{align*}
\end{proposition}

\begin{proof}
    Start by considering a fix $\delta>0$ and denote $h_n=(2\pi^2n/\mu)^{1/3}$ for the asymptotic height in this regime, recall Theorem \ref{thm:inter lln}.
    By observing that $\bT_n$ and $\T_n$ have the same distribution once we condition on their height, the numerator of the proposition can be transformed into
    \begin{align*}
        \Prob{\bT_n\in A_n,\left|\frac{h(\bT_n)}{h_n}-1\right|\leq\delta}
        &=\sum_{|k-h_n|\leq\delta h_n}\CProb{\bT_n\in A_n}{h(\bT_n)=k}\Prob{\big.h(\bT_n)=k\big.}\\
        &=\sum_{|k-h_n|\leq\delta h_n}\CProb{\T_n\in A_n}{h(\T_n)=k}\Prob{\big.h(\bT_n)=k\big.}\,.
    \end{align*}
    Moreover, using again the uniform distribution of $\bT_n$ and $\T_n$ on trees of a given height, we have that
    \begin{align*}
        \Prob{\big.h(\bT_n)=k\big.}=\frac{e^{-\mu k}}{Z_n}C_{n-1}\Prob{\big.h(\T_n)=k\big.}\,.
    \end{align*}
    Using that the sum is taken over all $|k-h_n|\leq\delta h_n$, we see that $\mu k=\mu h_n+O(\delta(\mu^2n)^{1/3})$, leading to
    \begin{align*}
        \Prob{\bT_n\in A_n,\left|\frac{h(\bT_n)}{h_n}-1\right|\leq\delta}
        &=\sum_{|k-h_n|\leq\delta h_n}\CProb{\T_n\in A_n}{h(\T_n)=k}\frac{e^{-\mu k}}{Z_n}C_{n-1}\Prob{\big.h(\T_n)=k\big.}\\
        &=e^{-\mu h_n+O(\delta(\mu^2n)^{1/3})}\frac{C_{n-1}}{Z_n}\Prob{\T_n\in A_n,\left|\frac{h(\T_n)}{h_n}-1\right|\leq\delta}\\
        &=e^{-\mu h_n+O(\delta(\mu^2n)^{1/3})}\frac{C_{n-1}}{Z_n}\CProb{\T_n\in A_n}{\left|\frac{h(\T_n)}{h_n}-1\right|\leq\delta}\Prob{\left|\frac{h(\T_n)}{h_n}-1\right|\leq\delta}\,.
    \end{align*}

    To conclude the proof, we first observe that
    \begin{align*}
        C_{n-1}\Prob{\left|\frac{h(\T_n)}{h_n}-1\right|\leq\delta}=H_{n,\floor{h_n(1+\delta)}+1}-H_{n,\ceil{h_n(1-\delta)}}\sim\frac{4^n\pi^2}{(1+\delta)^3h_n^3}\frac{1}{\left(1+\tan^2\frac{\pi}{\floor{h_n(1+\delta)}+1}\right)^n}
    \end{align*}
    where $H_{n,m}$ is defined in~\eqref{eq:Hnm} and its asymptotic behaviour is given in Theorem~\ref{thm:ht-bdd tree}.
    Denoting $\tilde{h_n}=\floor{h_n(1+\delta)}+1$ and observing that $\mu h_n=\mu\tilde{h}_n+O(\delta(\mu^2n)^{1/3})$, it follows that
    \begin{align*}
        \frac{\Prob{\bT_n\in A_n,\left|\frac{h(\bT_n)}{h_n}-1\right|\leq\delta}}{\CProb{\T_n\in A_n}{\left|\frac{h(\T_n)}{h_n}-1\right|\leq\delta}}
        &\sim\frac{e^{-\mu\tilde{h}_n+O(\delta(\mu^2n)^{1/3})}}{Z_n}\frac{4^n\pi^2}{(1+\delta)^3h_n^3}\frac{1}{\left(1+\tan^2\frac{\pi}{\tilde{h_n}}\right)^n}\\
        &=\frac{e^{-n\lambda_{\mu/n}(\tilde{h}_n)+O(\delta(\mu^2n)^{1/3})}}{Z_n}\frac{4^n\pi^2}{(1+\delta)^3h_n^3}\,,
    \end{align*}
    where $\lambda_x$ is defined in~\eqref{eq:func_lambdax_minvalue}.
    By considering $\delta>0$ small enough, we observe that Proposition~\ref{prop:asympt lx}~\numitem{3} tells us here that
    \begin{align*}
        n\lambda_{\mu/n}(\tilde{h}_n)&=\lambda_{\mu/n}(t_{\mu/n})+O\left(\delta^2(\mu^2n)^{1/3}\right)
    \end{align*}
    Finally, the statement of the proposition follows from the different cases of the asymptotic of $Z_n$ from Theorem~\ref{thm:part-fun}.
    Indeed, in the case $1/\sqrt{n}\ll\mu\ll n^{1/4}$, we have that
    \begin{align*}
        \frac{\Prob{\bT_n\in A_n,\left|\frac{h(\bT_n)}{h_n}-1\right|\leq\delta}}{\CProb{\T_n\in A_n}{\left|\frac{h(\T_n)}{h_n}-1\right|\leq\delta}}
        &\sim e^{-n\lambda_{\mu/n}(t_{\mu/n})+O(\delta(\mu^2n)^{1/3})}\frac{4^n\pi^2}{(1+\delta)^3h_n^3}\frac{1}{4^ne^\mu(e^\mu-1)}\frac{2^{1/3}3^{1/2}}{\pi^{5/6}}\frac{n^{5/6}}{\mu^{1/3}}e^{3\left(\frac{\pi^2\mu^2n}{4}\right)^{1/3}}\\
        &=e^{O(\delta(\mu^2n)^{1/3})}\,,
    \end{align*}
    and the same result applies to $\mu\asymp n^{1/4}$ and to $n^{1/4}\ll\mu\ll n$.
    This concludes the proof as, for any $\epsilon>0$, we are then able to find $\delta>0$ small enough such that $O(\delta(\mu^2n)^{1/3})<\epsilon(\mu^2n)^{1/3}$.
\end{proof}

\subsection{Bounding the width of a bridge}
\begin{proposition}[Subgaussian width tails for bridges] \label{prop:bridge_subg_tail}
    Let $B^{(n)}_x=(B^{(n)}_{x;t}:0\leq t \leq n)$ be a uniform bridge of length $n$ ending at $x$, that is a uniform element from $\cB^{(n)}_x$.
    Then, for any $A>0$, there exists $R_1,R_2>0$ such that, for any $|x|\leq A\sqrt{n}$ and $T\geq0$, we have
    \begin{align*}
        \Prob{\frac{w\big(B^{(n)}_x\big)}{\sqrt{n}}>T}\leq R_1e^{-R_2T^2}\,.
    \end{align*}
\end{proposition}

\begin{proof}
We first observe that the reflection $\phi_x$ defined in Lemma~\ref{lem:cut_bij} can at most double the width of the walk, meaning that
\begin{align*}
    \frac{1}{2}w\big(B^{(n)}_x\big)\leq w\big(\phi_x\big(B^{(n)}_x\big)\big)\leq 2w\big(B^{(n)}_x\big) \, .
\end{align*}
Moreover, using that $\phi_x$ is a bijection and thus pushes forward uniform measures to uniform measures, and letting $e=-\mathbbm{1}_{\textrm{$n$ is odd}}$ for the end-point of the new bridge, we have 
\begin{align*}
    \Prob{\frac{w\big(B^{(n)}_x\big)}{\sqrt{n}}>T} &\leq \CProb{\frac{w\big(B^{(n)}_e\big)}{\sqrt{n}}>\frac{T}{2}}{\exists t,B^{(n)}_{e;t}=\floor{x/2}} 
    \leq \frac{\Prob{w\big(B^{(n)}_e\big)>\sqrt{n}T/2}}{\Prob{\exists t,B^{(n)}_{e;t}=\floor{x/2}}} \, .
\end{align*}
Using the distributional convergence of $B^{(n)}_e$ towards the Brownian bridge along with the fact that $|\floor{x/2}|\leq A\sqrt{n}$, we see that the denominator asymptotically satisfies
\begin{align*}
    \Prob{\exists t,B^{(n)}_{e;t}=\floor{x/2}}
    &\geq\Prob{\max_{0\leq t\leq n}\frac{B^{(n)}_{e;t}}{\sqrt{n}}>A,\min_{0\leq t\leq n}\frac{B^{(n)}_{e;t}}{\sqrt{n}}<-A}\\
    &=\big(1+o(1)\big)\Prob{\max_{0\leq t\leq 1}B_t>A,\min_{0\leq t\leq 1}B_t<-A}\,,
\end{align*}
where $(B_t:0\leq t\leq 1)$ is a standard Brownian bridge from $0$ to $0$. Furthermore, the right-hand side corresponds to a strictly positive function $K(A)>0$, see for instance \cite[Formula 1.15.8, Pg 180]{borodin2012handbook}.
To conclude the proof, first observe that if $n$ is even, we can extend $B^{(n)}_0$ to $B^{(n+1)}_{-1}$ and add at most $+1$ to the width.
Then, by observing that the canonical map $\phi$ from Lemma~\ref{lem:cycle} does not modify the width, the numerator
\begin{align*}
    \Prob{\Big.w\big(B^{(n)}_e\big)>\sqrt{n}T/2}\leq\Prob{\Big.w\big(E^{(n)}\big)>\sqrt{n}T/2-1}\,,
\end{align*}
where $E^{(n)}$ is a uniform excursion of size either $n$ or $n+1$ (depending on the parity of $n$).
Finally, using that the width of an excursion is equal to the width of the corresponding tree, and that the distribution of the corresponding tree is a critical Bienamyé tree with $\mathrm{Geo}(1/2)$ offspring distribution, conditioned to have size $\ceil{n/2}$, we can use~\cite[Theorem 1.1]{LAB_LD_SJ_13} which states the existence of $c_1, C_1$ positive constants independent of $T$ and $n$ such that
\begin{align*}
    \Prob{\Big.w\big(E^{(n)}\big)>\sqrt{n}T/2-1}\leq C_1e^{-c_1(\sqrt{n}T/2-1)^2/n}\,.
\end{align*}
Combining all previous inequalities lead to the existence of $R_1,R_2>0$ such that, for $n$ large enough, we have
\begin{align*}
    \Prob{\frac{w\big(B^{(n)}_x\big)}{\sqrt{n}}>T} &\leq \big(1+o(1)\big)\frac{C_1e^{-c_1(\sqrt{n}T/2-1)^2/n}}{K(A)}\leq R_1e^{-R_2T^2}\,.
\end{align*}
This can then be extended to the desired result by choosing $R_1$ and $R_2$ such that the inequality trivially holds when $n$ is small.
\end{proof}

\begin{corollary}[Exponential moments of the width of uniform bridges]\label{cor:exp_mom_local_time}
    We denote by $B^{(n)}_{x,y}$ a uniform bridge of length $n$ from $x$ to $y$, that is a uniform element from $\cB^{(n)}_{x,y}$.
    Then, for any $s\in\R$ and $A>0$, there exists $M>0$ such that, for all $n\geq1$ and all $0\leq x,y\leq A\sqrt{n}$, we have
    \begin{align*}
        \Exp{\exp\left(s\frac{w\big(B^{(n)}_{x,y}\big)}{\sqrt{n}}\right)}\leq M\,.
    \end{align*}
\end{corollary}

\begin{proof}
    We start by observing that $B^{(n)}_{x,y}$ can be mapped to $B^{(n)}_{y-x}$ a bridge from $0$ to $y-x$, by simply translating the walk.
    Moreover, shifting a walk does not affect its width, so we have
    \begin{align*}
        w\big(B^{(n)}_{x,y}\big)\overset{d}{=}w\big(B^{(n)}_{y-x}\big)\,.
    \end{align*}
    Further using that $0<w(B^{(n)}_{y-x})\leq n$, we see that the desired expected value re-writes as
    \begin{align*}
        \Exp{\exp\left(s\frac{w\big(B^{(n)}_{x,y}\big)}{\sqrt{n}}\right)}
        &=\Exp{\exp\left(s\frac{w\big(B^{(n)}_{y-x}\big)}{\sqrt{n}}\right)}
        =\sum_{k=1}^n\Exp{\exp\left(s\frac{w\big(B^{(n)}_{y-x}\big)}{\sqrt{n}}\right)\mathbbm{1}_{\{\sqrt{k-1}<w(B^{(n)}_{y-x})/\sqrt{n}\leq\sqrt{k}\}}}\,,
    \end{align*}
    
    Now, thanks to Proposition~\ref{prop:bridge_subg_tail} and since $|y-x|\leq A\sqrt{n}$ whenever $0\leq x,y\leq A\sqrt{n}$, we have that
    \begin{align*}
        \Exp{\exp\left(s\frac{w\big(B^{(n)}_{x,y}\big)}{\sqrt{n}}\right)\mathbbm{1}_{\{\sqrt{k-1}<w(B^{(n)}_{y-x})/\sqrt{n}\leq\sqrt{k}\}}}
        &\leq e^{s\sqrt{k}}\Prob{\frac{w\big(B^{(n)}_{y-x}\big)}{\sqrt{n}}>\sqrt{k-1}}
        \leq R_1e^{s\sqrt{k}-R_2(k-1)}\,.
    \end{align*}
    Since the sum of the right-hand side above is finite even if we sum over all $k \in \N$, the corollary follows.
\end{proof}

\subsection{Width in the non-Brownian regime}

\begin{figure}[H]
    \centering
    \includegraphics[width=0.98\textwidth]{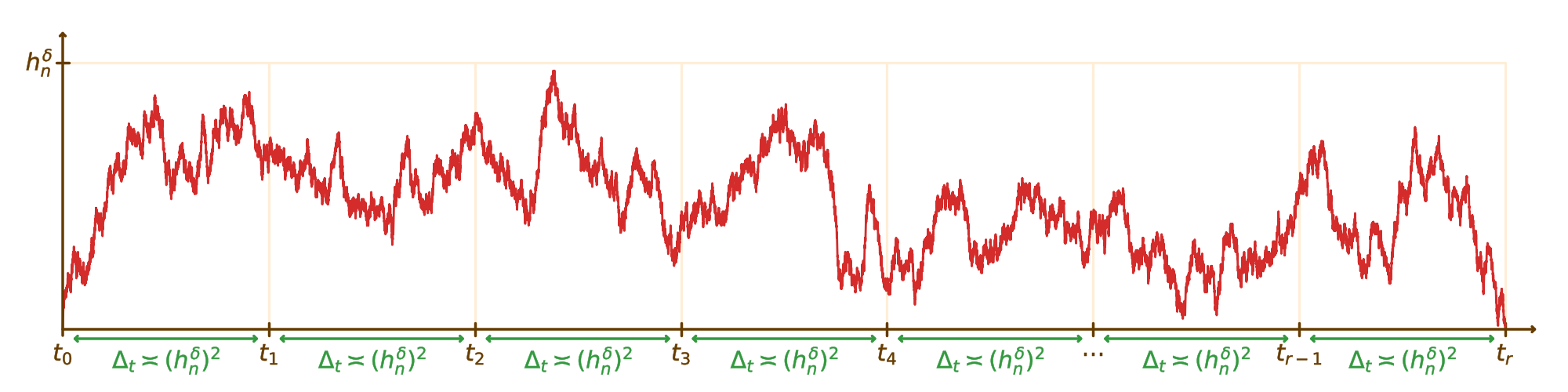}
    \caption{A representation of the multiscale decomposition used for the proof of Theorem~\ref{thm:width}. We bound the height with $h_n^\delta=(1+\delta)(2\pi^2n/\mu)^{1/3}$ and split the contour of a uniform tree into $r\sim(\mu^2n)^{1/3}$ bridges of length $\Delta_t\sim2(n/\mu)^{2/3}$. This way the length of each bridge is of the same order as the square of its fluctuations, thus are objects from the \emph{Brownian world}, admitting good enough width tails for our purposes, thanks to Proposition~\ref{prop:bridge_subg_tail}.}
    \label{fig:broscale}
\end{figure}

To conclude this section, we now prove Theorem~\ref{thm:width} stating that the width of a height-biased tree when $1/\sqrt{n}\ll\mu\ll n$ is of order $(\mu n^2)^{1/3}$, although without providing the exact constant.
The proof relies on splitting the contour of a tree with bounded height into standard bridges whose height is the square of their length.
A representation of this multiscale idea can be found in Figure~\ref{fig:broscale}.

\begin{proof}[Proof of Theorem~\ref{thm:width}]
    To start the proof, first note that, since the width counts the number of nodes at a certain height in the tree, by summing over all possible heights we directly have the bound
    \begin{align*}
        w(\bT_n)\cdot h(\bT_n)\geq n\,.
    \end{align*}
    Using the result of Theorem~\ref{thm:inter lln} which states that $h(\bT_n)$ is asymptotically of order $(n/\mu)^{1/3}$, we directly see that $w(\bT_n)$ is at least of order $(\mu n^2)^{1/3}$ as desired.
    This proves the lower bound of the theorem and we now focus on proving that $w(\bT_n)$ is not of higher order than $(\mu n^2)^{1/3}$.
    Following the proof of Proposition~\ref{prop:contour_random_narrow_exc_pos}, we let $h_n=(2\pi^2n/\mu)^{1/3}$ be the asymptotic height of $\bT_n$.

    Fix for now some $K>0$ and $\epsilon>0$.
    We start by applying Proposition~\ref{prop:contour_random_narrow_exc_pos} to the event $A_n=\{w(T)>K(\mu n^2)^{}\}$ to find $\delta>0$ such that
    \begin{align*}
        \limsup_{n\rightarrow\infty}\frac{1}{(\mu^2n)^{1/3}}\log\left(\frac{\Prob{w(\bT_n)>K(\mu n^2)^{1/3},\left|\frac{h(\bT_n)}{h_n}-1\right|\leq\delta}}{\CProb{w(\T_n)>K(\mu n^2)^{1/3}}{\left|\frac{h(\T_n)}{h_n}-1\right|\leq\delta}}\right)\leq\epsilon\,.
    \end{align*}
    By further using Theorem~\ref{thm:inter lln}, we know that $h(\bT_n)$ converges in probability to $h_n$.
    Combining this with the previous result leads to
    \begin{align}\label{eq:bound width}
        \Prob{\Big.w(\bT_n)>K(\mu n^2)^{1/3}}\leq e^{O(\epsilon(\mu^2n)^{1/3})}\CProb{w(\T_n)>K(\mu n^2)^{1/3}}{\left|\frac{h(\T_n)}{h_n}-1\right|\leq\delta}+o(1)\,.
    \end{align}
    Proving the theorem now boils down to proving an upper bound on the width of a uniform tree.

    Given a uniform tree $\T_n$, recall we let $C^{(\T_n)}=(C^{(\T_n)}_t:0\leq t\leq 2n-1)$ to be its contour walk.
    We are going to split this walk into a sequence of bridges of length $2(n/\mu)^{2/3}$ so that, once conditioned to be within the range $[0,h_n)$, their length is exactly this range squared. 
    We refer to Figure~\ref{fig:broscale} for a representation of this idea.

    Let $r=r_n=\floor{(\mu^2n)^{1/3}}$ and set $t_0,\ldots,t_r$ to evenly spaced time points from the interval $2n-1$, meaning that
    \begin{align*}
        t_i=\floor{\frac{i(2n-1)}{r}}\,.
    \end{align*}
    We observe here that $r\sim(\mu^2n)^{1/3}$ and that $\Delta_t=\max\{t_i-t_{i-1}\}\sim2(n/\mu)^{2/3}$.
    Finally, we denote by $B_i=(C^{(\T_n)}_{t+t_{i-1}}:0\leq t\leq t_i-t_{i-1})$ the walk corresponding to the portion of $C^{(n)}$ between time $t_{i-1}$ and time $t_i$ (rescaled in time to start at $0$).
    We now make two important observations.
    \begin{enumerate}
        \item Since $C^{(\T_n)}$ is obtained by concatenating $B_1,\ldots,B_r$, its width is bounded by the sum of their widths:
        \begin{align*}
            w\big(C^{(\T_n)}\big)\leq\sum_{i=1}^rw(B_i)\,.
        \end{align*}
        \item By the Markov property $B_1,B_2,\ldots,B_r$ are independent and uniform bridges conditionally given $C^{(\T_n)}_{t_i}$ for $0\leq i\leq r$.
    \end{enumerate}
    To use these observations, start by applying a standard exponential Markov's inequality to see that, for any $\rho>0$, we have
    \begin{align*}
        \CProb{w(\T_n)>K(\mu n^2)^{1/3}}{\left|\frac{h(\T_n)}{h_n}-1\right|\leq\delta}
        &\leq e^{-\rho K(\mu n^2)^{1/3}}\Exp{e^{\rho w(\T_n)}~\bigg|~\left|\frac{h(\T_n)}{h_n}-1\right|\leq\delta}\,.
    \end{align*}
    We then combine the first observation and a nested expectation to obtain that
    \begin{align*}
        &\hspace{-0.5cm}\CProb{w(\T_n)>K(\mu n^2)^{1/3}}{\left|\frac{h(\T_n)}{h_n}-1\right|\leq\delta}\\
        &\leq e^{-\rho K(\mu n^2)^{1/3}}\Exp{\Exp{e^{\rho w(\T_n)}~\bigg|~\left|\frac{h(\T_n)}{h_n}-1\right|\leq\delta,C^{(\T_n)}_{t_0},C^{(\T_n)}_{t_1},\ldots,C^{(\T_n)}_{t_r}}}\\
        &\leq e^{-\rho K(\mu n^2)^{1/3}}\Exp{\Exp{e^{\rho\left[w(B_1)+\ldots+w(B_r)\right]}~\bigg|~\left|\frac{h(\T_n)}{h_n}-1\right|\leq\delta,C^{(\T_n)}_{t_0},C^{(\T_n)}_{t_1},\ldots,C^{(\T_n)}_{t_r}}}\,.
    \end{align*}
    The condition that $h(\T_n)$ is bounded can be reduced to simply having $0\leq C^{\T_n}_{t_i}\leq(1+\delta)h_n$ by dividing with the probability that a walk of length $2n-1$ is a contour whose height is of order $h_n$, which itself is of order $e^{O(\log n)+O((\mu^2n)^{1/3})}$.
    Finally, the second observation states that, under the above conditioning, all bridges are independent from each other and their endpoints both belong to the interval $[0,(1+\delta)h_n]$.
    Further using that the length of the bridges is bounded by $\Delta_t\sim2(n/\mu)^{2/3}$ and that extending a walk can only increase its width, the previous bound can be transformed into
    \begin{align*}
        \CProb{w(\T_n)>K(\mu n^2)^{1/3}}{\left|\frac{h(\T_n)}{h_n}-1\right|\leq\delta}
        &\leq e^{-\rho K(\mu n^2)^{1/3}+O(\log n)+O((\mu^2n)^{1/3})}\left(\max_{0\leq x,y\leq A\sqrt{\Delta_t}}\Exp{e^{\rho w\left(B^{(\Delta_t)}_{x,y}\right)}}\right)^r\,,
    \end{align*}
    where $A>0$ is some constant such that $(1+\delta)h_n\leq A\sqrt{\Delta_t}$ and $B^{(\Delta_t)}_{x,y}$ is a uniform random bridge of length $\Delta_t$ going from $x$ to $y$.

    For the final part of this proof, start by replacing $\rho$ with $1/\sqrt{\Delta_t}$ and apply Corollary~\ref{cor:exp_mom_local_time} with the same $A$ and $s=1$ to find some $M>0$ such that, for any $\Delta_t>0$
    \begin{align*}
        \max_{0\leq x,y\leq A\sqrt{\Delta_t}}\Exp{e^{w\left(B^{(\Delta_t)}_{x,y}\right)/\sqrt{\Delta_t}}}\leq M\,.
    \end{align*}
    Finally, by recalling that $\mu\gg(\log n)^{3/2}/\sqrt{n}$, that $r\sim(\mu^2n)^{1/3}$, and that $\Delta_t\sim2(n/\mu)^{2/3}$, combining the previous two bounds lead to
    \begin{align*}
        \CProb{w(\T_n)>K(\mu n^2)^{1/3}}{\left|\frac{h(\T_n)}{h_n}-1\right|\leq\delta}&\leq\exp\left(-\big(1+o(1)\big)\left[\frac{K}{\sqrt{2}}-\log M+O(1)\right](\mu^2n)^{1/3}\right)\,.
    \end{align*}
    Plugging this back into~\eqref{eq:bound width} finally leads to
    \begin{align*}
        \Prob{\Big.w(\bT_n)>K(\mu n^2)^{1/3}}\leq\exp\left(-\big(1+o(1)\big)\left[\frac{K}{\sqrt{2}}-\log M+O(\epsilon)+O(1)\right](\mu^2n)^{1/3}\right)+o(1)
    \end{align*}
    Observe that the definitions of $M$, $O(\epsilon)$, and $O(1)$ only depend on $\epsilon$, $\delta$, and $A$ and are thus independent of $K$.
    This means that we can take $K$ large enough such that the exponential converges to $0$, thus proving the desired upper bound and concluding the proof.
\end{proof}

\section{Local results around the root: local convergence and condensation}\label{sec:loc_results}
In this section we prove Theorems \ref{thm:loclim} and \ref{thm:root deg}, completing the picture of local behavior around the root.

\subsection{Local limits}\label{sec:local_conv_prelims}

Before moving on to the results, we briefly discuss the precise topology that the convergence in Theorem \ref{thm:loclim} takes place in. Even though the following definitions are standard, we include it for the sake of completeness. The notion of local convergence of sparse graphs goes back to works of Benjamini and Schramm \cite{BSc01}, and Aldous and Steele \cite{AS04}. In the context of rooted planar trees, we use the notion of local convergence as presented by Abraham and Delmas in \cite{abraham2015introduction}, which respects the ordering of the children of every vertex.

Let $T_1, T_2$ be two rooted plane trees (possibly infinite). We equip the set of all plane trees $\mathbbm{t}$ with the metric

\begin{align*}
    d_{\mathrm{loc}}(T_1,T_2)= \inf \left\{\frac{1}{1+R} : R \in \mathbb{N}, B_R(T_1)=B_R(T_2) \right\}
\end{align*} where $B_R(T)$ denotes the subtree of $T$ spanned by the vertices at distance at most $R$ from the root. 

\begin{definition}[Balls around a rooted tree]
  For $T$ a rooted tree and $a>0$, define by $\mathcal{B}_{a}(T)$ the open ball of radius $a$ around $T$ in the local topology, i.e.,
  \begin{align*}
      \mathcal{B}_{a}(T):=\{T:d_{\mathrm{loc}}(T,T') < a \}.
  \end{align*}
  In other words, for $r\in \mathbb{N}^*$, $\mathcal{B}_{1/r}(T)$ consists of those rooted trees that look identical with $T$ up to generation $r$.
\end{definition}

\begin{remark} \label{rem:loclim}
The set of finite trees $\mathbbm{t}_{f}$ is a countable dense subset of $\mathbbm{t} = \mathbbm{t}_f \cup \mathbbm{t}_\infty$, where $\mathbbm{t}_\infty$ stands for the set of infinite trees. It is also an easy exercise to check that the balls $\mathcal{B}_{1/r}$ are both closed and open, and two balls are either disjoint or one is contained in the other.
\end{remark}

It is a well known result of Kesten \cite{kesten1986subdiffusive} that the local limit of the uniform tree $\T_n$ is the the so-called Kesten's tree (KT) associated to the critical $\mathrm{Geo}(1/2)$ offspring distribution with mass function 
$$p(n)=\left( \frac 1 2 \right)^{n+1} \, ,n=0,1,2,\dots$$
which is a multi-type Biénaymé tree with two types of individuals: normal and special. The special individuals reproduce according to the size-biased distribution $p^*(n)=np(n), n=0,1,\dots,$ whereas the normal ones reproduce according to $p(n)$. The root is special and reproduces according to $p^*(n)=np(n)$ where one of the children, chosen uniformly, is special and the others are normal. This is an infinite tree with a single spine, i.e. with a unique infinite path emanating from the root. When we write Kesten's tree or simply KT, we always refer to the one associated to the critical geometric offspring distribution {as described above}, since this is the only case we encounter in this work. 

\medskip

We now move on to prove \cref{thm:loclim}, stating that the local limit of the tree $\mathbf{T}_n$ when $0 \leq \mu \ll 1$ is also KT, i.e., the height bias is not strong enough to change the local behavior. We know from \cite{durhuus2023trees} that when $\mu$ is a positive real constant, the local limit is the Poisson tree with parameter $\mu$ of \cite{abraham2020very}. This is an infinite tree which is not an instance of Kesten's tree, but it still has two types, with a backbone without leaves corresponding to
individuals having an infinite progeny on which are
grafted independent critical (finite) Bienaymé trees. Since $\mu = \alpha$ with $\alpha>0$ is already covered in \cite{durhuus2023trees}, and the $\mu \to \alpha$ case can be treated with the same techniques found therein up to minor modifications, we do not give further details on the explicit properties of this limiting tree and refer the reader to \cite{abraham2020very,durhuus2023trees}.

For the $0 \leq \mu \ll 1$ case,
we will need the following lemma for the proof of \cref{thm:loclim}. Below, for notational convenience, we write for any event $A$ of trees,
\begin{align*}
    \Prob{\bT_n \in A}=\nu_n^\mu(A)\,.
\end{align*}

\begin{lemma} \label{lem:ballvol}
Assume that ~ $0 \leq \mu \ll 1$ and $T_0$ a finite tree of height $r$ with $K$ vertices at level $r$. For any $M \in \mathbb{N}$, it holds that
\begin{equation} \label{ballvol}
\nu_n ^\mu \big(\mathcal{B}_{1/r} (T_0) \big) \geq\big(1 + o(1)\big)K \frac{e^{-\mu r}}{4^{|T_0| -K}} \left( \sum\limits_{S=1} ^M  C_{S-1} 4^{-S} \right)^{K-1} \, ,
\end{equation}
where $\mathcal{B}_{1/r}(T_0)$ is the open ball of radius $1/r$ around the tree $T_0$ in the local topology.
\end{lemma}

\begin{proof}

\noindent \textbf{Case $1/\sqrt{n}\ll \mu \ll 1$.} Since $\mu \to 0$ and $\mu^2 n \to \infty$ for the window in question, we can assume without loss of generality that $\mu=\mu_n$ and $\mu^2 n$ are respectively decreasing and increasing functions of $n$.
Note that any tree  $ T \in \mathcal{B}_{1/r} (T_0)$ can be obtained by grafting $K$ trees on the top $K$ vertices of $T_0$, i.e., the \emph{leaves} of $T_0$ at maximal height. The subtrees rooted at the top vertices of $T_0$ are referred to as \emph{branches} in the rest of the proof.
Let $M$ be some fixed large constant and consider the subset of $\mathcal{B}_{1/r} (T_0)$ where all branches except one have size $< M$. Assuming that $n > KM + |T_0|$, {exactly} one branch has size $\geq M$ which we call the large one. We further restrict to the case where the height of the subtree that is grafted on the leaf corresponding to the one with the largest subtree, is at least $M$. In other words, the largest branch is also the highest with height at least $M$, and the total height of the tree $T$ thus formed after the grafting procedure is at least $r+M$.
 
Let $\Omega(T_0, r)$ be the subset of $\mathcal{B}_{1/r} (T_0)$ consisting of trees of size $n$ satisfying the properties above and whose large branch is precisely the first one, where we hide from notation the $n$-dependence of $\Omega$. Then, {by using the symmetry between all branches and only considering the case where the first branch is the largest one, we have}
$$
\nu^\mu _n \big(\mathcal{B}_{1/r} (T_0)\big) \geq K \cdot \nu^\mu _n \big(\Omega(T_0,r)\big) \,.
$$
We proceed to estimate the right hand side. Denoting the sizes of the subbranches grafted on the top $K$ vertices by $N_i, ~i\in \{1, \dots , K\},$ we have
\begin{align*}
\nu^\mu _n \big(\Omega(T_0,r)\big) & \geq  \frac{1}{Z_n^\mu}   \sum \limits_{\substack{n_2,\dots, n_K <M\\ n_1 = n- n_2 - \dots - n_K}} C_{n_2-1} \dots C_{n_K-1} \sum_{m=M} ^n e^{-\mu (m+r)} (H_{n_1,m+1}-H_{n_1,m}) \, ,
\end{align*}
where $H_{n,m}$ is the number of trees of height $<m$, as defined in~\eqref{eq:Hnm}.

Recall that $1/\sqrt{n} \ll \mu \ll 1$. For this window, we employ \cref{lem:z and w} and \cref{prop:W main}, where we choose $\epsilon_n$ to be $1/(\mu^2n)^{1/6 -\delta}$ for any fixed $0< \delta < 1/6$ to estimate the second sum on the right hand side. Using $n_1 < n$, $\mu_{n_1}> \mu_n$ and $\mu^2_{n_1}n_1 < \mu^2 _n n$ we have, for $n$ large,
\begin{align} \label{eq: low_bound}
\sum_{m=M} ^n e^{-\mu_n m} \big(H_{n_1,m+1}-H_{n_1,m}\big) &\geq \sum_{m=M} ^n e^{-\mu_{n_1} m} \big(H_{n_1,m+1}-H_{n_1,m}\big)\nonumber \\
&\geq\big(1+ o(1)\big)  \frac{4^{n_1}\pi ^{5/6}}{2^{1/3} 3^{1/2}} \frac{\mu_n}{n} \exp{\left( -3\left( \frac{\pi^2 \mu_n ^2 n}{4} \right)^{1/3} \right)}\,,
\end{align}
which together with \cref{thm:part-fun} leads to
\bb \label{eq: ball_vol}
\nu^\mu _n \big(\Omega(T_0,r)\big) \geq \big(1 + o(1)\big) \frac{e^{-\mu r }}{4^{|T_0|-K}} \left( \sum_{s=1} ^M C_{s-1} 4^{-s} \right)^{K-1}\,,
\ee
which concludes the proof for the regime in question once multiplied by $K$.

\smallskip

\noindent \textbf{Case $\mu \sqrt{n} \to \alpha \in [0,\infty)$.} Recall $\T_n$ is a uniform plane tree, i.e., sampled from \eqref{eq:T dist} corresponding to $\mu=0$. It follows from the fact that $h(\T_n)$ admits sub-Gaussian tails \cite[Theorem 1.2]{LAB_LD_SJ_13},

$$
Z_n ^\mu = \E\left[e^{-\mu h(\T_n)}\right] \sim C_{n-1} \E \left[e^{-\alpha ||\mathbf{e}||}\right] 
$$
where recall $||\mathbf{e}||= \sup_{t \in[0,1]} \mathbf{e}(t)$ with $\mathbf{e}$ being a standard Brownian excursion. We can then repeat the arguments above and analogously to \eqref{eq: low_bound} obtain a lower bound such that for any $\epsilon >0$ and $n$ sufficiently large
$$\sum_{m=M} ^n e^{-\mu m} \big(H_{n_1,m+1}-H_{n_1,m}\big) \geq (1 - \epsilon)C_{n_1 -1} ~ \E \left[ e^{-\alpha ||\mathbf{e}||} \right]  \,  $$
which, using the well-known fact $C_n \sim \pi^{-1/2} 4^n n^{-3/2}$, simply reproduces the lower bound found in \eqref{eq: ball_vol}. This concludes the proof once multiplied by $K$.
\end{proof}

\begin{proof}[Proof of \cref{thm:loclim}]
As mentioned before, when $\mu$ is constant, the local limit was already identified in \cite{durhuus2023trees}. When $\mu \to \alpha>0$, the same analysis found in \cite[Lemma 4.2]{durhuus2023trees} for the lower bound of the ball volumes can be repeated: for $\delta>0$ and $n$ sufficiently large, $\mu$ can be replaced by $\alpha-\delta$ or $\alpha+\delta$ wherever needed. Since $\delta$ is arbitrary, in the $n \to \infty$ limit one obtains the same ball volume which characterizes the local limit to be the same as $\mu = \alpha$ constant case.

In the $0\leq \mu \ll 1$ case, we apply Lemma \ref{lem:ballvol}. Note that taking $n \rightarrow \infty$ then $M \rightarrow \infty$, the right hand side of \eqref{ballvol} is given by

$$
2 K \cdot 2^K \cdot 4^{-|T_0|} \, 
$$

which is equal to $\nu^{(0)}(\mathcal{B}_{1/r}(T_0))$ where $\nu^{(0)}$ is the law of the uniform {tree $\T_n$}. \cref{rem:loclim} implies that any open subset of $\mathbbm{t}$ can be written as a disjoint union of balls. Then it follows from the Portemanteau theorem that the local limit when $0 \leq \mu \ll 1$ is Kesten's tree.

\end{proof}

\subsection{Condensation at the root away from the local regime}\label{sec:root deg}

The goal of this section is to prove Theorem~\ref{thm:root deg} regarding the asymptotic behaviour of the root degree in a random tree $\mathbf{T}_n$ distributed according to~\eqref{eq:T dist}.
To do so, we start by providing the following proposition, showing the asymptotic probability of the root degree to be a given value in the range of possibilities.

\begin{proposition}\label{prop:root degree}
    Let $\mu=\mu_n$ such that $1\ll\mu\ll n^{1/4}$.
    Let $\mathbf{T}_n$ be a random tree sampled according to~\eqref{eq:T dist}.
    Then, for any fixed $[-R,R]$, and $\rho\in [-R,R]$, we have that
    \begin{align*}
\Prob{\Big.\deg_{\mathbf{T}_n}(\varnothing)=\lfloor2\mu+\rho\sqrt{\mu}\rfloor}\geq\big(1+o(1)\big)\frac{e^{-\rho^2/12}}{\sqrt{12\pi\mu}}\,,
    \end{align*}
    where the $o(1)$ term is uniformly small over $\rho \in [-R,R]$.
\end{proposition}

\begin{proof}
    Assume that $R>0$ is fixed and $\rho\in[-R,R]$.
    We further write $r=\lfloor2\mu+\sqrt{\mu}\rho\rfloor$ for the rest of the proof and assume without loss of generality that $n$ is large enough so that $r<n-1$.
    We recall that $H_{n,m}$ is the number of trees of size $n$ and height $<m$, as defined in~\eqref{eq:Hnm}.

    To start the proof, note that a tree of height $\leq m$ and root degree $r$ can be bijectively mapped with a forest of $r$ trees of height $<m$ and whose total size equals $n-1$, by simply looking at the subtrees or branches emanating from the offsprings of the root. By the last observation, the total number of trees of height $m$ and root degree $r$ is
    \begin{align*}
        \Big|\Big\{T\in\bbt_n:h(T)=m,\deg_{T}(\varnothing)=r\Big\}\Big|=\sum_{n_1+\cdots+n_r=n-1}\Big(H_{n_1,m}\cdots H_{n_r,m}-H_{n_1,m-1}\cdots H_{n_r,m-1}\Big)\,,
    \end{align*}

    where the sum is taken over $n_i\geq1$.
    Using this formula, we obtain that
    \begin{align*}
        \Prob{\deg_{\mathbf{T}_n}(\varnothing)=r}
        &=\frac{1}{Z_n}\sum_{m=2}^{n-1}e^{-\mu m}\sum_{n_1+\ldots+n_r=n-1}\Big(H_{n_1,m}\cdots H_{n_r,m}-H_{n_1,m-1}\cdots H_{n_r,m-1}\Big)\\
        &=\frac{1}{Z_n}\sum_{n_1+\ldots+n_r=n-1}\left[(e^\mu-1)\sum_{m=3}^{n-1}e^{-\mu m}\prod_{i=1}^rH_{n_i,m-1}+e^{-\mu(n-1)}\prod_{i=1}^rC_{n_i-1}\right]\,.
    \end{align*}
    Moreover, by considering an arbitrary $\epsilon=\epsilon_n$ such that $1\gg\epsilon\gg1/(\mu^2n)^{1/6}\asymp (n/\mu^4)^{1/6}/t_n$, Theorem~\ref{thm:inter clt} tells us that we can reduce the previous sum over $m$ by the sum over $|m-t_n|\leq\epsilon t_n$.
    Further using that $\mu\gg1$, we obtain the first approximate lower bound as
    \begin{align*}
        \Prob{\deg_{\mathbf{T}_n}(\varnothing)=r}
        &\geq\big(1+o(1)\big)\frac{e^\mu}{Z_n}\sum_{|m-t_n|\leq\epsilon t_n}e^{-\mu m}\sum_{n_1+\ldots+n_r=n-1}\prod_{i=1}^rH_{n_i,m-1}\,.
    \end{align*}
    In the case of this sum, since $m\geq(1-\epsilon)t_n\gg1$, we can apply Theorem~\ref{thm:ht-bdd tree} stating that, for any $\gamma=\gamma_n\gg1$ and uniformly over any $n_i>\gamma m^2$, we have that
    \begin{align*}
        H_{n_i,m-1}
        &=\frac{4^n}{m}\frac{\tan^2\frac{\pi}{m}}{\left(1+\tan^2\frac{\pi}{m}\right)^n}+O\left(\frac{4^n}{m^3}\frac{1}{\left(1+\tan^2\frac{2\pi}{m}\right)^n}\right)\\
        &=\frac{4^n}{m}\frac{\tan^2\frac{\pi}{m}}{\left(1+\tan^2\frac{\pi}{m}\right)^n}\left[1+O\left(e^{-\frac{3n\pi^2}{m^2}}\right)\right]\,.
    \end{align*}
    Since this approximation only works when $n_i$ is large enough, we split the sum over all $n_i$ in the approximation of the probability according to whether $n_i\leq\gamma m^2$ or $n_i>\gamma m^2$.
    Without loss of generality, and to simplify the upcoming formulas, we further assume that $\gamma$ is chosen so that $\gamma m^2$ is an integer.
    We also observe that there is no constraint on the lower bound of the rate at which $\gamma$ diverges.
    With that in mind, we can add finitely many diverging upper bounds on $\gamma$ and the result would still apply.
    For any such upper bound that we add to $\gamma$, we mention that ``$\gamma$ can be chosen to diverge arbitrarily slowly'', having in mind that this means that the related result holds as long as some undisclosed upper bound holds for $\gamma$.
    
    To start the proof, observe that $m\sim t_n$ here, so that
    \begin{align*}
        \frac{rm^2}{n}\asymp \frac{\mu}{n}\left(\frac{n}{\mu}\right)^{2/3}=\left(\frac{\mu}{n}\right)^{1/3}\,.
    \end{align*}
    Thus, by choosing $\gamma$ diverging slowly enough (here meaning that $\gamma(\mu/n)^{1/3}\ll1$), we can see that it is impossible for all $n_i$ to be less than $\gamma m^2$.
    By further re-ordering the $n_i$ so that the first ones are larger than $\gamma m^2$, we see that
    \begin{align*}
        \sum_{n_1+\ldots+n_r=n-1}\prod_{i=1}^rH_{n_i,m-1}
        &=\sum_{\ell=1}^r\binom{r}{\ell}\sum_{\left\{1\leq n_{\ell+1},\dotsc,n_r\leq\gamma m^2\right\}}\sum_{\left\{n_1+\ldots+n_r=n-1:n_1,\dotsc,n_\ell>\gamma m^2\right\}}\prod_{i=1}^rH_{n_i,m-1}\,.
    \end{align*}
    Write for now $N=n_{\ell+1}+\ldots+n_r=O(\gamma(n^2\mu)^{1/3})=o(n)$.
    Using the approximation for $H_{n_i,m-1}$ when $n_i>\gamma m^2$, we obtain that
    \begin{align*}
        \sum_{\left\{n_1+\ldots+n_r=n-1:n_1,\dotsc,n_\ell>\gamma m^2\right\}}\prod_{i=1}^\ell H_{n_i,m-1}
        &=\sum_{\left\{n_1+\ldots+n_\ell=n-1-N:n_i>\gamma m^2\right\}}\prod_{i=1}^\ell\frac{4^{n_i}}{m}\frac{\tan^2\frac{\pi}{m}}{\left(1+\tan^2\frac{\pi}{m}\right)^{n_i}}\left[1+O\left(e^{-\frac{3n\pi^2}{m^2}}\right)\right]\\
        &=\frac{4^{n-1-N}\pi^{2\ell}}{m^{3\ell}}\frac{1}{\left(1+\tan^2\frac{\pi}{m}\right)^{n-1-N}}\left[1+O\left(e^{-\frac{3n\pi^2}{m^2}}\right)\right]^\ell\binom{n-2-N-\ell\gamma m^2}{\ell-1}\,.
    \end{align*}
    It is now worth adding a constraint on $\ell$ to simplify the computations down the line.

    The previous sum over $\ell$ covers the range $1,\dotsc,r\sim2\mu$.
    However, computations show that only the terms $\ell\sim\mu$ with a window of size $\sqrt{\mu}$ contribute significantly to the sum.
    For this reason, and since we are still considering a lower bound, thus allowing us to remove terms at our own discretion, we now focus on the terms $|\ell-\mu|\leq\gamma\sqrt{\mu}$ (assuming that $\gamma\ll\sqrt{\mu}$, so that this set of $\ell$ is contained within $[1,r]$).
    For any such $\ell$, we now see that
    \begin{align*}
        2+N+\ell\gamma m^2=O(\gamma(n^2\mu)^{1/3})=o(n)\,,
    \end{align*}
    from which we can deduce that
    \begin{align*}
        \binom{n-2-N-\ell\gamma m^2}{\ell-1}
        &=\frac{n^{\ell-1}}{(\ell-1)!}\prod_{i=0}^{\ell-2}\left(1-\frac{2+N+\ell\gamma m^2+i}{n}\right)\\
        &=\frac{n^{\ell-1}}{(\ell-1)!}\exp\left(O\left(\frac{\gamma\ell(n^2\mu)^{1/3}}{n}\right)\right)\,.
    \end{align*}
    Using that we consider $\ell\sim\mu$ here, the $O(\cdot)$ term becomes $O(\gamma(\mu^4/n)^{1/3})$ which, for $\gamma$ diverging slowly enough, is $o(1)$.
    Further observing that $\ell e^{-3n\pi^2/m^2}\ll1$, the above product on $H_{n_i,m-1}$ simplifies to
    \begin{align*}
        \sum_{n_1+\ldots+n_r=n-1}\prod_{i=1}^rH_{n_i,m-1}
        &\geq\big(1+o(1)\big)\sum_{|\ell-\mu|\leq\gamma\sqrt{\mu}}\binom{r}{\ell}\sum_{1\leq n_{\ell+1},\dotsc,n_r\leq\gamma m^2}\prod_{i=\ell+1}^rH_{n_i,m-1}\frac{4^{n-1-N}\pi^{2\ell}}{m^{3\ell}}\frac{1}{\left(1+\tan^2\frac{\pi}{m}\right)^{n-1-N}}\frac{n^{\ell-1}}{(\ell-1)!}\\
        &=\big(1+o(1)\big)\left(\frac{4}{1+\tan^2\frac{\pi}{m}}\right)^{n-1}\sum_{|\ell-\mu|\leq\gamma\sqrt{\mu}}\binom{r}{\ell}\frac{\pi^{2\ell}n^{\ell-1}}{m^{3\ell}(\ell-1)!}\left(\sum_{s=1}^{\gamma m^2}H_{s,m-1}\left(\frac{1+\tan^2\frac{\pi}{m}}{4}\right)^s\right)^{r-\ell}\,,
    \end{align*}
    where we recall that $N=n_{\ell+1}+\ldots+n_r$.

    Applying~\cite[Eq(9)]{de1972average} for $\theta\rightarrow0$, we know that
    \begin{align*}
        \sum_{s=1}^\infty\frac{1}{4^s}H_{s,m-1}=\frac{m-1}{2m}\,.
    \end{align*}
    Further recalling from Theorem~\ref{thm:ht-bdd tree} that
    \begin{align*}
        H_{s,m-1}=O\left(\frac{4^s}{m^3}\right)\,,
    \end{align*}
    we can use the infinite sum to approximate the finite one and obtain
    \begin{align*}
        \sum_{s=1}^{\gamma m^2}H_{s,m-1}\left(\frac{1+\tan^2\frac{\pi}{m}}{4}\right)^s
        &=\sum_{s=1}^{\gamma m^2}\frac{H_{s,m-1}}{4^s}\left[\left(1+\tan^2\frac{\pi}{m}\right)^s-1\right]-\sum_{s>\gamma m^2}\frac{H_{s,m-1}}{4^s}+\frac{m-1}{2m}\\
        &=\sum_{s=1}^{\gamma m^2}\frac{H_{s,m-1}}{4^s}\left[\left(1+\tan^2\frac{\pi}{m}\right)^s-1\right]+\frac{1}{2}+O\left(\frac{1}{m}\right)\,.
    \end{align*}
    Now, in order to bound the difference between the two finite sums, we can use that
    \begin{align*}
        s\mapsto\left[\left(1+\tan^2\frac{\pi}{m}\right)^s-1\right]\frac{m^2e^{-\gamma}}{s}=\left[e^{\frac{s\pi^2}{m^2}+O\left(\frac{s}{m^4}\right)}-1\right]\frac{m^2e^{-9.9\gamma}}{s}
    \end{align*}
    is bounded by some constant $C$ as long as $1\leq s\leq\gamma m^2$ and $\gamma$ diverges slowly enough (we used here that $\pi^2<9.9$).
    It follows that
    \begin{align*}
        \sum_{s=1}^{\gamma m^2}\frac{H_{s,m-1}}{4^s}\left[\left(1+\tan^2\frac{\pi}{m}\right)^s-1\right]\leq Ce^{9.9\gamma}\sum_{s=1}^{\gamma m^2}\frac{H_{s,m-1}}{4^s}\frac{s}{m^2}=O\left(e^{9.9\gamma}\frac{\gamma^2m^4}{m^5}\right)=O\left(\frac{e^{10\gamma}}{m}\right)\,,
    \end{align*}
    where we used again that $H_{s,m-1}=O(4^s/m^3)$ and further that $e^{9.9\gamma}\gamma^2\ll e^{10\gamma}$.
    Pluggin this result back into the approximation for the finite sum, we eventually obtain that
    \begin{align*}
        \sum_{s=1}^{\gamma m^2}H_{s,m-1}\left(\frac{1+\tan^2\frac{\pi}{m}}{4}\right)^s
        &=\frac{1}{2}+O\left(\frac{e^{10\gamma}}{m}\right)\,.
    \end{align*}
    Recall now that $r-\ell\asymp \mu$ along with $m\asymp (n/\mu)^{1/3}$ to see that $(r-\ell)/m\asymp (\mu^4/n)^{1/3}\ll1$, allowing us to choose $\gamma$ diverging slowly enough such that
    \begin{align*}
        \sum_{n_1+\ldots+n_r=n-1}\prod_{i=1}^rH_{n_i,m-1}
        &\geq\big(1+o(1)\big)\left(\frac{4}{1+\tan^2\frac{\pi}{m}}\right)^{n-1}\sum_{|\ell-\mu|\leq\gamma\sqrt{\mu}}\binom{r}{\ell}\frac{\pi^{2\ell}n^{\ell-1}}{m^{3\ell}(\ell-1)!}\frac{1}{2^{r-\ell}}\,.
    \end{align*}
    Moreover, using that $|m-t_n|\leq\epsilon t_n$ along with Proposition~\ref{prop:asympt lx}~\numitem{1}, we see that
    \begin{align*}
        m=\left(\frac{2\pi^2n}{\mu}\right)^{1/3}\left[1+O\left(\left(\frac{\mu}{n}\right)^{2/3}\right)+O(\epsilon)\right]\,.
    \end{align*}
    Using that $\epsilon\gg1/(\mu^2n)^{1/6}$ which implies that $\mu\epsilon\gg(\mu^4/n)^{1/6}$, we can choose $\epsilon$ small enough so that it satisfies the previous assumptions and further $\mu\epsilon\ll1$.
    In that scenario, and using that $\ell\sim\mu$, it follows that
    \begin{align*}
        m^{3\ell}\sim\left(\frac{2\pi^2n}{\mu}\right)^\ell\,,
    \end{align*}
    which finally leads to
    \begin{align*}
        \sum_{n_1+\ldots+n_r=n-1}\prod_{i=1}^rH_{n_i,m-1}
        &\geq\big(1+o(1)\big)\left(\frac{4}{1+\tan^2\frac{\pi}{m}}\right)^{n-1}\sum_{|\ell-\mu|\leq\gamma\sqrt{\mu}}\binom{r}{\ell}\frac{\pi^{2\ell}n^{\ell-1}}{\left(\frac{2\pi^2n}{\mu}\right)^\ell(\ell-1)!}\frac{1}{2^{r-\ell}}\\
        &=\big(1+o(1)\big)\left(\frac{4}{1+\tan^2\frac{\pi}{m}}\right)^{n-1}\frac{1}{n2^r}\sum_{|\ell-\mu|\leq\gamma\sqrt{\mu}}\binom{r}{\ell}\frac{\mu^\ell}{(\ell-1)!}\,.
    \end{align*}
    Relating this formula to the probability of the degree of the root, we eventually obtain that
    \begin{align*}
        \Prob{\deg_{\mathbf{T}_n}(\varnothing)=r}
        &\geq\big(1+o(1)\big)\frac{e^\mu}{Z_n}\sum_{|m-t_n|\leq\epsilon t_n}e^{-\mu m}\left(\frac{4}{1+\tan^2\frac{\pi}{m}}\right)^{n-1}\frac{1}{n2^r}\sum_{|\ell-\mu|\leq\gamma\sqrt{\mu}}\binom{r}{\ell}\frac{\mu^\ell}{(\ell-1)!}\\
        &=\big(1+o(1)\big)\frac{e^\mu}{Z_n}\frac{4^{n-1}}{n2^r}\sum_{|\ell-\mu|\leq\gamma\sqrt{\mu}}\binom{r}{\ell}\frac{\mu^\ell}{(\ell-1)!}\sum_{|m-t_n|\leq\epsilon t_n}e^{-n\lambda_{\mu/n}(m)}\,.
    \end{align*}
    We now use two comparisons to the integral to conclude this proof.

    As a first step, we refer to the proof of Proposition~\ref{prop:W main} to see that
    \begin{align*}
        \sum_{|m-t_n|\leq\epsilon t_n}e^{-n\lambda_{\mu/n}(m)}
        &\sim e^{n\lambda_{\mu/n}(t_n)}\sqrt{2\pi}\frac{1}{\sqrt{3}}\left(\frac{2\pi^2n}{\mu^4}\right)^{1/6}\,,
    \end{align*}
    which along with Theorem~\ref{thm:part-fun} and Proposition~\ref{prop:asympt lx}~\numitem{2} leads to
    \begin{align*}
        \Prob{\deg_{\mathbf{T}_n}(\varnothing)=r}
        &\geq\big(1+o(1)\big)\frac{e^\mu}{4^ne^{2\mu}\frac{\pi^{5/6}\mu^{1/3}}{2^{1/3}3^{1/2}n^{5/6}}}\frac{4^{n-1}}{n2^r}\sqrt{2\pi}\frac{1}{\sqrt{3}}\left(\frac{2\pi^2n}{\mu^4}\right)^{1/6}\sum_{|\ell-\mu|\leq\gamma\sqrt{\mu}}\binom{r}{\ell}\frac{\mu^\ell}{(\ell-1)!}\\
        &=\big(1+o(1)\big)\frac{1}{2e^\mu2^r}\sum_{|\ell-\mu|\leq\gamma\sqrt{\mu}}\binom{r}{\ell}\frac{\mu^{\ell-1}}{(\ell-1)!}\,.
    \end{align*}
    For the final term, recall that $r=\lfloor2\mu+\rho\sqrt{\mu}\rfloor$ with $\rho\in[-R,R]$ and let $\lambda=(\ell-\mu)/\sqrt{\mu}\in[-\gamma,\gamma]$.
    Further assume for now that $2\mu+\rho\sqrt{\mu}$ is an integer (this can be obtained by replacing $\rho$ with $\tilde{\rho}=\rho+O(1/\sqrt{\mu})$, but we prefer avoiding the use of tilde for clarity of the formulas).
    Then, Stirling's formula tells us that
    \begin{align*}
        \binom{r}{\ell}\frac{\mu^{\ell-1}}{(\ell-1)!}
        &=\big(1+o(1)\big)\sqrt{\frac{4\pi\mu}{(2\pi\mu)^3}}\exp\Bigg((2\mu+\rho\sqrt{\mu})\log\Big(2\mu+\rho\sqrt{\mu}\Big)-(\mu+\lambda\sqrt{\mu}-1)\log\left(\frac{\mu+\lambda\sqrt{\mu}-1}{e\mu}\right)\\
        &\hspace{0.5cm}-(\mu+\lambda\sqrt{\mu})\log\Big(\mu+\lambda\sqrt{\mu}\Big)-(\mu+(\rho-\lambda)\sqrt{\mu})\log\Big(\mu+(\rho-\lambda)\sqrt{\mu}\Big)\Bigg)\,.
    \end{align*}
    Observe that, for any constant $c>0$ and $x>0$ possibly diverging, we have
    \begin{align*}
        (c\mu+x\sqrt{\mu})\log\Big(c\mu+x\sqrt{\mu}\Big)
        &=(c\mu+x\sqrt{\mu})\left[\log(c\mu)+\frac{x\sqrt{\mu}}{c\mu}-\frac{x^2\mu}{2c^2\mu^2}+O\left(\frac{x^3}{\mu^{3/2}}\right)\right]\\
        &=(c\mu+x\sqrt{\mu})\log(c\mu)+x\sqrt{\mu}+\frac{x^2}{2c}+O\left(\frac{x^3}{\sqrt{\mu}}\right)\,.
    \end{align*}
    For $x$ being equal to either $\rho$ or $\lambda$, it is bounded by $\gamma$. Thus, we can choose $\gamma$ diverging slowly enough so that the $O(\cdot)$ term in the RHS above is always $o(1)$.
    Plugging this formula in the previous equation, and simplifying the terms, we obtain that
    \begin{align*}
        \binom{r}{\ell}\frac{\mu^{\ell-1}}{(\ell-1)!}
        &=\big(1+o(1)\big)\frac{1}{\sqrt{2}\pi\mu}\exp\Bigg((2\mu+\rho\sqrt{\mu})\log2+\frac{\rho^2}{4}-\frac{\lambda^2}{2}-\frac{(\rho-\lambda)^2}{2}+\mu-\frac{\lambda^2}{2}\Bigg)\\
        &=\big(1+o(1)\big)\frac{e^\mu2^r}{\sqrt{2}\pi\mu}\exp\Bigg(-\frac{\rho^2}{4}-\frac{3\lambda^2}{2}+\rho\lambda\Bigg)\\
        &=\big(1+o(1)\big)\frac{e^\mu2^r}{\sqrt{2}\pi\mu}\exp\Bigg(-\frac{\rho^2}{12}-\frac{3}{2}\left(\lambda-\frac{\rho}{3}\right)^2\Bigg)\,.
    \end{align*}
    Observe that, while we assumed here that $\rho$ was such that $2\mu+\rho\sqrt{\mu}$ was an integer, this assumption can now be removed without loss of generality, since any $O(1/\sqrt{\mu})$ error terms on $\rho$ would not affect the asymptotic behaviour.
    Plugging this back into the lower bound for the probability $\Prob{\deg_{\mathbf{T}_n}(\varnothing)=r}$, we now see that
    \begin{align*}
        \Prob{\deg_{\mathbf{T}_n}(\varnothing)=r}
        &\geq\big(1+o(1)\big)\frac{1}{2e^\mu2^r}\sum_{|\ell-\mu|\leq\gamma\sqrt{\mu}}\frac{e^\mu2^r}{\sqrt{2}\pi\mu}\exp\Bigg(-\frac{\rho^2}{12}-\frac{3}{2}\left(\lambda-\frac{\rho}{3}\right)^2\Bigg)\\
        &=\big(1+o(1)\big)\frac{e^{-\rho^2/12}}{\sqrt{8}\pi\mu}\sum_{|\ell-\mu|\leq\gamma\sqrt{\mu}}\exp\Bigg(-\frac{3}{2}\left(\frac{\ell-\mu}{\sqrt{\mu}}-\frac{\rho}{3}\right)^2\Bigg)\,.
    \end{align*}
    To finally conclude this proof, observe that the Riemann approximation of the integral applies here, leading to
    \begin{align*}
        \sum_{|\ell-\mu|\leq\gamma\sqrt{\mu}}\exp\Bigg(-\frac{3}{2}\left(\frac{\ell-\mu}{\sqrt{\mu}}-\frac{\rho}{3}\right)^2\Bigg)\sim\sqrt{\mu}\int_\R e^{-\frac{3(x-\rho)^2}{2}}dx=\sqrt{\frac{2\mu\pi}{3}}\,.
    \end{align*}
    The desired result then follows from plugging this formula in the previous equation.
\end{proof}

\begin{proof}[Proof of Theorem~\ref{thm:root deg}]
    With this result the proof of the asymptotic distribution of the root is quite straightforward.
    Indeed, for any $t\in\R$, we see that
    \begin{align*}
        \Prob{\frac{\deg_{\bT_n}(\varnothing)-2\mu}{\sqrt{6\mu}}\leq t}
        &=\Prob{\deg_{\bT_n}(\varnothing)\leq\lfloor2\mu+t\sqrt{6\mu}\rfloor\Big.}\\
        &\geq\sqrt{6\mu}\int_{-R}^t\Prob{\deg_{\bT_n}(\varnothing)=\lfloor2\mu+u\sqrt{6\mu}\rfloor\Big.}du\,,
    \end{align*}
    for any $R>t$. 
    Thus, by applying Proposition~\ref{prop:root degree}, it follows that
    \begin{align*}
        \Prob{\frac{\deg_{\bT_n}(\varnothing)-2\mu}{\sqrt{6\mu}}\leq t}
        &\geq\big(1+o(1)\big)\sqrt{6\mu}\int_{-R}^t\frac{e^{-(\sqrt{6}u)^2/12}}{\sqrt{12\pi\mu}}du.
    \end{align*}
    By first letting $n \to \infty$ using that the $o(1)$ is uniformly small in $[-R,R]$, and then letting $R \to \infty$, the right-hand side above equals $\frac{1}{\sqrt{2\pi}}\int_{-\infty}^te^{-u^2/2}du$, which is the distribution function of the standard normal random variable. The Portmanteau theorem now finishes the proof.
\end{proof}

\begin{remark}[Local CLT]
    In fact, the same argument above shows that the inequality in the statement of Proposition~\ref{prop:root degree} is in fact an equality, yielding a \emph{local CLT}. 
\end{remark}

\medskip

\newpage

\appendix

\section{Bounding functions with polynomial factor and exponent}\label{sec:func_poly_fac}

This section focuses on proving the single result below, which we later use for the proof of both Theorem~\ref{thm:ht-bdd tree} and Theorem~\ref{thm:W asympt}.
The proof uses a method where exponentials are iterated until reaching the desired behaviour.

\begin{lemma}\label{lem:useful-bound}
    For any real number $\alpha\in\R$ and any two strictly positive numbers $\beta,\gamma>0$, consider the function
    \begin{align*}
        F_t(x,y)=\frac{t^\alpha}{x^\beta}\exp\left(-\frac{t}{y^\gamma}\right)\,.
    \end{align*}
    Let $a_+=a_+(t)$ be strictly positive and bounded, and $\epsilon=\epsilon(t)$ converging to $0$ as $t$ diverges. 
    Then there exists $k=k(t)$ an integer depending on $t$ and a sequence of increasing values $(a_0\leq\ldots\leq a_k)$ (implicitly depending on $t$) such that $a_0=\epsilon$, $a_k=a_+$, and satisfying
    \begin{align*}
        \lim_{t\rightarrow\infty}\sum_{i=1}^kF_t(a_{i-1},a_i)=0\,.
    \end{align*}
\end{lemma}

\begin{proof}
    In order to prove the theorem, we directly provide a proper sequence and show that the desired sum converges to $0$.
    Before doing so, note that requiring $a_0=\epsilon$ is equivalent to requiring $a_0\leq\epsilon$ as the sum is a decreasing function of $a_0$.
    Similarly, by observing that $F(a,b)$ converges to $0$ for any fix $a,b>0$ and that the sum is increasing with $a_k$ (at fixed $k$), we can assume without loss of generality that $a_+=1$.
    Consider now the sequence defined inductively by $b_0=1$ and
    \begin{align*}
        b_{i+1}=\exp\left(\frac{\gamma t}{2\beta}b_i\right)\,.
    \end{align*}
    Further define $k=\min\{i:b_i\geq1/\epsilon^\gamma\}$ and let $a_i=b_{k-i}^{-1/\gamma}$ for $0\leq i\leq k$.
    We prove that $(a_0,\ldots,a_k)$ satisfies all the assumptions of the lemma.

    Before proving that the desired sum converges to $0$, let us show that the initial assumptions are met.
    By definition, we see that $a_k=b_0^{-1/\gamma}=1$ since $b_0=1$.
    Moreover, by the definition of $k$, we know that $b_k\geq1/\epsilon^\gamma$ and so $a_0=b_k^{-1/\gamma}\leq\epsilon$.
    Finally, since $e^x\geq x$ for any $x\in\R$ and using that $t$ can be chosen arbitrarily large, for example such that $t\geq2\beta/\gamma$, we see that $b_{i+1}\geq b_i$ and so $a_{i+1}\geq a_i$.
    This proves that all initial assumptions are met and we now focus on proving the convergence of the sum.

    Start by observing that, for any $1\leq i\leq k$, we have
    \begin{align*}
        F_t(a_{i-1},a_i)&=\frac{t^\alpha}{a_{i-1}^\beta}\exp\left(-\frac{t}{a_i^\gamma}\right)=t^\alpha b_{k-i+1}^{\beta/\gamma}\exp\Big(-tb_{k-i}\Big)=t^\alpha\exp\left(-\frac{t}{2}b_{k-i}\right)\,.
    \end{align*}
    It follows that the sum can be re-written as
    \begin{align*}
        \sum_{i=1}^kF_t(a_{i-1},a_i)=\sum_{i=0}^{k-1}t^\alpha\exp\left(-\frac{t}{2}b_i\right)\,.
    \end{align*}
    But now, use $e^x\geq x+1$ for any $x\in\R$ and a simple induction to see that $b_i\geq i$, which leads to
    \begin{align*}
        \sum_{i=1}^kF_t(a_{i-1},a_i)\leq\sum_{i=0}^{k-1}t^\alpha\exp\left(-\frac{t}{2}i\right)\leq\frac{t^\alpha e^{-t/2}}{1-e^{-t/2}}\,,
    \end{align*}
    and the desired convergence to $0$ follows.
\end{proof}

\begin{corollary}\label{cor:useful-bound}
    Let $\alpha\in\R$ be a real number and $B,\Gamma$ be two continuous functions.
    Define the function
    \begin{align*}
        F_t(x,y)=\frac{t^\alpha}{B(x)}\exp\left(-\frac{t}{\Gamma(y)}\right)\,.
    \end{align*}
    Assume that $B,\Gamma$ are strictly positive on $(0,a_+]$ and that there exists $C_\beta,C_\gamma>0$ and $\beta,\gamma>0$ such that $B(x)\sim C_\beta x^\beta$ and $\Gamma(x)\sim C_\gamma x^\gamma$ when $x\rightarrow0$.
    Then, for any $\epsilon=\epsilon(t)$ converging to $0$ as $t$ diverges and $\Tilde{a}_+=\Tilde{a}_+(t)\leq a_+$, there exists $k=k(t)$ an integer depending on $t$ and a sequence of increasing values $(a_0\leq\ldots\leq a_k)$ (implicitly depending on $t$) such that $a_0=\epsilon$, $a_k=\Tilde{a}_+$, and satisfying
    \begin{align*}
        \lim_{t\rightarrow\infty}\sum_{i=1}^kF_t(a_{i-1},a_i)=0\,.
    \end{align*}
\end{corollary}

\begin{proof}
    From the assumption on $B$ and $\Gamma$, we can find $C>0$ such that, for any $0\leq x\leq\tilde{a}_+$, we have $B(x)\geq C^{-1}x^\beta$ and $\Gamma(x)\leq Cx^\gamma$.
    Now, let $\epsilon=a_0\leq\cdots\leq a_k=\tilde{a}_+$ be the sequence satisfying
    \begin{align*}
        \lim_{t\rightarrow\infty}\sum_{i=1}^k\frac{t^\alpha}{a_{i-1}^\beta}\exp\left(-\frac{t}{a_i^\gamma}\right)=0\,,
    \end{align*}
    which exists according to Lemma~\ref{lem:useful-bound}.
    Then, we observe that this allows to bound the sum of $F_t$, as
    \begin{align*}
        \sum_{i=1}^{k+1}F_t(a_{i-1},a_i)
        &\leq\sum_{i=1}^k\frac{Ct^\alpha}{a_{i-1}^\beta}\exp\left(-\frac{t}{Ca_i^\gamma}\right)=C^{1+\alpha}\sum_{i=1}^k\frac{(t')^\alpha}{a_{i-1}^\beta}\exp\left(-\frac{t'}{a_i^\gamma}\right)\,,
    \end{align*}
    where we used the change of variable $t'=t/C$.
    Using that $t\rightarrow\infty$ if and only if $t'\rightarrow\infty$ along with the definition of the $a_i$ proves the desired convergence to $0$.
\end{proof}

\section{Trees of bounded height}\label{sec:ht-bdd tree}

In this section, we provide some useful properties regarding the number of trees of bounded height.
Recall from \eqref{eq:Hnm} we let $H_{n,m}$ be the number of trees of size $n$ and height less than $m$, that is
\begin{align*}
    H_{n,m}=\Big|\big\{T\in\bbt_n:h(t)<m\big\}\Big|\,.
\end{align*}
A closed formula for $H_{n,m}$ is already known~\cite[(14)]{de1972average} and is given by
\begin{align}\label{eq:ht-bdd tree}
    H_{n,m}=\frac{4^n}{m+1}\sum_{k=1}^{\lfloor m/2\rfloor}\sin^2\left(\frac{\pi k}{m+1}\right)\cos^{2n-2}\left(\frac{\pi k}{m+1}\right)=\frac{4^n}{m+1}\sum_{k=1}^{\lfloor m/2\rfloor}\frac{\tan^2\left(\frac{\pi k}{m+1}\right)}{\left(1+\tan^2\left(\frac{\pi k}{m+1}\right)\right)^n}\,.
\end{align}
Using this formula, we now prove the following theorem, bounding $H_{n,m}$ for some choices of $m$.

\begin{theorem}\label{thm:ht-bdd tree}
    Let $h_n$ be such that $4\leq h_n\ll\sqrt{n}$.
    Then, for any $4\leq m\leq h_n$, we have
    \begin{align*}
        H_{n,m}&=\frac{4^n}{m+1}\frac{\tan^2\left(\frac{\pi}{m+1}\right)}{\left(1+\tan^2\left(\frac{\pi}{m+1}\right)\right)^n}+O\left(\frac{4^n}{m^3}\frac{1}{\left(1+\tan^2\left(\frac{2\pi}{m+1}\right)\right)^n}\right)\\
        &=\big(1+o(1)\big)\frac{4^n}{m+1}\frac{\tan^2\left(\frac{\pi}{m+1}\right)}{\left(1+\tan^2\left(\frac{\pi}{m+1}\right)\right)^n}
    \end{align*}
    where the $o(\cdot)$ and $O(\cdot)$ terms only depend on $h_n$, and are uniformly bounded over all such $m$.
    Moreover, the second asymptotic equality trivially holds when $m\in\{2,3\}$.
\end{theorem}

\begin{proof}
    Let $h_n$ satisfy the conditions and consider $4\leq m\leq h_n$.
    Observe that the result is equivalent to saying that the first two terms of~\eqref{eq:ht-bdd tree} dominate the sum, with the first term being the dominant one.
    The fact that the second term is negligible with respect to the first one directly follows from their ratio being an exponential with a diverging term (recalling that $m\leq h_n$ and $n/h_n^2\gg1$).
    Thus, we now prove that the term $k=2$ dominates all the terms corresponding to $k\geq2$ in~\eqref{eq:ht-bdd tree}.

    Consider first the case $m$ finite.
    In that case $H_{n,m}$ is a finite sum of exponential terms in $n$ with negative exponents.
    Moreover all exponents are distinct, so only the exponent of minimal value contributes to the asymptotic behaviour.
    Further using that $m\geq4$ so that the sum has at least two terms, we see that the given asymptotic behaviour holds.
    Since this is true for $m$ fixed and $n\rightarrow\infty$, we can also directly extend it to $\Tilde{m}=\Tilde{m}_n$ diverging slowly enough (for example $\tilde{m}=(1-\epsilon)\sqrt{3\pi^2n/\log n}$ for any $\epsilon>0$ works, but the exact rate is not relevant to the rest of the proof).
    Thus the first statement of the lemma is proven for all $4\leq m\leq\Tilde{m}$ and we now focus on $\Tilde{m}\leq m\leq h_n$ (assuming that $\Tilde{m}\leq h_n$).
    For the rest of the proof we repeatedly use that, under the assumption that $\Tilde{m}\leq m\leq h_n\ll\sqrt{n}$, any $o(\cdot)$ and $O(\cdot)$ terms related to the fact that $m$ or $n/m^2$ both diverge to infinity is uniform over all such $m$.

    To extend the result of the lemma to all $m\leq h_n$, we first prove that the terms of the sum with $k>\log m$ are negligible.
    Let us observe that, for any $\alpha,\beta>0$ such that $\alpha+\beta<\pi/2$, we have
    \begin{align*}
        1+\tan^2(\alpha+\beta)
        &=1+\left(\frac{\tan\alpha+\tan\beta}{1-\tan\alpha\tan\beta}\right)^2=\frac{\big(1+\tan^2\alpha\big)\big(1+\tan^2\beta\big)}{\big(1-\tan\alpha\tan\beta\big)^2}\geq\big(1+\tan^2\alpha\big)\big(1+\tan^2\beta\big)\,.
    \end{align*}
    Following this inequality, a simple induction leads to
    \begin{align*}
        1+\tan^2\left(\frac{\pi k}{m+1}\right)\geq\left(1+\tan^2\left(\frac{\pi}{m+1}\right)\right)^k\,.
    \end{align*}
    Using this bound along with the fact that $k\mapsto\tan^2(\pi k/(m+1))$ is increasing, we observe that
    \begin{align*}
        \sum_{\log m\leq k\leq m/2}\frac{\tan^2\left(\frac{\pi k}{m+1}\right)}{\left(1+\tan^2\left(\frac{\pi k}{m+1}\right)\right)^n}&\leq\tan^2\left(\frac{\pi m}{2(m+1)}\right)\sum_{\log m\leq k\leq m/2}\frac{1}{\left(1+\tan^2\left(\frac{\pi}{m+1}\right)\right)^{kn}}\\
        &\leq\tan^2\left(\frac{\pi}{2}-\frac{\pi}{2(m+1)}\right)\frac{1}{\left(1+\tan^2\left(\frac{\pi}{m+1}\right)\right)^{n\log m}}\frac{1}{1-\left(1+\tan^2\left(\frac{\pi}{m+1}\right)\right)^{-n}}\\
        &=\big(1+o(1)\big)\frac{4m^2}{\pi^2}e^{-\frac{\pi n}{(m+1)^2}\log m+O\left(\frac{n\log m}{m^4}\right)}
    \end{align*}
    Recall that the second term is of order
    \begin{align*}
        \frac{\tan^2\left(\frac{2\pi}{m+1}\right)}{\left(1+\tan^2\left(\frac{2\pi}{m+1}\right)\right)^n}=\big(1+o(1)\big)\frac{4\pi^2}{m^2}e^{-\frac{4n\pi^2}{(m+1)^2}+O\left(\frac{n}{m^4}\right)}
    \end{align*}
    to see that
    \begin{align*}
        \frac{\left(1+\tan^2\left(\frac{2\pi}{m+1}\right)\right)^n}{\tan^2\left(\frac{2\pi}{m+1}\right)}\sum_{\log m\leq k\leq m/2}\frac{\tan^2\frac{\pi k}{m+1}}{\left(1+\tan^2\frac{\pi k}{m+1}\right)^n}&\leq\big(1+o(1)\big)\frac{m^4}{\pi^4}e^{-(1+o(1))\frac{n\pi^2}{m^2}\log m}\,.
    \end{align*}
    The fact that the sum for $k\geq\log m$ is of smaller order than the second term simply follows from the fact that $n/m^2$ diverges and that $m^4e^{-x\log m}$ converges to $0$ as $m,x\rightarrow\infty$.
    We now prove that the rest of the sum is negligible.

    Observe that, uniformly over any $\Tilde{m}\leq m\leq h_n$ and $2\leq k<\log m$, we have
    \begin{align*}
        \frac{\tan^2\left(\frac{\pi k}{m+1}\right)}{\left(1+\tan^2\left(\frac{\pi k}{m+1}\right)\right)^n}&=\big(1+o(1)\big)\frac{\pi^2k^2}{m^2}e^{-\frac{\pi^2k^2n}{(m+1)^2}+O\left(\frac{nk^4}{m^4}\right)}=e^{-(1+o(1))\frac{\pi^2k^2n}{m^2}}\,.
    \end{align*}
    From this behaviour, we see that the terms are decreasing in $k$.
    Thus, for any $3\leq a\leq b<\log m$, using the minimal term to bound the sum along with the fact that $a\leq b$, it follows
    \begin{align*}
        \frac{\left(1+\tan^2\left(\frac{2\pi}{m+1}\right)\right)^2}{\tan^2\left(\frac{2\pi}{m+1}\right)}\sum_{k=a}^b\frac{\tan^2\left(\frac{\pi k}{m+1}\right)}{\left(1+\tan^2\left(\frac{\pi k}{m+1}\right)\right)^n}&\leq\big(1+o(1)\big)b^3e^{-(1+o(1))\frac{\pi^2(a^2-4)n}{(m+1)^2}}\,,
    \end{align*}
    where all $o(\cdot)$ terms are uniform over all the parameters considered here.
    Apply now Corollary~\ref{cor:useful-bound} with $\alpha=0$, $B(x)=x^3$, $\Gamma(x)=x^2/(4-x^2)$, $a_+=1/3$, mapping $t$ with $\pi^2n/(m+1)^2$, and $\epsilon=1/\log n$ to find $a_0\leq\cdots\leq a_K$ such that
    \begin{align*}
        \lim_{n\rightarrow\infty}\sum_{i=1}^K\frac{1}{a_{i-1}^2}\exp\left(-\left(\frac{1}{a_i^2}-1\right)\frac{\pi^2n}{(m+1)^2}\right)=0\,.
    \end{align*}
    Using this sequence and combining with the previous inequality, we obtain that
    \begin{align*}
        \frac{\left(1+\tan^2\left(\frac{2\pi}{m+1}\right)\right)^2}{\tan^2\left(\frac{2\pi}{m+1}\right)}\sum_{3\leq k<\log m}\frac{\tan^2\left(\frac{\pi k}{m+1}\right)}{\left(1+\tan^2\left(\frac{\pi k}{m+1}\right)\right)^n}&=\sum_{i=1}^K\frac{\left(1+\tan^2\left(\frac{2\pi}{m+1}\right)\right)^2}{\tan^2\left(\frac{2\pi}{m+1}\right)}\sum_{1/a_i\leq k<1/a_{i-1}}\frac{\tan^2\left(\frac{\pi k}{m+1}\right)}{\left(1+\tan^2\left(\frac{\pi k}{m+1}\right)\right)^n}\\
        &\leq\big(1+o(1)\big)\sum_{i=1}^K\frac{1}{a_{i-1}^2}\exp\left(-\left(\frac{1}{a_i^2}-1\right)\frac{\pi^2n}{(m+1)^2}+O\left(\frac{n\log^4n}{m^4}\right)\right)\,,
    \end{align*}
    where we used the monotonicity of the upper bound (increasing with $a_i$ and decreasing with $a_{i-1}$) to relax the integer constraints.
    By the definition of the $a_i$ terms, the upper bound converges to $0$ and so the sum of the terms from $3$ to $\log n$ are negligible with respect to the second term.
    Combined with the fact that the sum after $\log n$ is also negligible proves the desired result.
\end{proof}

The number of trees of bounded height is very useful in order to characterize the behaviour of the partition function $Z_n$ introduced in~\eqref{eq:part fun}, as its definition implies that
\begin{align*}
    Z_n^\mu=\sum_{T\in\bbt_n}e^{-\mu h(T)}=\sum_{m=1}^n\big(H_{n,m+1}-H_{n,m}\big)e^{-\mu m}\,.
\end{align*}
Moreover, using that $H_{n,m}=C_{n-1}$ whenever $m\geq n$ and that $H_{n,1}=0$ (assuming $n>1$), we see that the previous formula re-writes as
\begin{align}\label{eq:part fun & W}
    Z_n^\mu=C_{n-1}e^{-\mu(n-1)}+e^\mu(e^\mu-1)\sum_{m=2}^{n-1}H_{n,m}e^{-\mu(m+1)}\,.
\end{align}
Thus, it is rather natural to introduce the sum
\begin{align}\label{eq:def W}
    W_n=W_n^\mu=\frac{1}{4^n}\sum_{m=3}^{n}H_{n,m-1}e^{-\mu m}\,,
\end{align}
which itself contains all of the non-trivial asymptotic behaviour of $Z_n$, as stated by the following lemma.

\begin{lemma}\label{lem:z and w}
    Let $\mu=\mu_n$ be such that $\mu\gg1/\sqrt{n}$.
    Then, $Z_n$ as defined in~\eqref{eq:part fun} is related to $W_n$ introduced in~\eqref{eq:def W} by
    \begin{align*}
        Z_n=\big(1+o(1)\big)4^ne^\mu(e^\mu-1)W_n
    \end{align*}
\end{lemma}

\begin{proof}
    The statement of this lemma simply says that the term $C_{n-1}e^{-\mu(n-1)}$ can be ignored in the formula~\eqref{eq:part fun & W}.
    To see that it is the case, first note that, for any $m \leq n$,
    \begin{align*}
        W_n\geq\frac{1}{4^n}H_{n,m}e^{-\mu m}\geq\frac{\tan^2\frac{\pi}{m}}{m}\frac{e^{-\mu m}}{\left(1+\tan^2\frac{\pi}{m}\right)^n}\,,
    \end{align*}
    from which it follows that, for any $m\gg1$, we have
    \begin{align*}
        \frac{C_{n-1}e^{-\mu(n-1)}}{4^ne^\mu(e^\mu-1)W_n}\leq\big(C+o(1)\big)\frac{m^3}{n^{3/2}}\frac{1}{e^\mu-1}\exp\left(-\mu(n-m)+\frac{n\pi^2}{m^2}+O\left(\frac{n}{m^4}\right)\right)\,.
    \end{align*}
    Use now $m=\lfloor\sqrt{n}\rfloor$ to obtain that
    \begin{align*}
        \frac{C_{n-1}e^{-\mu(n-1)}}{4^ne^\mu(e^\mu-1)W_n}\leq\big(C+o(1)\big)\frac{e^{-(1+o(1))\mu n}}{e^\mu-1}\,.
    \end{align*}
    Recalling that $\mu\gg1/\sqrt{n}$, we see that the right-hand side necessarily converges to $0$ and so the asymptotic behaviour of $Z_n$ is characterized by the term $W_n$.
    Observe that the convergence occurs further away than $\mu\gg1/\sqrt{n}$ and can be pushed as far as $\mu\gg(\log n)/n$ with the same method.
    However, we do not need this result to go further than $1/\sqrt{n}$ and so keep it as is here.
\end{proof}

\section{Exponent of the partition function}\label{sec:on part fun}

Thanks to Lemma~\ref{lem:z and w}, we see that characterizing the behaviour of $Z_n$ is similar to characterizing the behaviour of $W_n$.
Furthermore, combining the formula of $W_n$ from~\eqref{eq:def W} along with Theorem~\ref{thm:ht-bdd tree}, we see that, approximately,
\begin{align*}
    W_n\asymp \sum_{m=3}^n\frac{\tan^2\left(\frac{\pi}{m}\right)}{m}\frac{e^{-\mu m}}{\left(1+\tan^2\left(\frac{\pi}{m}\right)\right)^n}\,.
\end{align*}
Observe that this approximation only works as long as the trees contributing to the partition function all have their height of smaller order than $\sqrt{n}$.
While this might seem like a circular argument, we in fact prove that only the trees of height $\ll\sqrt{n}$ contribute to $W_n$, and only intend to provide a high level view of the argument here.

In order to characterize this sum, we see that the first fraction in the sum behaves as a polynomial, whereas the second one has exponential terms.
It is thus natural to focus on the behaviour of the exponential terms in order to characterize $W_n$, as we expect them to be the main contributor to the sum.
We now use this observation to define and study the exponent function.

{At this point, it is useful to redefine for convenience the function $\lambda_x$ for any $x>0$ and its minimizer $t_x$, already defined in \eqref{eq:func_lambdax_minvalue}},
\begin{align}\label{eq:lx}
    \lambda_x(t)=xt+\log\left(1+\tan^2\frac{\pi}{t}\right)\,,
\end{align}
whose domain is $t\in(2,\infty)$.
From this definition, we directly see that $\lambda_x$ diverges to $\infty$ when $t$ approaches any of its boundary, implying that it admits a minimum (non necessarily unique).
Further noting that
\begin{align*}
    \lambda_x'(t)=x-\frac{2\pi}{t^2}\tan\frac{\pi}{t}\,,
\end{align*}
we see that $\lambda_x$ is strictly convex and thus the minimum $t_x$ is unique. Observe that $t_x$ satisfies
\begin{align*}
    x=\frac{2\pi}{t_x^2}\tan\frac{\pi}{t_x}\,.
\end{align*}
We now provide an important result regarding the asymptotic behaviour of $t_x$ and $\lambda_x$ when $x$ converges to $0$.

\begin{proposition}\label{prop:asympt lx}
    Let $x>0$ and $\lambda_x$ be the function defined in~\eqref{eq:lx}.
    Let $t_x$ be the minimizer of $\lambda_x$.
    Then, as $x$ converges to $0$, we have the following asymptotic behaviour.
    \begin{itemize}
        \item[\numitem{1}] The minimizer $t_x$ of the function satisfies
        \begin{align*}
            t_x=\left(\frac{2\pi^2}{x}\right)^{1/3}\left[1+\frac{1}{9}\left(\frac{\pi x}{2}\right)^{2/3}+O\left(x^{4/3}\right)\right]=\left(\frac{2\pi^2}{x}\right)^{1/3}+O\left(x^{1/3}\right)
        \end{align*}
        \item[\numitem{2}] The minimum of the function satisfies
        \begin{align*}
            \lambda_x(t_x)=3\left(\frac{\pi x}{2}\right)^{2/3}+\frac{1}{6}\left(\frac{x}{2\pi}\right)^{4/3}+O\left(x^2\right)\,,
        \end{align*}
        \item[\numitem{3}] For any $\epsilon=\epsilon_x$ converging to $0$ as $x$ goes to $0$ (allowed to be negative), the asymptotic behaviour of $\lambda_x$ at $(1+\epsilon)t_x$ is given by
        \begin{align*}
            \lambda_x\big((1+\epsilon)t_x\big)-\lambda_x(t_x)&=3\left(\frac{\pi x}{2}\right)^{2/3}\epsilon^2+O\left(\epsilon^3x^{2/3}+\epsilon^2x^{4/3}\right)\,,
        \end{align*}
        and the $O(\cdot)$ term is uniform over any window of the form $[-\delta,\delta]$ for $\delta=\delta_x\ll1$.
    \end{itemize}
\end{proposition}

\begin{proof}
    We start with the behaviour of $t_x$ as given in~\numitem{1}.
    Since the second equality directly follows from the first one, we only show that the first asymptotic holds. 
    While the method to obtain this result in the first place requires inductively deducting the terms of $t_x$, we now simply verify that the given formula satisfies our assumption.
    Indeed, by first observing that $t_x$ necessarily diverges to $\infty$ when $x$ converges to $0$, and by using the Taylor expansion of the tangent function, we obtain
    \begin{align*}
        \frac{2\pi}{t_x^2}\tan\left(\frac{\pi}{t_x}\right)&=\frac{2\pi}{t_x^2}\left(\frac{\pi}{t_x}+\frac{\pi^3}{3t_x^3}+O\left(\frac{1}{t_x^5}\right)\right)=\frac{2\pi^2}{t_x^3}+\frac{2\pi^4}{3t_x^5}+O\left(\frac{1}{t_x^7}\right)\,.
    \end{align*}
    Thus, replacing $t_x$ by the given formula leads to
    \begin{align*}
        \frac{2\pi}{t_x^2}\tan^2\left(\frac{\pi}{t_x}\right)&=\frac{2\pi^2}{\frac{2\pi^2}{x}}\left[1-\frac{1}{3}\left(\frac{\pi x}{2}\right)^{2/3}+O\left(x^{4/3}\right)\right]+\frac{2\pi^4}{3\left(\frac{2\pi^2}{x}\right)^{5/3}}\left[1+O\left(x^{2/3}\right)\right]+O\left(x^{7/3}\right)\\
        &=x+\left[\frac{2\pi^4}{3(2\pi^2)^{5/3}}-\frac{1}{3}\left(\frac{\pi}{2}\right)^{2/3}\right]x^{5/3}+O\left(x^{7/3}\right)\,.
    \end{align*}
    Observe that the $x^{5/3}$ order term in the previous formula is actually multiplied by $0$.
    This implies that the approximation is correct and furthermore we see that the $O(\cdot)$ terms are exactly the right order as it matches the terms from the expansion of the tangent function.
    We now move on to the asymptotic behaviour of $\lambda_x$.

    Starting with the definition of $\lambda_x$ from~\eqref{eq:lx}, since $t_x$ diverges, we have
    \begin{align*}
        \lambda_x(t_x)&=xt_x+\log\left(1+\left(\frac{\pi}{t_x}+\frac{\pi^3}{3t_x^3}+O\left(\frac{1}{t_x^5}\right)\right)^2\right)\\
        &=xt_x+\log\left(1+\frac{\pi^2}{t_x^2}\left(1+\frac{2\pi^2}{3t_x^2}+O\left(\frac{1}{t_x^4}\right)\right)\right)\\
        &=xt_x+\left[\frac{\pi^2}{t_x^2}\left(1+\frac{2\pi^2}{3t_x^2}+O\left(\frac{1}{t_x^4}\right)\right)-\frac{1}{2}\frac{\pi^4}{t_x^4}\left(1+O\left(\frac{1}{t_x^2}\right)\right)+O\left(\frac{1}{t_x^6}\right)\right]\\
        &=xt_x+\frac{\pi^2}{t_x^2}+\frac{\pi^4}{6t_x^4}+O\left(\frac{1}{t_x^6}\right)\,.
    \end{align*}
    To prove~\numitem{2}, plug in this equation the previous asymptotic behaviour of $t_x$ to obtain
    \begin{align*}
        \lambda_x(t_x)&=x\left(\frac{2\pi^2}{x}\right)^{1/3}\left[1+\frac{1}{9}\left(\frac{\pi x}{2}\right)^{2/3}+O\left(x^{4/3}\right)\right]+\frac{\pi^2}{\left(\frac{2\pi^2}{x}\right)^{2/3}}\left[1-\frac{2}{9}\left(\frac{\pi x}{2}\right)^{2/3}+O\left(x^{4/3}\right)\right]\\
        &\hspace{0.5cm}+\frac{\pi^4}{6\left(\frac{2\pi^2}{x}\right)^{4/3}}\left[1+O\left(x^{2/3}\right)\right]+O\left(x^2\right)\\
        &=3\left(\frac{\pi x}{2}\right)^{2/3}+\frac{1}{6}\left(\frac{x}{2\pi}\right)^{4/3}+O\left(x^2\right)\,,
    \end{align*}
    as desired.

    Finally, for~\numitem{3}, start with the definition of $\lambda_x$ from~\eqref{eq:lx} to see that
    \begin{align*}
        \lambda_x\big((1+\epsilon)t_x\big)-\lambda_x(t_x)&=x\epsilon t_x+\log\left(\frac{1+\tan^2\left(\frac{\pi}{(1+\epsilon)t_x}\right)}{1+\tan^2\left(\frac{\pi}{t_x}\right)}\right)\,.
    \end{align*}
    Now, using the trigonometric identity
    \begin{align*}
        1+\tan^2(\alpha+\beta)&=1+\left(\frac{\tan\alpha+\tan\beta}{1-\tan\alpha\tan\beta}\right)^2=\frac{\big(1+\tan^2\alpha\big)\big(1+\tan^2\beta\big)}{\big(1-\tan\alpha\tan\beta\big)^2}
    \end{align*}
    with $\alpha=\pi/t_n$ and $\beta=-\pi\epsilon_n/((1+\epsilon_n)t_n)$ leads to
    \begin{align*}
        \lambda_x\big((1+\epsilon)t_x\big)-\lambda_x(t_x)&=x\epsilon t_x+\log\left(\frac{1+\tan^2\left(\frac{\pi\epsilon}{(1+\epsilon)t_x}\right)}{\left(1+\tan\frac{\pi}{t_x}\tan\left(\frac{\pi\epsilon}{(1+\epsilon)t_x
        }\right)\right)^2}\right)\,.
    \end{align*}
    Using that $\epsilon\ll1$ and $t_x\gg1$, all the terms within the logarithm can be extended using Taylor approximation.
    Moreover, we recall that $t_x$ being the minimum of $\lambda_x$, it satisfies
    \begin{align*}
        x=\frac{2\pi}{t_x^2}\tan\left(\frac{\pi}{t_x}\right)\,,
    \end{align*}
    which leads to the following approximation
    \begin{align*}
        \lambda_x\big((1+\epsilon)t_x\big)-\lambda_x(t_x)&=x\epsilon t_x+\left[\frac{\pi^2\epsilon^2}{(1+\epsilon)^2t_x^2}-2\tan\left(\frac{\pi}{t_x}\right)\frac{\pi\epsilon}{(1+\epsilon)t_x}+O\left(\frac{\epsilon^2}{t_x^4}\right)\right]\\
        &=x\epsilon t_x+\frac{\pi^2\epsilon^2}{(1+\epsilon)^2t_x^2}-\frac{x\epsilon t_x}{1+\epsilon}+O\left(\frac{\epsilon^2}{t_x^4}\right)\\
        &=\left[(1+\epsilon)xt_x+\frac{\pi^2}{t_x^2}\right]\frac{\epsilon^2}{(1+\epsilon)^2}+O\left(\frac{\epsilon^2}{t_x^4}\right)\,.
    \end{align*}
    To conclude the proof, simply use the asymptotic value of $t_x$ obtained above to find that
    \begin{align*}
        \lambda_x\big((1+\epsilon)t_x\big)-\lambda_x(t_x)&=\left(\frac{\pi x}{2}\right)^{2/3}\bigg[2(1+\epsilon)\left[1+O\left(x^{2/3}\right)\right]+\left[1+O\left(x^{2/3}\right)\right]\bigg]\frac{\epsilon^2}{(1+\epsilon)^2}+O\left(\epsilon^2x^{4/3}\right)\\
        &=3\left(\frac{\pi x}{2}\right)^{2/3}\epsilon^2+O\left(\epsilon^3x^{2/3}+\epsilon^2x^{4/3}\right)\,.
    \end{align*}
    The uniformity of the $O(\cdot)$ terms directly comes from the fact that all terms depending on $x$ are used within the $O(\cdot)$ and not as constants.
\end{proof}

\section{Asymptotic of the partition function: part I}\label{sec:part fun 1}

The goal of this section is to show the following theorem, stating the asymptotic behaviour of the term $W_n$ introduced in~\eqref{eq:def W}.

\begin{theorem}\label{thm:W asympt}
    Let $\mu=\mu_n$ be such that $1/\sqrt{n}\ll\mu\ll n^{1/4}$.
    Then, as $n$ diverges to infinity, the sum $W_n$ defined in~\eqref{eq:def W} satisfies
    \begin{align*}
        W_n=\big(1+o(1)\big)\frac{\pi^{5/6}}{2^{1/3}3^{1/2}}\frac{\mu^{1/3}}{n^{5/6}}\exp\left(-3\left(\frac{\pi^2\mu^2n}{4}\right)^{1/3}\right)\,.
    \end{align*}
\end{theorem}

Due to its proof being rather intricate, we split it into three sub-sections.
We start by focusing on the terms of the sum which end up being the main contributors (in Section~\ref{sec:main W}), then prove that the remaining terms are negligible (in Section~\ref{sec:negligible W}), and finally combine all the different results together to prove that Theorem~\ref{thm:W asympt} holds (in Section~\ref{sec:W asymptotic}).

\subsection{Main contribution}\label{sec:main W}

In this section we focus on the part of the sum that contributes the most in the sum of $W_n$ as defined in~\eqref{eq:def W}.
This corresponds to the terms close to $t_n$ where $t_n$ is the minimizer of the function $\lambda_{\mu/n}$, defined in~\eqref{eq:lx}.
Before proving the desired convergence, we provide a quick lemma regarding the integral approximation of a Gaussian sum.

\begin{lemma}\label{lem:int approx}
    Let $\mu_n,\sigma_n$ be arbitrary (with $\sigma_n>0$) and $\ell_n,r_n\gg1$ as well as $\ell_n/\sigma_n,r_n/\sigma_n\gg1$.
    Then, we can approximate the gaussian integral with the sum from $\mu_n-\ell_n$ to $\mu_n+r_n$ as such
    \begin{align*}
        \left|\sum_{\mu_n-\ell_n<m<\mu_n+r_n}\exp\left(-\frac{(m-\mu_n)^2}{2\sigma_n^2}\right)-\sqrt{2\pi}\sigma_n\right|&\leq1+o(\sigma_n)\,.
    \end{align*}
\end{lemma}

\begin{proof}
    Using that $u\mapsto e^{-u^2/2\sigma_n^2}$ is decreasing on $[0,\infty)$, for any $m\geq\lceil\mu_n\rceil+1\geq\mu_n+1$, the comparison to the integral gives us that
    \begin{align*}
        \exp\left(-\frac{(m-\mu_n)^2}{2\sigma_n^2}\right)&\leq\int_{m-1}^m\exp\left(-\frac{(u-\mu_n)^2}{2\sigma_n^2}\right)du\leq\exp\left(-\frac{(m-\mu_n-1)^2}{2\sigma_n^2}\right)\,.
    \end{align*}
    Thus, by summing on $m$ from $\lceil\mu_n\rceil+1$ to $\lceil\mu_n+r_n\rceil-1$, a telescopic sum appears and leads to
    \begin{align*}
        0&\leq\int_{\lceil\mu_n\rceil}^{\lceil\mu_n+r_n\rceil-1}\exp\left(-\frac{(u-\mu_n)^2}{2\sigma_n^2}\right)du-\sum_{m=\lceil\mu_n\rceil+1}^{\lceil\mu_n+r_n\rceil-1}\exp\left(-\frac{(m-\mu_n)^2}{2\sigma_n^2}\right)\\
        &\leq\exp\left(-\frac{(\lceil\mu_n\rceil-\mu_n)^2}{2\sigma_n^2}\right)-\exp\left(-\frac{(\lceil\mu_n+r_n\rceil-\mu_n-1)^2}{2\sigma_n^2}\right)\\
        &\leq\exp\left(-\frac{(\lceil\mu_n\rceil-\mu_n)^2}{2\sigma_n^2}\right)\,.
    \end{align*}
    Observe that by subtracting the last term everywhere and reversing all inequality by multiplying them with $-1$, this is the same as
    \begin{align*}
        0&\leq\sum_{m=\lceil\mu_n\rceil}^{\lceil\mu_n+r_n\rceil-1}\exp\left(-\frac{(m-\mu_n)^2}{2\sigma_n^2}\right)-\int_{\lceil\mu_n\rceil}^{\lceil\mu_n+r_n\rceil-1}\exp\left(-\frac{(u-\mu_n)^2}{2\sigma_n^2}\right)du\\
        &\leq\exp\left(-\frac{(\lceil\mu_n\rceil-\mu_n)^2}{2\sigma_n^2}\right)
    \end{align*}
    In a similar manner, using that $u\mapsto e^{-u^2/2\sigma_n^2}$ is increasing on $(-\infty,0]$ and by summing over $m$ from $\lfloor\mu_n-\ell_n\rfloor+1$ to $\lceil\mu_n\rceil-1$, we obtain
    \begin{align*}
        0&\leq\sum_{m=\lfloor\mu_n-\ell_n\rfloor+1}^{\lceil\mu_n\rceil-1}\exp\left(-\frac{(m-\mu_n)^2}{2\sigma_n^2}\right)-\int_{\lfloor\mu_n-\ell_n\rfloor}^{\lceil\mu_n\rceil-1}\exp\left(-\frac{(u-\mu_n)^2}{2\sigma_n^2}\right)du\\
        &\leq\exp\left(-\frac{(\lceil\mu_n\rceil-\mu_n-1)^2}{2\sigma_n^2}\right)-\exp\left(-\frac{(\lfloor\mu_n-\ell_n\rfloor-\mu_n)^2}{2\sigma_n^2}\right)\\
        &\leq\exp\left(-\frac{(\lceil\mu_n\rceil-\mu_n-1)^2}{2\sigma_n^2}\right)\,.
    \end{align*}
    Before concluding, observe that, as $m$ is an integer, we have that $m<\mu_n+r_n$ is equivalent to $m\leq\lceil\mu_n+r_n\rceil-1$ and $m>\mu_n-\ell_n$ is equivalent to $m\geq\lfloor\mu_n-\ell_n\rfloor+1$.
    Thus, if we combine the previous two displayed inequalities we obtain on one hand that
    \begin{align*}
        \sum_{\mu_n-\ell_n<m<\mu_n+r_n}\exp\left(-\frac{(m-\mu_n)^2}{2\sigma_n^2}\right)-\int_{\lfloor\mu_n-\ell_n\rfloor}^{\lceil\mu_n+r_n\rceil-1}\exp\left(-\frac{(u-\mu_n)^2}{2\sigma_n^2}\right)du
        &\geq-\int_{\lceil\mu_n\rceil-1}^{\lceil\mu_n\rceil}e^{-(u-\mu_n)^2/2\sigma_n^2}du\,.
    \end{align*}
    Observing that the exponential term is less than $1$, the lower bound is itself larger than $-1$.
    On the other hand, we have that
    \begin{align*}
        &\hspace{-0.5cm}\sum_{\mu_n-\ell_n<m<\mu_n+r_n}\exp\left(-\frac{(m-\mu_n)^2}{2\sigma_n^2}\right)-\int_{\lfloor\mu_n-\ell_n\rfloor}^{\lceil\mu_n+r_n\rceil-1}\exp\left(-\frac{(u-\mu_n)^2}{2\sigma_n^2}\right)du\\
        &\leq\exp\left(-\frac{(\lceil\mu_n\rceil-\mu_n)^2}{2\sigma_n^2}\right)+\exp\left(-\frac{(\lceil\mu_n\rceil-\mu_n-1)^2}{2\sigma_n^2}\right)-\int_{\lceil\mu_n\rceil-1}^{\lceil\mu_n\rceil}e^{-(u-\mu_n)^2/2\sigma_n^2}du\,,
    \end{align*}
    and by observing that the integral is larger than the minimum of the two other values and using that the exponential terms are less than $1$, we can replace the upper bound by $1$.
    The two previous results thus imply that
    \begin{align*}
        \left|\sum_{\mu_n-\ell_n<m<\mu_n+r_n}\exp\left(-\frac{(m-\mu_n)^2}{2\sigma_n^2}\right)-\int_{\lfloor\mu_n-\ell_n\rfloor}^{\lceil\mu_n+r_n\rceil-1}\exp\left(-\frac{(u-\mu_n)^2}{2\sigma_n^2}\right)du\right|\leq1\,.
    \end{align*}
    To conclude the proof, simply recall $\ell_n/\sigma_n$ and $r_n/\sigma_n$ both diverge to see that
    \begin{align*}
        \int_{\lfloor\mu_n-\ell_n\rfloor}^{\lceil\mu_n+r_n\rceil-1}\exp\left(-\frac{(u-\mu_n)^2}{2\sigma_n^2}\right)du&=\sigma_n\int_\frac{\lfloor \mu_n-\ell_n\rfloor-\mu_n}{\sigma_n}^\frac{\lceil\mu_n+r_n\rceil-1-\mu_n}{\sigma_n}\exp\left(-\frac{u^2}{2}\right)du\\
        &=\sigma_n\Big(\sqrt{2\pi}+o(1)\Big)\,,
    \end{align*}
    leading to the desired result.
\end{proof}

With this lemma and the results from Section~\ref{sec:on part fun} regarding the exponent function $\lambda_x$, we now have all the tools to prove the asymptotic behaviour of the main terms of the sum in~\eqref{eq:def W}.

\begin{proposition}\label{prop:W main}
    {Recall} $t_{\mu/n}$ is the minimizer of $\lambda_{\mu/n}$.
    Let $\epsilon_n$ be such that $\epsilon_n\ll1$ and $\epsilon_nt_n\gg1$.
    Further assume that $\epsilon_n(\mu^2n)^{1/6}\gg1$.
    Then, under these assumptions, we have that
    \begin{align*}
        \frac{1}{4^n}\sum_{|m-t_{\mu/n}|<\epsilon_nt_{\mu/n}}H_{n,m-1}e^{-\mu m}=\big(1+o(1)\big)\frac{\pi^{5/6}}{2^{1/3}3^{1/2}}\frac{\mu^{1/3}}{n^{5/6}}\exp\left(-3\left(\frac{\pi^2\mu^2n}{4}\right)^{1/3}\right)\,.
    \end{align*}
\end{proposition}

\begin{proof}
    We start by recalling that we assume here that $1/\sqrt{n}\ll\mu\ll n^{1/4}$.
    Moreover, we observe that this implies that $\mu/n\ll1$ and so we can apply the results of Proposition~\ref{prop:asympt lx} to $x=\mu/n$.
    In that scenario, and since $t_{\mu/n}$ is of order $(n/\mu)^{1/3}$ and that $\epsilon_n\ll1$, we can see that $1\ll(1-\epsilon)t_n\leq(1+\epsilon_n)t_{\mu/n}\ll\sqrt{n}$.
    Thus we can uniformly apply the approximation of $H_{n,m-1}$ from Theorem~\ref{thm:ht-bdd tree} and see that
    \begin{align*}
        \frac{1}{4^n}\sum_{|m-t_{\mu/n}|<\epsilon_nt_{\mu/n}}H_{n,m-1}e^{-\mu m}&=\big(1+o(1)\big)\frac{1}{4^n}\sum_m\frac{4^n}{m}\frac{\tan^2\left(\frac{\pi}{m}\right)}{\left(1+\tan^2\left(\frac{\pi}{m}\right)\right)^n}e^{-\mu m}\\
        &=\big(1+o(1)\big)\frac{\pi^2}{t_n^3}\sum_me^{-n\lambda_{\mu/n}(m)}\,,
    \end{align*}
    where we used the fact that $\epsilon_n\ll1$ and $t_n\gg1$ to remove the polynomial terms out of the sum.
    Now, for any sequence $m=m_n$ between $(1-\epsilon_n)t_{\mu/n}$ and $(1+\epsilon_n)t_{\mu/n}$, the sequence $\gamma_n=m/t_{\mu/n}-1\ll1$ and we can apply Proposition~\ref{prop:asympt lx}~\numitem{3} to obtain that
    \begin{align*}
        e^{-n\lambda_{\mu/n}(m)}&=e^{-n\lambda_{\mu/n}(t_{\mu/n})}\exp\left(-n\big(3+o(1)\big)\left(\frac{\pi\mu}{2n}\right)^{2/3}\gamma_n^2\right)\,,
    \end{align*}
    where the $o(\cdot)$ term is uniform over such $m$.
    Thus, plugging this formula back into the previous sum and replacing $\gamma_n$ by its value in relation to $m$, we obtain
    \begin{align*}
        \frac{1}{4^n}\sum_{|m-t_{\mu/n}|<\epsilon_nt_{\mu/n}}H_{n,m-1}e^{-\mu m}&=\big(1+o(1)\big)\frac{\pi^2}{t_{\mu/n}^3}e^{-n\lambda_{\mu/n}(t_{\mu/n})}\sum_me^{-n\frac{3+o(1)}{t_{\mu/n}^2}\left(\frac{\pi\mu}{2n}\right)^{2/3}(m-t_{\mu/n})^2}\,.
    \end{align*}
    Consider now the statement of Lemma~\ref{lem:int approx} with $\mu=t_{\mu/n}$, $2\sigma_n^2=\frac{1}{n}\frac{t_{\mu/n}^2}{3+o(1)}\left(\frac{2n}{\pi\mu}\right)^{2/3}$, and $\ell_n,r_n=\epsilon_nt_{\mu/n}$.
    First of all, thanks to the asymptotic behaviour of $t_{\mu/n}$ following from Proposition~\ref{prop:asympt lx}~\numitem{1}, we establish that
    \begin{align*}
        \sigma_n^2=\big(1+o(1)\big)\frac{1}{n}\left(\frac{2\pi^2n}{\mu}\right)^{2/3}\frac{1}{6}\left(\frac{2n}{\pi\mu}\right)^{2/3}=\big(1+o(1)\big)\frac{1}{3}\left(\frac{2\pi^2n}{\mu^4}\right)^{1/3}\,.
    \end{align*}
    It is worth noting that the assumption that $\mu\ll n^{1/4}$ implies that $\sigma_n\rightarrow\infty$.
    Now, we already know that $\ell_n=r_n\gg1$ since $\epsilon_nt_{\mu/n}\gg1$.
    For the second assumption, observe that
    \begin{align*}
        \frac{\ell_n}{\sigma_n}&=\big(1+o(1)\big)\epsilon_n\left(\frac{2\pi^2n}{\mu}\right)^{1/3}\sqrt{3}\left(\frac{\mu^4}{2\pi^2n}\right)^{1/6}=\big(C+o(1)\big)\epsilon_n(\mu^2n)^{1/6}\,,
    \end{align*}
    which diverges once again thanks to the assumption on $\epsilon_n$.
    This means that we can apply Lemma~\ref{lem:int approx} to the sum above and, further using that $\sigma_n\rightarrow\infty$, obtain
    \begin{align*}
        \frac{1}{4^n}\sum_{|m-t_{\mu/n}|<\epsilon_nt_{\mu/n}}H_{n,m-1}e^{-\mu m}&=\big(1+o(1)\big)\frac{\pi^2}{t_{\mu/n}^3}e^{-n\lambda_{\mu/n}(t_{\mu/n})}\sqrt{2\pi}\sigma_n\,.
    \end{align*}
    To conclude this proof recall the asymptotic behaviour of $t_{\mu/n}$ and $\sigma_n$ as well as that of $\lambda_{\mu/n}(t_{\mu/n})$ from Proposition~\ref{prop:asympt lx}~\numitem{1} and~\numitem{2} to see that
    \begin{align*}
        \frac{1}{4^n}\sum_{|m-t_{\mu/n}|<\epsilon_nt_{\mu/n}}H_{n,m-1}e^{-\mu m}&=\big(1+o(1)\big)\frac{\pi^2}{\frac{2\pi^2n}{\mu}}e^{-3n\left(\frac{\pi\mu}{2n}\right)^{2/3}+O\left(n\left(\frac{\mu}{n}\right)^{4/3}\right)}\sqrt{2\pi}\frac{1}{\sqrt{3}}\left(\frac{2\pi^2n}{\mu^4}\right)^{1/6}\\
        &=\big(1+o(1)\big)\frac{\pi^{5/6}}{2^{1/3}3^{1/2}}\frac{\mu^{1/3}}{n^{5/6}}e^{-3\left(\frac{\pi^2\mu^2n}{4}\right)^{1/3}}\,,
    \end{align*}
    where we used that $O(n(\mu/n)^{4/3})=O((\mu^4/n)^{1/3})=o(1)$.
    This result is exactly the claimed asymptotic behaviour.
\end{proof}

\subsection{Negligible terms}\label{sec:negligible W}

\def\smt{X}
We now prove that Theorem~\ref{thm:W asympt} holds by showing that the terms away from $t_{\mu/n}$ in the sum~\eqref{eq:def W} are negligible, where $t_{\mu/n}$ is the location of the minimum of the function $\lambda_{\mu/n}$ from~\eqref{eq:lx}. Before doing so, we consider the sequence of terms defined by
\begin{align}\label{eq:def smt}
    \smt_{n,m}=\frac{1}{4^n}\frac{n^{5/6}}{\mu^{1/3}}H_{n,m-1}\exp\left(-\mu m+3\left(\frac{\pi^2\mu^2n}{4}\right)^{1/3}\right)\,,
\end{align}
which exactly correspond the terms of the definition of $W_n$ in~\eqref{eq:def W} divided by the desired asymptotic behaviour of $W_n$ from Theorem~\ref{thm:W asympt} (up to constant terms).
Thus, thanks to Proposition~\ref{prop:W main} which shows that the main contribution of the terms in $W_n$ has the desired asymptotic behaviour, it suffices to show that the sum of the $\smt_{n,m}$ when $|m-t_{\mu/n}|\geq\epsilon_nt_{\mu/n}$ converges to $0$ for some appropriate $\epsilon_n$.
We divide this proof in two parts, depending on the absolute value of $m-t_{\mu/n}$.
We start with the easiest one of the two, corresponding to $m-t_{\mu/n}>0$.

\begin{proposition}\label{prop:negligible W+}
    {Recall $t_{\mu/n}=\argmin_{t \in (0,2)}\lambda_{\mu/n}(t)$, and} $\smt_{n,m}$ be as defined as in~\eqref{eq:def smt}, and $\epsilon_n\ll1$ be such that $\epsilon_nt_n\gg1$ and $\epsilon_n^2(\mu^2n)^{1/3}\gg\log(\mu^2n)$.
    Then, the sum over $m\geq(1+\epsilon_n)t_{\mu/n}$ of $\smt_{n,m}$ is negligible, as given by
    \begin{align*}
        \lim_{n\rightarrow\infty}\sum_{(1+\epsilon_n)t_{\mu/n}\leq m\leq n}\smt_{n,m}=0\,.
    \end{align*}
\end{proposition}

\begin{proof}
    We prove this bound by splitting the sum in two parts.
    For the highest values of $m$, namely $m>ct_{\mu/n}$ for some fixed $c>0$, we crudely use that $H_{n,m-1}\leq C_{n-1}$ where $C_{n-1}$ is the $(n-1)$-th Catalan number.
    This implies that
    \begin{align*}
        \sum_{m\geq ct_{\mu/n}}\smt_{n,m}&\leq\frac{C_{n-1}}{4^n}\frac{n^{5/6}}{\mu^{1/3}}\exp\left(3\left(\frac{\pi^2\mu^2n}{4}\right)^{1/3}\right)\sum_{m\geq ct_{\mu/n}}e^{-\mu m}\\
        &\leq\frac{C_{n-1}}{4^n}\frac{n^{5/6}}{\mu^{1/3}}\frac{1}{1-e^{-\mu}}\exp\left(3\left(\frac{\pi^2\mu^2n}{4}\right)^{1/3}-\mu ct_{\mu/n}\right)\,.
    \end{align*}
    We now recall that $C_{n-1}/4^n\sim 1/(4\sqrt{\pi}n^{3/2})$ and that Proposition~\ref{prop:asympt lx}~\numitem{1} tells us that
    \begin{align*}
        t_{\mu/n}=\left(\frac{2\pi^2n}{\mu}\right)^{1/3}+O\left(\left(\frac{\mu}{n}\right)^{1/3}\right)\,.
    \end{align*}
    Moreover, using that $e^\mu\geq1+\mu$, we can check that $1/(1-e^{-\mu})\leq(1+\mu)/\mu$.
    Finally, combining the previous remarks, it follows that
    \begin{align*}
        \sum_{m\geq ct_{\mu/n}}\smt_{n,m}&\leq\frac{1}{4\sqrt{\pi}}\frac{(1+\mu)}{n^{2/3}\mu^{4/3}}\exp\left(3\left(\frac{\pi^2\mu^2n}{4}\right)^{1/3}-c\big(2\pi^2\mu^2n\big)^{1/3}+O\left(\left(\frac{\mu^4}{n}\right)^{1/3}\right)\right)\\
        &=\frac{1+o(1)}{4\sqrt{\pi}}\frac{(1+\mu)}{(\mu^2n)^{2/3}}\exp\left(-(2c-3)\left(\frac{\pi^2\mu^2n}{4}\right)^{1/3}\right)\,.
    \end{align*}
    Thus, for any $2c-3>0$, the right-hand side converges to $0$ as desired.
    Without loss of generality, we now fix $c=2$.

    For the rest of the sum, we aim to bound $\smt_{n,m}$ when $(1+\epsilon_n)t_{\mu/n}\leq m<2t_{\mu/n}$.
    Observe that $2t_{\mu/n}\ll\sqrt{n}$ so that the result from Theorem~\ref{thm:ht-bdd tree} applies and we have
    \begin{align*}
        \sum_{(1+\epsilon_n)t_{\mu/n}\leq m<2t_{\mu/n}}\smt_{n,m}&=\frac{1+o(1)}{4^n}\frac{n^{5/6}}{\mu^{1/3}}\sum_m\frac{4^n}{m}\frac{\tan^2\left(\frac{\pi}{m}\right)}{\left(1+\tan^2\left(\frac{\pi}{m}\right)\right)^n}\exp\left(-\mu m+3\left(\frac{\pi^2\mu^2n}{4}\right)^{1/3}\right)\,.
    \end{align*}
    Now recalling the definition of $\lambda_x$ from~\eqref{eq:lx} and using that $\tan^2\pi/m\sim\pi^2/m^2$ uniformly for all $m$ {in the limit of} the sum, we see that the previous equation re-writes as
    \begin{align*}
        \sum_{(1+\epsilon_n)t_{\mu/n}\leq m<2t_{\mu/n}}\smt_{n,m}&=\big(1+o(1)\big)\frac{n^{5/6}}{\mu^{1/3}}\sum_m\frac{\pi^2}{m^3}\exp\left(-n\lambda_{\mu/n}(m)+3\left(\frac{\pi^2\mu^2n}{4}\right)^{1/3}\right)\,.
    \end{align*}
    Using that $m\mapsto e^{-n\lambda_{\mu/n}(m)}$ and $m\mapsto1/m^3$ are both decreasing when $m\geq t_{\mu/n}$, it follows that
    \begin{align*}
        &\sum_{(1+\epsilon_n)t_{\mu/n}\leq m<2t_{\mu/n}}\smt_{n,m}\\
        &\hspace{0.5cm}\leq\big(1+o(1)\big)\frac{n^{5/6}}{\mu^{1/3}}\big(2t_{\mu/n}-(1+\epsilon_n)t_{\mu/n}\big)\frac{\pi^2}{((1+\epsilon_n)t_{\mu/n})^3}\exp\left(-n\lambda_{\mu/n}\big((1+\epsilon_n)t_{\mu/n}\big)+3\left(\frac{\pi^2\mu^2n}{4}\right)^{1/3}\right)\,.
    \end{align*}
    To conclude this proof, first recall that $\epsilon_n\ll1$ and $t_{\mu/n}\sim(2\pi^2n/\mu)^{1/3}$, along with the asymptotic behaviour of $\lambda_{\mu/n}(t_{\mu/n})$ and $\lambda_{\mu/n}((1+\epsilon_n)t_{\mu/n})$ from Proposition~\ref{prop:asympt lx}~\numitem{2} and~\numitem{3}, to see that
    \begin{align*}
        \sum_{(1+\epsilon_n)t_{\mu/n}\leq m<2t_{\mu/n}}\smt_{n,m}&\leq\big(1+o(1)\big)\frac{n^{5/6}}{\mu^{1/3}}\frac{\pi^2}{\left(\frac{2\pi^2n}{\mu}\right)^{2/3}}\exp\left(-\big(3+o(1)\big)\left(\frac{\pi^2\mu^2n}{4}\right)^{1/3}\epsilon_n^2\right)\,.
    \end{align*}
    To conclude the proof, simply note that
    \begin{align*}
        \frac{n^{5/6}}{\mu^{1/3}}\frac{\pi^2}{\left(\frac{2\pi^2n}{\mu}\right)^{2/3}}=\left(\frac{\pi^4\mu^2n}{16}\right)^{1/6}
    \end{align*}
    and recall that we assume that $\epsilon_n^2(\mu^2n)^{1/3}\gg\log(\mu^2n)$, which is exactly the desired assumption to prove the convergence to $0$ of the above right-hand side.
    This proves that the second part of the sum over $m$ also converges to $0$, thus proving the statement of the lemma.
\end{proof}

\begin{proposition}\label{prop:negligible W-}
    {Recall $t_{\mu/n}=\argmin_{t \in (0,2)}\lambda_{\mu/n}(t)$, and} let $\smt_{n,m}$ be defined as in~\eqref{eq:def smt}, and $\epsilon_n\ll1$ be such that $\epsilon_nt_{\mu/n}\gg1$ and $\epsilon_n^2(\mu^2n)^{1/3}\gg\log(\mu^2n)$.
    Then, the sum over $m\leq(1-\epsilon_n)t_{\mu/n}$ of $\smt_{n,m}$ is negligible, i.e.,
    \begin{align*}
        \lim_{n\rightarrow\infty}\sum_{3\leq m\leq(1-\epsilon_n)t_{\mu/n}}\smt_{n,m}=0\,.
    \end{align*}
\end{proposition}

\begin{proof}
    Similar to the proof of Proposition~\ref{prop:negligible W+}, we need to split the sum in two parts in order to prove the desired convergence.
    Fix for now $\delta\in(0,1)$ arbitrary.
    For the values of $m$ closest to $t_{\mu/n}$, use the asymptotic behaviour of $H_{n,m}$ from Theorem~\ref{thm:ht-bdd tree}, along with the divergence of $m$ and the definition of $\lambda_{\mu/n}$ from~\eqref{eq:lx} to see that
    \begin{align*}
        \sum_{\delta t_{\mu/n}\leq m\leq(1-\epsilon_n)t_{\mu/n}}\smt_{n,m}&=\big(1+o(1)\big)\frac{n^{5/6}}{\mu^{1/3}}\exp\left(3\left(\frac{\pi^2\mu^2n}{4}\right)^{1/3}\right)\sum_m\frac{\pi^2}{m^3}e^{-n\lambda_{\mu/n}(m)}\,.
    \end{align*}
    Further using the fact that $m\mapsto\lambda_{\mu/n}(m)$ and $m\mapsto1/m^3$ are decreasing on the set of $m$ considered here, we can bound the terms of the sum by taking the extreme values and obtain
  
    \begin{align*}
        \sum_{\delta t_{\mu/n}\leq m\leq(1-\epsilon_n)t_{\mu/n}}\smt_{n,m}&\leq\big(1+o(1)\big)\frac{n^{5/6}}{\mu^{1/3}}\frac{\pi^2}{\delta^3t_n^2}\exp\left(3\left(\frac{\pi^2\mu^2n}{4}\right)^{1/3}-n\lambda_{\mu/n}\big((1-\epsilon_n)t_{\mu/n}\big)\right)\,.
    \end{align*}
    To conclude this first part of the proof, simply recall the asymptotic behaviour of $t_{\mu/n}$ and of $\lambda_{\mu/n}$ from Proposition~\ref{prop:asympt lx}~\numitem{1},~\numitem{2} and~\numitem{3} to obtain
    \begin{align*}
        \sum_{\delta t_{\mu/n}\leq m\leq(1-\epsilon_n)t_{\mu/n}}\smt_{n,m}&\leq\big(1+o(1)\big)\frac{n^{5/6}}{\mu^{1/3}}\frac{\pi^2}{\delta^3\left(\frac{2\pi^2n}{\mu}\right)^{2/3}}\exp\left(-\big(3+o(1)\big)\left(\frac{\pi^2\mu^2n}{4}\right)^{1/3}\epsilon_n^2\right)\\
        &\leq\big(1+o(1)\big)\frac{(\pi^4\mu^2n)^{1/6}}{\delta^3}\exp\left(-\big(3+o(1)\big)\left(\frac{\pi^2\mu^2n}{4}\right)^{1/3}\epsilon_n^2\right)\,.
    \end{align*}
    Thus, thanks to the assumption that $\epsilon_n^2(\mu^2n)^{1/3}\gg\log(\mu^2n)$, the right-hand side converges to $0$ as desired.
    We now prove that the second part of the sum converges to $0$ as well.

    Start by using the result from Theorem~\ref{thm:ht-bdd tree} along with the definition of $\smt_{n,m}$ from~\eqref{eq:def smt} to see that, for any $3\leq m\leq t_{\mu/n}\ll\sqrt{n}$ we have
    \begin{align*}
        \smt_{n,m}\leq\big(1+o(1)\big)\frac{n^{5/6}}{\mu^{1/3}}\frac{\tan^2\left(\frac{\pi}{m}\right)}{m}\exp\left(3\left(\frac{\pi^2\mu^2n}{4}\right)^{1/3}-n\lambda_{\mu/n}(m)\right)\,,
    \end{align*}
    where the $o(\cdot)$ term does not depend on $m$.
    Moreover, we observe that $m^2\tan^2(\pi/m)$ is bounded since $m\geq3$.
    Finally, recalling the definition of $\lambda_{\mu/n}$ from~\eqref{eq:lx}, we see that for any $m\geq3$, we have
    \begin{align*}
        m^2\lambda_{\mu/n}(m)\geq m^2\log\left(1+\tan^2\left(\frac{\pi}{m}\right)\right)\,,
    \end{align*}
    and the right hand side is bounded away from $0$.
    Thus, there exists $C>0$ such that, for any $n$ large enough and $3\leq m\leq t_{\mu/n}$ we have
    \begin{align*}
        \smt_{n,m}\leq C\frac{(2\pi^2)^{2/3}n^{5/6}}{m^3\mu^{1/3}}\exp\left(3\left(\frac{\pi^2\mu^2n}{4}\right)^{1/3}-\frac{1}{C}\frac{\pi^2n}{m^2}\right)\,,
    \end{align*}
    It follows that, for any $x\leq y$ such that $3\leq xt_{\mu/n}\leq yt_{\mu/n}\leq  t_{\mu/n}$ we have
    \begin{align*}
        \sum_{xt_{\mu/n}\leq m<yt_{\mu/n}}\smt_{n,m}&\leq\sum_{xt_{\mu/n}\leq m<yt_{\mu/n}}C\frac{(2\pi^2)^{2/3}n^{5/6}}{m^3\mu^{1/3}}\exp\left(3\left(\frac{\pi^2\mu^2n}{4}\right)^{1/3}-\frac{1}{C}\frac{\pi^2n}{m^2}\right)\\
        &\leq(y-x)t_{\mu/n}\cdot C\frac{(2\pi^2)^{2/3}n^{5/6}}{(xt_n)^3\mu^{1/3}}\exp\left(3\left(\frac{\pi^2\mu^2n}{4}\right)^{1/3}-\frac{1}{C}\frac{\pi^2n}{(yt_{\mu/n})^2}\right)\,.
    \end{align*}
    Note that $y\leq 1$ and use the asymptotic behaviour of $t_{\mu/n}$ from Proposition~\ref{prop:asympt lx}~\numitem{1} to see that
    \begin{align*}
        \sum_{xt_{\mu/n}\leq m<yt_{\mu/n}}\smt_{n,m}&\leq C\frac{(\mu^2n)^{1/6}}{x^3}\exp\left(-\left(\frac{1}{Cy^2}-3\right)\left(\frac{\pi^2\mu^2n}{4}\right)^{1/3}\right)\,.
    \end{align*}
    But now we observe that this upper bound is exactly of the same form as $F_t$ from Corollary~\ref{cor:useful-bound} where we match $t$ with $(\mu^2n)^{1/3}$ and have $\alpha=1/2$, $B(x)=x^3/C$, and $\Gamma(y)=(1/Cy^2-3)^{-1}$.
    Observe that $\Gamma$ might diverge when $Cy^2$ approaches $1/3$.
    However, by setting $a_+<1/\sqrt{3C}$, we see that all hypothesis from Corollary~\ref{cor:useful-bound} apply and it follows that we can find $3/t_{\mu/n}=a_0\leq\cdots\leq a_k=a_+$ such that
    \begin{align*}
        \sum_{i=1}^kC\frac{(\mu^2n)^{1/6}}{a_{i-1}^3}\exp\left(-\left(\frac{1}{Ca_i^2}-3\right)\left(\frac{\pi^2\mu^2n}{4}\right)^{1/3}\right)\longrightarrow0\,.
    \end{align*}
    This then implies that
    \begin{align*}
        \sum_{3\leq m<a_+t_{\mu/n}}\smt_{n,m}\longrightarrow0\,.
    \end{align*}
    To conclude the proof, simply recall that the sum from $\delta t_{\mu/n}$ to $(1-\epsilon_n)t_{\mu/n}$ converges to $0$ as well, for any arbitrary $\delta>0$.
    Thus, by setting $\delta=a_+$, the sum over $m$ from $3$ to $(1-\epsilon_n)t_{\mu/n}$ converges to $0$, as desired.
\end{proof}

\subsection{Combining all results}\label{sec:W asymptotic}

In here, we combine all the previous results in order to prove Theorem~\ref{thm:W asympt}.

\begin{proof}[Proof of Theorem~\ref{thm:W asympt}]
    Before starting the proof, we recall that we assume that $1/\sqrt{n}\ll\mu\ll n^{1/4}$.
    Under these assumptions, observe that it is possible to find $\epsilon_n\ll1$ such that $\epsilon_nt_{\mu/n}\gg1$, $\epsilon_n(\mu^2n)^{1/6}\gg1$, and $\epsilon_n^2(\mu^2n)^{1/3}\gg\log(\mu^2n)$.
    Indeed, in that case, since $t_{\mu/n}$ is of order $(n/\mu)^{1/3}$ it diverges to $\infty$ so it is possible to have $\epsilon_nt_{\mu/n}\gg1$.
    Moreover, the second condition is actually a consequence of the third one which can also be satisfied since $\log(\mu^2n)/(\mu^2n)^{1/3}\ll1$.
    We now choose such an $\epsilon_n$ which, by definition, satisfies the conditions from Proposition~\ref{prop:W main},~\ref{prop:negligible W+} and~\ref{prop:negligible W-}.

    Recall now the definition of $W_n$ from~\eqref{eq:def W} and see that Proposition~\ref{prop:W main} states that
    \begin{align*}
        \frac{1}{4^n}\sum_{|m-t_{\mu/n}|<\epsilon_nt_{\mu/n}}H_{n,m-1}e^{-\mu m}=\big(1+o(1)\big)\frac{\pi^{5/6}}{2^{1/3}3^{1/2}}\frac{\mu^{1/3}}{n^{5/6}}\exp\left(-3\left(\frac{\pi^2\mu^2n}{4}\right)^{1/3}\right)\,.
    \end{align*}
    Moreover, using the definition of $\smt_{n,m}$ from~\eqref{eq:def smt} along with Proposition~\ref{prop:negligible W+} and~\ref{prop:negligible W-}, we see that
    \begin{align*}
        \frac{1}{4^n}\sum_{|m-t_{\mu/n}|\geq\epsilon_nt_{\mu/n}}H_{n,m-1}e^{-\mu m}=o\left(\frac{\mu^{1/3}}{n^{5/6}}\exp\left(-3\left(\frac{\pi^2\mu^2n}{4}\right)^{1/3}\right)\right)\,.
    \end{align*}
    Combining the previous two results directly implies that
    \begin{align*}
        W_n=\frac{1}{4^n}\sum_{m=3}^nH_{n,m-1}e^{-\mu m}=\big(1+o(1)\big)\frac{\pi^{5/6}}{2^{1/3}3^{1/2}}\frac{\mu^{1/3}}{n^{5/6}}\exp\left(-3\left(\frac{\pi^2\mu^2n}{4}\right)^{1/3}\right)\,,
    \end{align*}
    exactly as desired.
\end{proof}

\section{Asymptotic of the partition function: part II}\label{sec:part fun 2}

The goal of this section is to extend the results from Theorem~\ref{thm:W asympt} regarding the asymptotic of $W_n$ as defined in~\eqref{eq:def W}.
We actually start with a result overlapping with Theorem~\ref{thm:W asympt}, thus showing that the difficult part of the proof of the theorem is actually when $\mu$ is ``close'' to $1/\sqrt{n}$, or more precisely when $\mu$ is at most of order $(\log^{3/2}n)/\sqrt{n}$.
We then provide a second result, overlapping with the first one and covering the case $\mu\gg n^{1/4}$, for which we show that the sum of $W_n$ can be reduced to at most two terms.
The two sections below cover these two cases.

\subsection{Generic behaviour}

In the case where $\mu$ is in the intermediate range but not too close to $1/\sqrt{n}$, the proof of Theorem~\ref{thm:W asympt} can actually be simplified and does not require bounds as subtle as to prove the case close to $1/\sqrt{n}$.
We provide here the complete proof of the asymptotic of $W_n$ in that regime.

\begin{proposition}\label{prop:W as a sum}
    Let $\mu=\mu_n$ be such that $(\log^{3/2}n)/\sqrt{n}\ll\mu\ll n$ and {recall $t_{\mu/n}=\argmin_{t \in (0,2)}\lambda_{\mu/n}(t)$.}
    Then, for any $\epsilon_n\ll1$ such that $(\mu^2n)^{1/3}\epsilon_n^2\gg\log n$ and $\epsilon\gg(\mu/n)^{1/3}$ (which are both possible thanks to the assumption on $\mu$) the sum $W_n$ defined in~\eqref{eq:def W} satisfies
    \begin{align*}
        W_n=\big(1+o(1)\big)\frac{\mu}{2n}e^{-n\lambda_{\mu/n}(t_{\mu/n})}\sum_{|m-t_{\mu/n}|<\epsilon_nt_{\mu/n}}\exp\left(-\left(\frac{3}{2}+o(1)\right)\left(\frac{\mu^4}{2\pi^2n}\right)^{1/3}\big(m-t_{\mu/n}\big)^2\right)\,,
    \end{align*}
    where the $o(\cdot)$ in the exponential are uniform over all terms of the sum.
\end{proposition}

\begin{proof}
    To prove this result we split the sum into four parts.
    We start with the main contribution of the sum, which almost directly follows previous results.
    Indeed, for any $\epsilon_n\ll1$, Theorem~\ref{thm:ht-bdd tree} tells us that
    \begin{align*}
        \frac{1}{4^n}\sum_{|m-t_{\mu/n}|<\epsilon_nt_{\mu/n}}H_{n,m-1}e^{-\mu m}=\big(1+o(1)\big)\frac{\pi^2}{t_{\mu/n}^3}\sum_{|m-t_{\mu/n}|\leq\epsilon_nt_{\mu/n}}e^{-n\lambda_{\mu/n}(m)}\,,
    \end{align*}
    where $\lambda_{\mu/n}$ is defined in~\eqref{eq:lx}.
    But then using the asymptotic behaviour of $t_{\mu/n}$ along with the behaviour of $\lambda_{\mu/n}$ close to $t_{\mu/n}$ with $\gamma=(m/t_{\mu/n}-1)$ from Proposition~\ref{prop:asympt lx}~\numitem{1} and~\numitem{3}, we see that
    \begin{align*}
        \frac{1}{4^n}\sum_{|m-t_{\mu/n}|<\epsilon_nt_{\mu/n}}H_{n,m-1}e^{-\mu m}&=\big(1+o(1)\big)\frac{n}{2\mu}e^{-n\lambda_{\mu/n}(t_{\mu/n})}\sum_{|m-t_{\mu/n}|<\epsilon_nt_{\mu/n}}\exp\left(-\big(3+o(1)\big)\left(\frac{\pi^2\mu^2n}{4}\right)^{1/3}\gamma^2\right)\\
        &=\big(1+o(1)\big)\frac{n}{2\mu}e^{-n\lambda_{\mu/n}(t_{\mu/n})}\sum_{|m-t_{\mu/n}|<\epsilon_nt_{\mu/n}}\exp\left(-\left(\frac{3}{2}+o(1)\right)\left(\frac{\mu^4}{2\pi^2n}\right)^{1/3}(m-t_{\mu/n})^2\right)\,,
    \end{align*}
    which is exactly the desired asymptotic.
    Before focusing on the other parts of the sum, first recall the approximation for $n\lambda_{\mu/n}(t_{\mu/n})$ from Proposition~\ref{prop:asympt lx}:
    \begin{align*}
        n\lambda_{\mu/n}(t_{\mu/n})=3\left(\frac{\pi^2\mu^2n}{4}\right)^{1/3}+O\left(\left(\frac{\mu^4}{n}\right)^{1/3}\right)\,.
    \end{align*}
    We now use this asymptotic behaviour to show that the other parts of the sum are negligible.

    For the first terms of the sum, use Theorem~\ref{thm:ht-bdd tree} to obtain
    \begin{align*}
        \frac{1}{4^n}\sum_{3\leq m\leq(1-\epsilon_n)t_{\mu/n}}H_{n,m-1}e^{-\mu m}&=\big(1+o(1)\big)\sum_{3\leq m\leq(1-\epsilon_n)t_{\mu/n}}\frac{\tan^2\left(\frac{\pi}{m}\right)}{m}e^{-n\lambda_{\mu/n}(m)}\,,
    \end{align*}
    where $\lambda_{\mu/n}$ is defined in~\eqref{eq:lx}.
    Using that $m\mapsto\lambda_{\mu/n}(m)$ is decreasing on $(2,t_{\mu/n})$, we can bound the exponential term by its value at $m=(1-\epsilon_n)t_{\mu/n}$.
    Further bounding the other terms by their maximum, corresponding to $m=3$, we obtain
    \begin{align*}
        \frac{1}{4^n}\sum_{3\leq m\leq(1-\epsilon_n)t_{\mu/n}}H_{n,m-1}e^{-\mu m}&\leq\big(1+o(1)\big)\frac{\tan^2\left(\frac{\pi}{3}\right)}{3}t_{\mu/n}e^{-n\lambda_{\mu/n}((1-\epsilon_n)t_{\mu/n})}\,.
    \end{align*}
    Putting any constant term into a flexible $C$ and dividing by the desired asymptotic leads to
    \begin{align*}
        &\frac{n}{\mu}\exp\left(3\left(\frac{\pi^2\mu^2n}{4}\right)^{1/3}+O\left(\left(\frac{\mu^4}{n}\right)^{1/3}\right)\right)\frac{1}{4^n}\sum_{3\leq m\leq(1-\epsilon_n)t_{\mu/n}}H_{n,m-1}e^{-\mu m}\\
        &\hspace{0.5cm}\leq\big(C+o(1)\big)\left(\frac{n}{\mu}\right)^{4/3}\exp\left(-\big(3+o(1)\big)\left(\frac{\pi^2\mu^2n}{4}\right)^{1/3}\epsilon_n^2+O\left(\left(\frac{\mu^4}{n}\right)^{1/3}\right)\right)\,.
    \end{align*}
    Recall now that we assumed that $(\mu^2n)^{1/3}\epsilon_n^2\gg\log n$ 
    as well as $\epsilon\gg(\mu/n)^{1/3}$ which are exactly the two conditions needed for the first term in the exponential to be the dominant term, proving that this part of the sum is negligible with respect to the desired asymptotic behaviour.
    
    Consider now $r_n\gg1$ such that $r_nt_{\mu/n}\ll\sqrt{n}$ and focus on the sum from $(1+\epsilon_n)t_{\mu/n}$ to $r_nt_{\mu/n}$.
    Applying the same method as for the previous case, we see that
    \begin{align*}
        \frac{1}{4^n}\sum_{(1+\epsilon_n)t_{\mu/n}\leq m<r_nt_{\mu/n}}H_{n,m-1}e^{-\mu m}&\leq\big(1+o(1)\big)\frac{\pi^2r_n}{t_{\mu/n}^2}e^{-n\lambda_{\mu/n}((1+\epsilon_n)t_{\mu/n})}\,,
    \end{align*}
    leading again to
    \begin{align*}
        &\frac{n}{\mu}\exp\left(3\left(\frac{\pi^2\mu^2n}{4}\right)^{1/3}+O\left(\left(\frac{\mu^4}{n}\right)^{1/3}\right)\right)\frac{1}{4^n}\sum_{(1+\epsilon_n)t_{\mu/n}\leq m<r_nt_{\mu/n}}H_{n,m-1}e^{-\mu m}\\
        &\hspace{0.5cm}\leq\big(C+o(1)\big)r_n\left(\frac{n}{\mu}\right)^{1/3} \exp\left(-\big(3+o(1)\big)\left(\frac{\pi^2\mu^2n}{4}\right)^{1/3}\epsilon_n^2+O\left(\left(\frac{\mu^4}{n}\right)^{1/3}\right)\right)\,.
    \end{align*}
    For the same reasons as the previous case, the right-hand side converges to $0$ and this part of the sum is also negligible.

    For the final part of the sum, consider the remainging values for $m$ and upper bound $H_{n,m-1}$ with $C_{n-1}\sim4^n/(4\sqrt{\pi}n^{3/2})$ to obtain
    \begin{align*}
        \frac{1}{4^n}\sum_{r_nt_{\mu/n}\leq m\leq n}H_{n,m-1}e^{-\mu m}&\leq\big(1+o(1)\big)\frac{1}{4\sqrt{\pi}n^{3/2}}\frac{e^{-\mu r_nt_{\mu/n}}}{1-e^{-\mu}}\,.
    \end{align*}
    Using that $\mu t_{\mu/n}$ is of order $(\mu^2n)^{1/3}$, we see that the same arguments as before proves that this part of the sum is also negligible.
    This concludes the proof of the proposition.
\end{proof}

\subsection{Finite sum behaviour}

We conclude this section by covering the case where $\mu\gg n^{1/4}$, for which the sum of $W_n$ can be reduced to only a couple of terms.

\begin{proposition}\label{prop:W mu-n}
    Let $\mu=\mu_n$ be such that $\mu\gg n^{1/4}$. {Recall $t_{\mu/n}=\argmin_{t \in (0,2)}\lambda_{\mu/n}(t)$.}
    Then we have
    \begin{align*}
        W_n=\big(1+o(1)\big)\sum_{m=\max\{\lfloor t_{\mu/n}\rfloor,3\}}^{\lceil t_{\mu/n}\rceil}\frac{\tan^2\frac{\pi}{m}}{m}\exp\left(-\mu m-n\log\left(1+\tan^2\frac{\pi}{m}\right)\right)\,.
    \end{align*}
\end{proposition}

\begin{proof}
    We start by assuming that $\mu/n$ is bounded away from $0$.
    In that case, we fix an arbitrarily large value $M$.
    Using that $H_{n,m-1}\leq C_{n-1}$, we see that
    \begin{align*}
        \frac{1}{4^n}\sum_{m>M}H_{n,m-1}e^{-\mu m}&\leq\big(1+o(1)\big)\frac{e^{-\mu M}}{4\sqrt{\pi}\sqrt{n}}\,.
    \end{align*}
    Moreover, considering the term $m=3$ and observing the definition from~\eqref{eq:ht-bdd tree} implies that $H_{n,2}=1$, we see that
    \begin{align*}
        \frac{1}{4^n}H_{n,2}e^{-3\mu}=\frac{e^{-3\mu}}{4^n}\,,
    \end{align*}
    from which we observe that the sum over $m>M$ is negligible with respect to the first term as long as $c(M-3)>\log4$ where $c=\liminf\mu/n$.
    It follows that, for any such $M$, we have
    \begin{align*}
        W_n=\big(1+o(1)\big)\frac{1}{4^n}\sum_{m=3}^MH_{n,m-1}e^{-\mu m}\,,
    \end{align*}
    and since this is a finite sum, we can directly apply the result from Theorem~\ref{thm:ht-bdd tree} and replace $H_{n,m-1}$ by their equivalent asymptotic behaviour:
    \begin{align*}
        W_n=\big(1+o(1)\big)\sum_{m=3}^M\frac{\tan^2\left(\frac{\pi}{m}\right)}{m}\frac{e^{-\mu m}}{\left(1+\tan^2\left(\frac{\pi}{m}\right)\right)^n}=\big(1+o(1)\big)\sum_{m=3}^M\frac{\tan^2\left(\frac{\pi}{m}\right)}{m}e^{-n\lambda_{\mu/n}(m)}\,,
    \end{align*}
    where $\lambda_{\mu/n}$ is defined in~\eqref{eq:lx}.
    Observing that $\lambda_{\mu/n}(m)$ is bounded away from $0$ in that scenario, we see that all terms are at least of exponential decay and so only the exponent with minimal value can be kept.
    Moreover, using that $\lambda_{\mu/n}$ is strictly convex, we see that the difference between two distinct exponents, except possibly for $\lfloor t_{\mu/n}\rfloor$ and $\lceil t_{\mu/n}\rceil$, is bounded away from $0$ as well.
    Thus only these two values can be kept in the sum and further, since it is possible that $t_{\mu/n}<3$, we need to restrict it to $m\geq3$.

    Assume now that $\mu\ll n$.
    We recall that we also assume that $\mu\gg n^{1/4}$.
    Under these assumptions, we can apply Proposition~\ref{prop:W as a sum} and obtain that
    \begin{align*}
        W_n=\big(1+o(1)\big)\frac{\mu}{2n}e^{-\lambda_n(t_{\mu/n})}\sum_{|m-t_{\mu/n}|<\epsilon_nt_{\mu/n}}\exp\left(-\left(\frac{3}{2}+o(1)\right)\left(\frac{\mu^4}{2\pi^2n}\right)^{1/3}\big(m-t_{\mu/n}\big)^2\right)\,.
    \end{align*}
    But now, since $\mu^4/n\gg1$, a geometric bound on the sum tells us that the first term asymptotically characterizes the whole sum, leading to
    \begin{align*}
        \sum_{m\geq\lceil t_{\mu/n}\rceil+1}\exp\left(-\left(\frac{3}{2}+o(1)\right)\left(\frac{\mu^4}{2\pi^2n}\right)^{1/3}\big(m-t_{\mu/n}\big)^2\right)=o\left(\exp\left(-\left(\frac{3}{2}+o(1)\right)\left(\frac{\mu^4}{2\pi^2n}\right)^{1/3}\big(\lceil t_{\mu/n}\rceil-t_{\mu/n}\big)^2\right)\right)
    \end{align*}
    as well as
    \begin{align*}
        \sum_{m\leq\lfloor t_{\mu/n}\rfloor-1}^\infty\exp\left(-\left(\frac{3}{2}+o(1)\right)\left(\frac{\mu^4}{2\pi^2n}\right)^{1/3}\big(m-t_{\mu/n}\big)^2\right)=o\left(\exp\left(-\left(\frac{3}{2}+o(1)\right)\left(\frac{\mu^4}{2\pi^2n}\right)^{1/3}\big(\lfloor t_{\mu/n}\rfloor-t_{\mu/n}\big)^2\right)\right)\,.
    \end{align*}
    It follows that all terms except for $m\in\{\lfloor t_{\mu/n}\rfloor,\lceil t_{\mu/n}\rceil\}$ are negligible, thus proving the desired result.
\end{proof}

\paragraph{Acknowledgements.}BC and NM were supported by travel grants from the MSCA-RISE-2020 project RandNET, Contract 101007705, for a visit to McGill University. NM further acknowledges support from an AMS Simons travel grant. BC and NM gratefully acknowledge the hospitality of McGill University for these visits.

This work was partially accomplished when both BC and NM were affiliated with TU Eindhoven, where BC has received funding from the European Union's Horizon 2020 Research and Innovation Programme under the Marie Sk\l{}odowska-Curie Grant Agreement No.~101034253.
BC, NM and MÜ also acknowledge their respective current affiliations, the Department of Mathematics at the National University of Singapore, the Department of Mathematics at the University of Illinois, Urbana-Champaign and Laboratoire d'Informatique de Paris Nord (LIPN) at Université Sorbonne Paris-Nord.

LAB Acknowledges the support of NSERC and of the Canada Research Chairs program.

\bibliographystyle{abbrvnat}
\bibliography{exp_trees}

\end{document}